\documentclass[reqno,12pt,a4paper]{amsart}
\usepackage{amscd,amsfonts,amssymb,amsthm,latexsym}
\usepackage{srcltx}

\theoremstyle{plain}
\newtheorem{theorem}{Theorem}[section]
\newtheorem*{theorem*}{Theorem}
\newtheorem{corollary}{Corollary}[section]
\newtheorem{lemma}{Lemma}[section]
\newtheorem{proposition}{Proposition}[section]
\newtheorem{hypothesis}{Hypothesis}

\theoremstyle{definition}

\newtheorem*{definition*}{Definition}

\theoremstyle{remark}

\newtheorem*{remark*}{Remark}

\numberwithin{equation}{section}

\begin{document}
\raggedbottom 

\title[Consecutive primes in short intervals]{Consecutive primes in short intervals}

\author{Artyom Radomskii}

\begin{abstract} We obtain a lower bound for
\[
  \#\{x/2< p_{n}\leq x:\ p_n \equiv\ldots\equiv p_{n+m}\equiv a\text{ (mod $q$)},\ p_{n+m} - p_{n}\leq y\},
  \]where $p_{n}$ is the $n^{\text{\emph{th}}}$ prime.
\end{abstract}

 \address{Steklov Mathematical Institute of Russian Academy of Sciences\\
 8 Gubkina St., Moscow 119991, Russia}

\thanks{This work was performed at the Steklov International Mathematical Center
and supported by the Ministry of Science and Higher Education of the
Russian Federation (agreement no. 075-15-2019-1614).}

\email{artyom.radomskii@mi-ras.ru}

\maketitle

\section{Introduction}

Let $p_{n}$ denote the $n^{\text{\emph{th}}}$ prime. We prove the following result.

\begin{theorem}\label{T6}
 There are positive absolute constants $c$ and $C$ such that the following holds. Let $\varepsilon$ be a real number with $0<\varepsilon <1$. Then there is a number $c_{0}(\varepsilon)>0$, depending only on $\varepsilon$, such that if $x\in \mathbb{R}$, $y\in \mathbb{R}$, $m\in\mathbb{Z}$, $q\in\mathbb{Z}$, $a\in\mathbb{Z}$ are such that
\begin{gather*}
c_{0}(\varepsilon) \leq y\leq \ln x,\\
1\leq m \leq c\cdot\varepsilon \ln y,\quad 1\leq q \leq y^{1-\varepsilon},\quad (a,q)=1,
\end{gather*}
 then
\begin{align*}
&\#\{x/2< p_n\leq x:\ p_n\equiv\dots\equiv p_{n+m}\equiv a\ (\text{\textup{mod }}q),\ p_{n+m}-p_{n}\leq y\}\geq\\
&\geq\pi(x)\left(\frac{y}{2q\ln x}\right)^{\text{\textup{exp}}(Cm)}.
\end{align*}
\end{theorem}

 Theorem \ref{T6} extends a result of Maynard \cite[Theorem 3.3]{Maynard} which showed the same result but with $y=\varepsilon \ln x$.

From Theorem \ref{T6} we obtain

\begin{corollary}\label{C5}
There is an absolute constant $C>0$ such that if $m$ is a positive integer, $x$ and $y$ are real numbers satisfying $\text{\textup{exp}}(Cm)\leq y \leq \ln x$, then
\begin{align*}
\#\{x/2< p_n\leq x:\ &p_{n+m}-p_{n}\leq y\} \geq  \pi(x)\left(\frac{y}{2\ln x}\right)^{\text{\textup{exp}}(Cm)}.
\end{align*}
\end{corollary}

 Let us introduce some notation. The symbol $b|a$ means that $b$ divides $a$. For fixed $a$ the sum $\sum_{b|a}$ and the product $\prod_{b|a}$ should be interpreted as being over all positive divisors of $a$.

We shall use the notation of I.\,M. Vinogradov. The symbol $A\ll B$ means that $|A|\leq cB$, where $c$ is a positive absolute constant.

We reserve the letter $p$ for primes. In particular, the sum $\sum_{p\leq K}$ should be interpreted as being over all prime numbers not exceeding $K$.

We also use the following notation:

$\#A$ --- the number of elements of a finite set $A$.

$\mathbb{N}$ --- the set of all positive integers.

$\mathbb{Z}$ --- the set of all integers.

$\mathbb{R}$ --- the set of all real numbers.

$\mathbb{C}$ --- the set of all complex numbers.

$\mathbb{P}$ --- the set of all prime numbers.

$[x]$ --- the integer part of a number $x$, i.\,e. $[x]$ is the largest integer $n$ such that $n\leq x$.

$\{x\}$ --- the fractional part of a number $x$, i.\,e. $\{x\}=x-[x]$.

$\lceil x\rceil$ --- the smallest integer $n$ such that $n\geq x$.

$\textup{Re}(s)$ --- the real part of a complex number $s$.

$\textup{Im}(s)$ --- the imaginary part of a complex number $s$.

$(a_1,\ldots, a_n)$ --- the greatest common divisor of integers $a_1,\ldots, a_n$.

$[a_1,\ldots, a_n]$ --- the least common multiple of integers $a_1,\ldots, a_n$.

$\varphi(n)$ --- the Euler totient function:
\[
\varphi(n)= \#\{1\leq m \leq n:\ (m,n)=1\}.
\]

$\mu(n)$ --- the Moebius function, which is defined as follows:

i) $\mu (1)=1$;

ii) $\mu(n)=0$, if there is a prime $p$ such that $p^2| n$;

iii) $\mu (n)=(-1)^{s}$, if $n=q_1\cdots q_s$, where $q_1<\ldots < q_s$ are primes.

$\Lambda (n)$ --- the von Mangoldt function:
\[
\Lambda(n)=
\begin{cases}
\ln p, &\text{if $n=p^k$;}\\
0, &\text{if $n\neq p^k$.}
\end{cases}
\]

$P^{-}(n)$ --- the least prime factor of $n$ (by convention $P^{-}(1)=+\infty$).

$\binom{n}{k} = n!/(k! (n-k)!)$ --- the binomial coefficient.

If $(a,b)=1$, $a$ and $b$ are said to be prime to one another or \emph{coprime}. The numbers $a, b, c, \ldots, k$ are said to be coprime if every two of them are coprime.

 For real numbers $x$, $y$ we also use $(x,y)$  to denote the open interval and $[x,y]$ to denote the closed interval. Also by $(a_1,\ldots, a_n)$ we denote a vector. The usage of the notation should be clear from the context.

We put
\[
\sum_{\varnothing} = 0,\qquad \prod_{\varnothing}=1.
\]

We define
\[
\mathcal{M}= \{n\in \mathbb{N}:\ \mu(n)\neq 0\}.
\]

We use the following functions:

\begin{gather*}
li (x)=\int_{2}^{x}\frac{dt}{\ln t};\\
\Phi(x,z)=\#\{1\leq n\leq x:\ P^{-}(n)>z\};\\
\pi(x)=\sum_{p\leq x} 1;\\
\theta(x)=\sum_{p\leq x} \ln p;\\
 \psi(x)=\sum_{n\leq x} \Lambda(n);\\
\pi(x; q, a)=\sum_{\substack{p\leq x\\
p\equiv a\ (\textup{mod }q)}} 1;\\
\psi(x; q, a)=\sum_{\substack{n\leq x\\
n\equiv a\ (\textup{mod }q)}} \Lambda (n).
\end{gather*}

Let $m$ be an integer with $m>1$ and $a$ be an integer. If $(a,m)=1$, then $a^{\varphi(m)}\equiv 1$ (mod $m$) (the Fermat--Euler theorem). Let $d$ be the smallest positive value of $\gamma$ for which $a^{\gamma}\equiv 1$ (mod $m$). We call $d$ the \emph{order} of $a$ (mod $m$), and say that \emph{$a$ belongs to $d$} (mod $m$).

 Let $q$ be a positive integer. We recall that \emph{a Dirichlet character modulo $q$} is a function $\chi:\ \mathbb{Z}\to \mathbb{C}$ such that\\
 1) $\chi(n+q)=\chi(n)$ for all $n\in \mathbb{Z}$ (i.\,e. $\chi$ is a periodic function with the period $q$);\\
 2) $\chi(mn)=\chi(m)\chi(n)$ for all $m, n\in \mathbb{Z}$ (i.\,e. $\chi$ is a totally multiplicative function);\\
 3) $\chi(1) =1$;\\
 4) $\chi(n)= 0$ for all $n\in \mathbb{Z}$ such that $(n,q)>1$.

 By $X_q$ we denote the set of all Dirichlet characters modulo $q$. We recall that $\# X_q = \varphi(q)$ and that the \emph{principal character modulo $q$} is
  \[
  \chi_0 (n):=
  \begin{cases}
  1, &\text{if $(n,q)=1$;}\\
  0, &\text{if $(n,q)>1$.}
  \end{cases}
  \]Let $\chi\in X_q$. We say that the character \emph{$\chi$ restricted by $(n,q)=1$ has period $q_{1}$} if it has the property that $\chi (m)=\chi (n)$ for all $m, n \in \mathbb{Z}$ such that $(m,q)=1$, $(n,q)=1$ and $m\equiv n$ (mod $q_1$). Let $c(\chi)$ denote the \emph{conductor} of $\chi$, defined to be the least positive integer $q_{1}$ such that $\chi$ restricted by $(n,q)=1$ has period $q_{1}$. We say that $\chi$ is \emph{primitive} if $c(\chi)=q$, and \emph{imprimitive} precisely if $c(\chi)<q$. By $X^{*}_{q}$ we denote the set of all primitive characters modulo $q$. We observe that the principal character modulo 1 is primitive. On the other hand, any principal character modulo $q>1$ is imprimitive, since its conductor is clearly 1. For $\chi\in X_q$ we put
  \begin{gather*}
  E_{\chi_0}(\chi):=
  \begin{cases}
  1, &\text{if $\chi$ is the principal character modulo $q$;}\\
  0, &\text{otherwise;}
  \end{cases}\\
  \psi (x,\chi)=\sum_{n\leq x} \Lambda(n)\chi(n);\\
  \psi^{\prime}(x,\chi)=\psi(x,\chi)-E_{\chi_0}(\chi)x.
  \end{gather*}

   A character $\chi$ is called \emph{real} if $\chi(n)\in \mathbb{R}$ for all $n\in \mathbb{Z}$. A character $\chi$ is called \emph{complex} if there is an integer $n$ such that $\textup{Im}(\chi(n))\neq 0$.

  We say that characters $\chi_1$ and $\chi_2$ (modulo $q_1$ and modulo $q_2$ respectively) are equal and write $\chi_1 = \chi_2$ if $\chi_1 (n)=\chi_2 (n)$  for any integer $n$. Otherwise, we say that characters $\chi_1$ and $\chi_2$ are not equal and write $\chi_1 \neq \chi_2$.

  Let $\chi$ be a Dirichlet character modulo $q$. The corresponding \emph{$L$-function} is defined by series
    \[
  L(s,\chi)=\sum_{n=1}^{\infty}\frac{\chi(n)}{n^{s}}
  \]for $s\in \mathbb{C}$ such that $\textup{Re}(s)>1$. It is well-known, if $\chi$ is not the principal character modulo $q$, then $L(s,\chi)$ can be analytically continued to $\mathbb{C}$. If $\chi$ is the principal character modulo $q$, then $L(s,\chi)$ can be analytically continued to $\mathbb{C}\setminus \{1\}$ with a simple pole at $s=1$.

  We say that two linear functions with integer coefficients $L_1(n)=a_{1}n+b_{1}$ and $L_2(n)=a_{2}n+b_{2}$ are equal and write $L_1 = L_2$ if $a_1=a_2$ and $b_{1}=b_{2}$. Otherwise, we say that linear functions $L_1$ and $L_2$ are not equal and write $L_1 \neq L_2$.

 Let $\mathcal{L}=\{L_1,\ldots,L_k\}$ be a set of $k$ linear functions with integer coefficients:
\[
L_{i}(n)=a_i n+b_i,\ i=1,\ldots, k.
\]For $L(n)=a n+b,$ $a,b\in \mathbb{Z}$, we define
\[
\Delta_{L}=|a|\prod_{i=1}^{k}|a b_i - b a_i|.
\]We say that $L=an+b\in \mathcal{L}$ if there is $1 \leq i \leq k$ such that $L=L_i$. Otherwise, we say that $L=an+b\notin \mathcal{L}$.

This paper is organized as follows. In sections \ref{S_Lemmas_General}\,--\,\ref{S_Lemmas_Psi} we give necessary lemmas. In section  \ref{S_T6_C5} we prove Theorem \ref{T6} and Corollary \ref{C5}.

\section{Preparatory Lemmas}\label{S_Lemmas_General}
In this section we give some well-known lemmas which will be used in the following sections.

\begin{lemma}[The Fundamental Theorem of Arithmetic; see, for example, \mbox{\cite[Chapter 1]{Vinogradov}}]\label{L.Fundam.Arithm}
Let $n$ be an integer with $n>1$. Then $n$ is a product of primes and, apart from rearrangement of factors, $n$ can be expressed as a product of primes in one way only.
\end{lemma} Let $n$ be an integer with $n>1$. From Lemma \ref{L.Fundam.Arithm} we have
\[
n=\widetilde{q}_{1}\cdots \widetilde{q}_{l},
\]where $\widetilde{q}_{i}$, $i=1,\ldots, l$, are primes. The primes $\widetilde{q}_{i}$, $i=1,\ldots, l$, are not necessarily distinct, nor arranged in any particular order. If we arrange them in increasing order, associate sets of equal primes into single factors, we obtain
\[
n= q_{1}^{\alpha_{1}}\cdots q_{r}^{\alpha_{r}},
\]where $q_{1}<\ldots < q_{r}$ are primes and $\alpha_{1},\ldots, \alpha_{r}$ are positive integers. We then say that $n$ is expressed in \emph{standard form}. From Lemma \ref{L.Fundam.Arithm} it follows that the standard form of $n$ is unique.

\begin{lemma}[see, for example,  \mbox{\cite[Chapter 1]{Prachar}}]\label{L_prod(1-1/p)}
Let $x$ be a real number with $x\geq 2$. Then
\begin{gather*}
b_1\ln x \leq \prod_{p\leq x} \Bigl(1-\frac{1}{p}\Bigr)^{-1}\leq b_2 \ln x,\\
b_3 \ln x\leq \prod_{p\leq x} \Bigl(1+\frac{1}{p}\Bigr)\leq b_{4} \ln x,
\end{gather*}where $b_1$, $b_2$, $b_3$ and $b_{4}$ are positive absolute constants.
\end{lemma}

\begin{lemma}[see, for example, \mbox{\cite[Chapters 1, 2]{Ingham}}]\label{L_prime_theorem}
The limits $\lim_{x\to \infty}\psi(x)/x$, $\lim_{x\to \infty}\theta(x)/x$, $\lim_{x\to \infty}\pi(x)/(x/\ln x)$,
$\lim_{n\to \infty} p_{n}/(n\ln n)$ exist and
 \begin{gather*}
\lim_{x\to +\infty}\frac{\psi(x)}{x}= \lim_{x\to +\infty}\frac{\theta(x)}{x}=
\lim_{x\to +\infty}\frac{\pi(x)}{x/\ln x}=1,\\
\lim_{n\to +\infty} \frac{p_n}{n\ln n} = 1.
\end{gather*}
\end{lemma}From Lemma \ref{L_prime_theorem} we obtain

\begin{lemma}\label{L:pi_ineq}
The following holds
\begin{gather}
b_5 x \leq \psi (x) \leq b_6 x,\quad x\geq 2; \label{Notation:psi_est}\\
b_7 x \leq \theta(x) \leq b_8 x,\quad x\geq 2; \notag\\
b_9\frac{ x}{\ln x} \leq \pi(x) \leq b_{10}\frac{ x}{\ln x},\quad x\geq 2;\label{Notation:Pi_func_est}\\
b_{11} n\ln (n+2) \leq p_n \leq b_{12}n\ln (n+2),\quad n\geq 1,\notag
\end{gather}where $b_5$, $b_6$, $b_7$, $b_8$, $b_9$, $b_{10}$, $b_{11}$ and $b_{12}$ are positive absolute constants.
\end{lemma}

\begin{lemma}[see, for example, \mbox{\cite[Chapter 2]{Vinogradov}}]\label{L_about_Euler_func}
Let $n$ be an integer with $n>1$. Then
\[
\varphi(n)=n\prod_{p|n}\Bigl(1-\frac{1}{p}\Bigr).
\]
\end{lemma}
From Lemma \ref{L_about_Euler_func} we obtain the following two lemmas.
\begin{lemma}\label{L:Euler_func_ineq}
Let $m$ and $n$ be integers with $m\geq 1$ and $n\geq 1$. Then
\[
 \varphi(mn) \geq \varphi(m)\varphi(n).
\]
\end{lemma}

\begin{lemma}\label{L:Euler_func_multipl}
Let $m$ and $n$ be integers with $m\geq 1$, $n\geq 1$ and $(m,n)=1$. Then
\[
 \varphi(mn) = \varphi(m)\varphi(n).
\]
\end{lemma}

\begin{lemma}\label{L:Euler_func_representation}
Let $n$ be an integer with $n\geq 1$. Then
\begin{equation}\label{Euler_Tozhd}
\frac{n}{\varphi(n)} = \sum_{d|n} \frac{\mu^{2}(d)}{\varphi(d)}.
\end{equation}
\end{lemma}
\textsc{Proof of Lemma \ref{L:Euler_func_representation}.} For $n=1$ the equality \eqref{Euler_Tozhd} holds. Let $n>1$. Let $n$ be expressed in standard form
\[
n=q_{1}^{\alpha_{1}}\cdots q_{r}^{\alpha_{r}},
\]where $q_{1}<\ldots < q_{r}$ are prime numbers. Applying Lemmas \ref{L_about_Euler_func} and \ref{L:Euler_func_multipl}, we have
\begin{align*}
\frac{n}{\varphi(n)}&=\prod_{p|n}\Bigl(1-\frac{1}{p}\Bigr)^{-1}=
\prod_{p|n}\Bigl(1+\frac{1}{p-1}\Bigr)=\\
&=\Bigl(1+\frac{1}{q_{1}-1}\Bigr)\cdots \Bigl(1+\frac{1}{q_{r}-1}\Bigr)=\\
&=\Bigl(1+\frac{1}{\varphi(q_{1})}\Bigr)\cdots \Bigl(1+\frac{1}{\varphi(q_{r})}\Bigr)=\\
&=1+ \sum_{s=1}^{r}\sum_{1\leq i_{1}<\ldots< i_{s}\leq r} \frac{1}{\varphi(q_{i_{1}})\cdots \varphi(q_{i_{s}})}=\\
&=1+ \sum_{s=1}^{r}\sum_{1\leq i_{1}<\ldots< i_{s}\leq r} \frac{1}{\varphi(q_{i_{1}}\cdots q_{i_{s}})}=\\
&=\sum_{\substack{d|n\\ d\in \mathcal{M}}}\frac{1}{\varphi(d)}=
\sum_{d|n}\frac{\mu^{2}(d)}{\varphi(d)}.
\end{align*}Lemma \ref{L:Euler_func_representation} is proved.

\begin{lemma}[see, for example, \mbox{\cite[Chapter 1]{Prachar}}]\label{L:low_est_Euler}
Let $n$ be an integer with $n\geq 1$. Then
\[
\varphi(n)\geq c\,\frac{n}{\ln\ln (n+2)},
\]where $c>0$ is an absolute constant.
\end{lemma}

\begin{lemma}[see, for example, \mbox{\cite[Chapter 28]{Davenport}}]\label{L:series_1/Euler_func}
Let $x$ be a real number with $x \geq 2$. Then
\[
\sum_{1 \leq n \leq x}\frac{1}{\varphi (n)} \leq c\ln x,
\]where $c>0$ is an absolute constant.
\end{lemma}

\begin{lemma}[see, for example, \mbox{\cite[Chapter 5]{Halberstam_Richert}}]\label{Lemma_sum_ln_p}
Let $n$ be an integer with $n\geq 1$. Then
\[
\sum_{p|n}\frac{\ln p}{p}\leq c \ln\ln (3 n),
\]where $c>0$ is an absolute constant.
\end{lemma}

\begin{lemma}\label{L_Diophantine}
Let $a$, $b$ and $c$ be integers with $(a,b)|c$. Then the equation
\begin{equation}\label{Diophant_II}
a x+b y = c
\end{equation}has a solution in the integers.
\end{lemma}
\textsc{Proof of Lemma \ref{L_Diophantine}.} We put $d=(a,b)$. We have $c=dl$, where $l\in\mathbb{Z}$. It is well-known
(see, for example, \cite[Chapter 1, Exercise 1]{Vinogradov}), the equation
\begin{equation}\label{Diophant_I}
ax+by=d
\end{equation}has a solution in the integers. Let $x_{0}\in\mathbb{Z}$ and  $y_{0}\in\mathbb{Z}$ be a solution of \eqref{Diophant_I}. Then the integers  $lx_{0}$ and $ly_{0}$  satisfy \eqref{Diophant_II}. Lemma \ref{L_Diophantine} is proved.

\begin{lemma}\label{L_Binom_C_n_k}
Let $n$ and $k$ be integers with $1 \leq k \leq n$. Then
\begin{equation}\label{Binom_ineq_for_L22}
\binom{n}{k} \geq k^{-k} (n-k)^{k}.
\end{equation}
\end{lemma}
\textsc{Proof of Lemma \ref{L_Binom_C_n_k}.} For $k=n$ the inequality \eqref{Binom_ineq_for_L22} holds. Let $1\leq k < n$. Then
\begin{align*}
\binom{n}{k}&=\frac{n!}{k! (n-k)!}=\frac{n(n-1)\cdots(n-k+1)}{k!}\geq \frac{(n-k)^k}{k!}\geq\\
&\geq k^{-k}(n-k)^{k}.
\end{align*} Lemma \ref{L_Binom_C_n_k} is proved.

\begin{lemma}[see \mbox{\cite[Chapter 0]{Hall}}]\label{L_about_Phi(x,z)}
Let $x$ and $z$ be real numbers with $2 \leq z \leq x/2$. Then
\[
\Phi (x,z) \geq c_{0} \frac{x}{\ln z},
\]where $c_{0}>0$ is an absolute constant.
\end{lemma}

\section{Lemmas on Dirichlet characters}\label{S_Lemmas_Dir_Charact}

In this section we give some well-known lemmas on Dirichlet characters which will be used in the following sections.

\begin{lemma}\label{Lemma_A5}
Let $q$ be an integer with $q>1$ and $\chi$ be a Dirichlet character modulo $q$. Let $n$ be an integer with $(n,q)=1$. Then
\[
\chi(n)=\textup{exp}\Big(2\pi i \frac{b}{\varphi(q)}\Big),\quad b\in\{0,\ldots, \varphi(q)-1\}.
\]In particular, $|\chi (n)|=1$.
\end{lemma}
\textsc{Proof of Lemma \ref{Lemma_A5}.} By the Fermat--Euler theorem, we have
\[
n^{\varphi(q)} \equiv 1\ \  \text{(mod $q$)}.
\]Hence,
\[
(\chi(n))^{\varphi(q)}=1.
\]We obtain
\[
\chi(n)=\textup{exp}\Big(2\pi i \frac{b}{\varphi(q)}\Big),\quad b\in\{0,\ldots, \varphi(q)-1\}.
\]In particular, $|\chi (n)|=1$. Lemma \ref{Lemma_A5} is proved.

\begin{lemma}\label{Lemma_B5}
Let $q$ be a positive integer and $\chi$ be a real character modulo $q$. Let $n$ be an integer with $(n,q)=1$. Then $\chi(n)\in \{-1, 1\}$.
\end{lemma}
\textsc{Proof of Lemma \ref{Lemma_B5}.} If $q=1$, then $\chi(m)=1$ for any integer $m$. Let $q$ be an integer with $q>1$ and $\chi$ be a real character modulo $q$. Let $n$ be an integer with $(n,q)=1$. From Lemma \ref{Lemma_A5} we have $|\chi(n)|=1$. Since $\chi(n)$ is real, we see that $\chi(n)\in \{-1, 1\}$. Lemma \ref{Lemma_B5} is proved.

\begin{lemma}\label{Lemma_a_b}
Let $a$, $b$ and $n$ be integers with $1\leq a < b$, $a|b$ and $(n,a)=1$. Then there is an integer $t$ such that $(n+ta, b)=1$.
\end{lemma}
\textsc{Proof of Lemma \ref{Lemma_a_b}.} If $(n,b)=1$, we take $t=0$. Let $(n,b)>1$. Then the set
\[
\Omega = \{ p|b:\ p \nmid a\} \neq \varnothing.
\]Let
\[
\Omega = \{ q_1,\ldots, q_r\},\quad q_1<\ldots<q_r.
\]Let $1 \leq i \leq r$. Since $(a,q_i)=1$, the congruence
\[
n + t a\equiv 1\ (\textup{mod } q_i)
\]has a solution, i.\,e. there is an integer $m_i$ such that $n+ a m_i \equiv 1 \ (\textup{mod } q_i)$. We consider the system
\begin{equation}\label{L:systema}
\begin{cases}
t\equiv m_1\ (\textup{mod } q_1) \\
\ \ \ \vdots\\
t\equiv m_r\ (\textup{mod } q_r) .
\end{cases}
\end{equation}Since the numbers $q_1, \ldots, q_{r}$ are coprime, the system has a solution. Let an integer $t_0$ satisfy the system \eqref{L:systema}. We claim that the number $t_0$ is desired, i.\,e. we claim that $(n+t_0 a, b)=1$. Assume the converse: $(n+ t_0 a, b)>1$. Then there is a prime $p$ such that $p|b$ and $p| (n+t_0 a)$. If $p\nmid a$, then $p\in \Omega$, i.\,e. $p=q_i$ for some $1 \leq i \leq r$. But
\[
n+t_0 a \equiv 1\ (\textup{mod } q_i)
\]and hence $p\nmid (n+t_0a)$. We obtain a contradiction. Thus this case is impossible. Hence, $p|a$. Since $p| (n+t_0 a)$, we obtain that $p|n$. Hence, $(n,a)>1$. This contradicts the assumption of the lemma. Therefore, the assumption $(n+ t_0 a, b)>1$ is false. Hence, $(n+ t_0 a, b)=1$. Lemma \ref{Lemma_a_b} is proved.

\begin{lemma}\label{Lemma_about_(q,q1)}
 Let $q$ be an integer with $q \geq 2$, let $\chi\in X_q$. Suppose that $\chi$ restricted by $(n,q)=1$ has period $q_1$. Then $\chi$ restricted by $(n,q)=1$ has also period $(q,q_1)$.
\end{lemma}
\textsc{Proof of Lemma \ref{Lemma_about_(q,q1)}.} We put $\delta = (q, q_1)$. Let $m$ and $n$ be integers such that $(m,q)=1$, $(n,q)=1$ and $m\equiv n$ (mod $\delta$). We must prove that $\chi (m) = \chi (n)$. By Lemma \ref{L_Diophantine}, there are integers $k$ and $l$ such that
\[
m+ q_1 k= n + q l.
\]We put $A=m+ q_1 k= n+ql$. Since  $(n,q)=1$, we have $(n+ql, q)=1$. Hence, $(A,q) =1$. Since $\chi$ has period $q$, we have
\[
\chi(A)=\chi(n+ql)=\chi(n).
\]Since $(A,q)=1$, $(m,q)=1$ and $A\equiv m$ (mod $q_1$), we have $\chi(A) = \chi (m)$. Hence, $\chi (m)= \chi (n)$. Lemma \ref{Lemma_about_(q,q1)} is proved.

\begin{lemma}\label{L_c(chi)_delit_q}
Let $q\geq 1$ and $\chi\in X_q$. Then $c(\chi)$ divides $q$.
\end{lemma}
\textsc{Proof of Lemma \ref{L_c(chi)_delit_q}.} If $q=1$, then $c(\chi)=1$ and the statement is obvious. Let $q\geq 2$. By Lemma \ref{Lemma_about_(q,q1)}, $\chi$ restricted by $(n,q)=1$ has period $\delta = (c(\chi), q)$. If $c(\chi)$ is not a divisor of $q$, then $\delta< c(\chi)$, which contradicts the definition of the conductor. Lemma \ref{L_c(chi)_delit_q} is proved.

\begin{lemma}\label{L_induced}
Let $q\geq 1$ and $\chi\in X_q$. Then there exists a unique Dirichlet character $\chi_1\in X_{c(\chi)}$ such that
\begin{equation}\label{chi_induced_chi1_formula}
\chi(n)=\begin{cases}
\chi_{1}(n), &\text{if $(n, q)=1;$}\\
\ \ \,0, &\text{if $(n, q)>1$.}
\end{cases}
\end{equation}Furthermore, $\chi_1$ is primitive.
\end{lemma}
We say that $\chi_1$ \emph{induces} $\chi$.

\textsc{Proof of Lemma \ref{L_induced}.} I) Let $q=1$. Then $c(\chi)=1$, $\# X_1 = 1$, $\chi_1 = \chi$ and the statement is obvious.

II) Let $q \geq 2$ and $\chi$ be a primitive character modulo $q$. Then $c(\chi)=q$ and we can take $\chi_1=\chi$. Let us prove the uniqueness. Suppose that there are two different characters $\chi_1, \chi_2\in X_q$ satisfying \eqref{chi_induced_chi1_formula}. Then for any $n$ such that $(n,q)>1$, we have $\chi_1 (n)=0=\chi_2(n)$. For any $n$ such that $(n,q)=1$, we have $\chi_1 (n)=\chi(n)=\chi_2 (n)$. Hence, $\chi_1 (n)=\chi_2 (n)$ for any integer $n$, i.\,e. $\chi_1=\chi_2$. We obtain a contradiction.

III) Let $q\geq 2$ and $\chi$ be an imprimitive character modulo $q$. Then $1\leq c(\chi) < q$ and by Lemma \ref{L_c(chi)_delit_q} we have $c(\chi)| q$. We define $\chi_1$. Let $n\in \mathbb{Z}$. We consider several cases.\\
 1) $(n,c(\chi))>1$. Then we put $\chi_1 (n)=0$.\\
 2) $(n,c(\chi))=1$. By Lemma \ref{Lemma_a_b} there is an integer $t$ such that
 \begin{equation}\label{Lemma3.4.Super.t}
 (n+tc(\chi),q) = 1.
 \end{equation} We put $$\chi_1 (n) = \chi (n+tc(\chi)).$$ The choice of $t$, subject to \eqref{Lemma3.4.Super.t}, is immaterial, since $\chi$ restricted by $(n,q)=1$ has period $c(\chi)$. Thus, $\chi_1 (n)$ is defined for any integer $n$. We claim that $\chi_1$ is a character modulo $c(\chi)$. By construction,
 \[
 \chi_1 (n)=0 \text{ for any $n\in \mathbb{Z}$ such that $(n,c(\chi))>1$}.
 \]
  By Lemma \ref{Lemma_a_b}, there is an integer $t$ such that $(1+tc(\chi),q) = 1$. Since the choice of such $t$ is immaterial, we take $t=0$. We have
  \[
  \chi_1 (1)=\chi(1)=1.
  \]Now we prove that
 \begin{equation}\label{period_chi1}
 \chi_1 (n+c(\chi))= \chi_1 (n) \text{ for any $n\in \mathbb{Z}$}.
 \end{equation}If $(n, c(\chi))>1$, then we have $(n+c(\chi),c(\chi))>1$. Hence,
 \[
 \chi_1(n+c(\chi))=0=\chi_1 (n).
 \]Let $(n, c(\chi))=1$. Then we have $(n+c(\chi), c(\chi))=1$. By Lemma \ref{Lemma_a_b}, there are integers $t_1$ and $t_2$ such that $(n+t_{1}c(\chi),q) = 1$ and $(n+c(\chi)+t_{2}c(\chi),q) = 1$. By construction, we have
 \begin{gather*}
 \chi_1(n) = \chi (n+t_{1} c(\chi)),\\
 \chi_1(n+c(\chi)) = \chi (n+c(\chi)+t_{2} c(\chi)).
 \end{gather*}Since $\chi$ restricted by $(n,q)=1$ has period $c(\chi)$, we have $\chi (n+t_{1} c(\chi))= \chi (n+c(\chi)+t_{2} c(\chi))$. Hence, $\chi_1(n) = \chi_1(n+c(\chi))$ and \eqref{period_chi1} is proved. Now we prove that
 \begin{equation}\label{multiple_chi1}
 \chi_{1}(m n)=\chi_1 (m) \chi_1 (n)\text{ for any $m, n\in \mathbb{Z}$}.
 \end{equation}If $(m, c(\chi))>1$, then we have $(mn, c(\chi))>1$. Hence, $\chi_1 (mn)=0$, $\chi_1(m)=0$. Therefore \eqref{multiple_chi1} holds. Similarly, \eqref{multiple_chi1} holds if $(n, c(\chi))>1$. Let $(m, c(\chi))=1$, $(n, c(\chi))=1$. Then $(mn, c(\chi))=1$. By Lemma \ref{Lemma_a_b}, there are integers $t_1$, $t_2$ and $t_3$ such that $(m+t_{1}c(\chi),q) = 1$, $(n+t_{2}c(\chi),q) = 1$ and $(mn+t_{3}c(\chi),q) = 1$. We put $m_1 = m+t_{1}c(\chi)$, $n_1 = n+t_{2}c(\chi)$ and $u = mn+t_{3}c(\chi)$. By construction,
 \begin{gather*}
 \chi_1(mn)=\chi(u),\\
 \chi_1(m)=\chi(m_1),\\
 \chi_1 (n)=\chi(n_1).
 \end{gather*}Since $\chi$ is a totally multiplicative function, we have
 \[
 \chi_1(m) \chi_1 (n)=\chi(m_1)\chi(n_1)=\chi(m_{1} n_1).
 \]Since $(m_1, q)=1$, $(n_1, q)=1$, we have  $(m_1 n_1, q)=1$. It is clear that $m_1 n_1 \equiv u$ (mod $c(\chi)$). Since $\chi$ restricted by $(n,q)=1$ has period $c(\chi)$, we have $\chi (u)=\chi (m_1 n_1)$. Hence, $\chi_1 (mn)= \chi_1 (m)\chi_1 (n)$ and \eqref{multiple_chi1} is proved. Thus, we have proved that $\chi_1$ is a character modulo $c(\chi)$, i.\,e. $\chi_1 \in X_{c(\chi)}$.

  Now we prove that $\chi_1$ satisfies \eqref{chi_induced_chi1_formula}. It suffices to show that
  \begin{equation}\label{chi=chi_1}
\chi_1(n)=\chi(n), \text{ if $(n,q)=1$}.
\end{equation}
 Since $(n,q)=1$, we have $(n,c(\chi))=1$ (see Lemma \ref{L_c(chi)_delit_q}). By Lemma \ref{Lemma_a_b}, there is an integer $t$ such that $(n+t c(\chi),q) = 1$. By construction $\chi_1 (n) = \chi (n+t c(\chi))$. Since $(n+t c(\chi),q) = 1$, $(n,q)=1$ and $n+tc(\chi) \equiv n$ (mod $c(\chi)$), we have $\chi(n+t c(\chi))=\chi(n)$. Hence, $\chi_1(n)=\chi(n)$ and \eqref{chi=chi_1} is proved.

   Now we prove that $\chi_1$ is a primitive character. Suppose that a positive integer $q_{2}$ satisfies the property that $\chi_{1}$ restricted by $(n, c(\chi))=1$ has period $q_{2}$. Let $m$ and $n$ be integers such that $(m,q)=1$, $(n,q)=1$ and $m\equiv n$ (mod $q_2$). By Lemma \ref{L_c(chi)_delit_q}, we have $(m, c(\chi))=1$, $(n, c(\chi))=1$. We have (see \eqref{chi=chi_1})
\[
\chi(m)=\chi_1 (m)=\chi_1 (n)=\chi(n).
\]Hence, $\chi$ restricted by $(n, q)=1$ has period $q_2$. From the definition of a conductor it follows that $q_2 \geq c(\chi)$. Hence, $\chi_1$ is a primitive character.

 Now we prove the uniqueness. Suppose that there are two different characters $\chi_1, \chi_2 \in X_{c(\chi)}$ satisfying \eqref{chi_induced_chi1_formula}. If $(n, c(\chi))>1$, then $\chi_1 (n)=0=\chi_2 (n)$. Let $(n,c(\chi))=1$. By Lemma \ref{Lemma_a_b}, there is an integer $t$ such that $(n+t c(\chi),q) = 1$. Since $\chi_1$ and $\chi_2$ are periodic functions with period $c(\chi)$, we have
\[
\chi_1 (n) = \chi_1 (n+tc(\chi))=\chi(n+tc(\chi))=\chi_2 (n+tc(\chi))=\chi_2 (n).
\]Thus, $\chi_1 (n)=\chi_2 (n)$ for any $n\in \mathbb{Z}$ and hence $\chi_1=\chi_2$. We obtain a contradiction. The uniqueness is proved. Lemma \ref{L_induced} is proved.

\begin{lemma}\label{Lemma_char_razl}
Let $q$ be an integer with $q>1$ and let $q$ be expressed in standard form
\[
q=q_{1}^{\alpha_{1}}\cdots q_{r}^{\alpha_{r}},
\]where $q_{1}<\ldots < q_{r}$ are primes and $\alpha_{1},\ldots,\alpha_{r}$ are positive integers. Let $\chi$ be a Dirichlet character modulo $q$. Then there exist unique characters $\chi_{i}$ modulo $q_{i}^{\alpha_{i}}$, $i=1,\ldots, r,$ such that
\begin{equation}\label{L2.17_Predstavlenie_char}
\chi(n)=\chi_{1}(n)\cdots \chi_{r}(n)\quad \text{for all $n$}.
\end{equation}
Furthermore, if the character $\chi$ is real, then all characters $\chi_{i}$, $i=1,\ldots, r$, are real. If the character $\chi$ is primitive, then all characters $\chi_{i}$, $i=1,\ldots, r$, are primitive.
\end{lemma}
\textsc{Proof of Lemma \ref{Lemma_char_razl}.} For any $1\leq i\leq r$ we take $A_{i}$ such that
\begin{equation}\label{L2.17_Ai}
\begin{cases}
A_{i}\equiv 1\ \text{(mod $q_{i}^{\alpha_{i}}$)},\\
A_{i}\equiv 0\ \text{(mod $q_{j}^{\alpha_{j}}$)}\ \text{for any $1\leq j\leq r$, $j\neq i$}.
\end{cases}
\end{equation}Since the moduli of these congruences are coprime, the system has a solution (see, for example, \cite[Chapter 4]{Vinogradov}). Thus, integers $A_{1},\ldots, A_{r}$ are defined.

  Let $1\leq i \leq r$ and $n\in\mathbb{Z}$. We put
  \begin{equation}\label{L2.17_BASIC}
  \chi_{i}(n)=\chi(nA_{i}+\sum_{\substack{1\leq j\leq r\\ j\neq i}}A_{j}).
  \end{equation}We claim that $\chi_{i}$ is a Dirichlet character modulo $q_{i}^{\alpha_{i}}$.

1) It is clear that $\chi_{i}: \mathbb{Z}\to \mathbb{C}$.

2) Let us show that
\begin{equation}\label{L2.17_Periodic}
\chi_{i}(n+q_{i}^{\alpha_{i}})= \chi_{i}(n)\ \text{for all $n\in \mathbb{Z}$}.
\end{equation}Let $n\in \mathbb{Z}$. We put
\begin{gather*}
m_{1}=(n+q_{i}^{\alpha_{i}})A_{i}+\sum_{\substack{1\leq j \leq r\\ j\neq i}}A_{j},\qquad m_{2}=nA_{i}+\sum_{\substack{1\leq j \leq r\\ j\neq i}}A_{j}.
\end{gather*}We have
\begin{gather*}
\chi_{i}(n+q_{i}^{\alpha_{i}})=\chi(m_{1}),\\
\chi_{i}(n)=\chi(m_{2}),\\
m_{1}-m_{2}=q_{i}^{\alpha_{i}}A_{i}.
\end{gather*}The number  $A_{i}$ is divisible by $q_{j}^{\alpha_{j}}$ for any $1\leq j \leq r$, $j\neq i$. Hence, the number $m_{1}-m_{2}$ is divisible by $q$, i.\,e.
\[
m_{1}\equiv m_{2}\ \text{(mod $q$)}.
\]Therefore $\chi(m_{1})=\chi(m_{2})$ and, hence, $\chi_{i}(n+q_{i}^{\alpha_{i}})=\chi_{i}(n)$. Thus, \eqref{L2.17_Periodic} is proved.

3) Let us show that
\begin{equation}\label{L2.17_ZERO}
\chi_{i}(n)=0\ \text{for all $n\in \mathbb{Z}$ such that $(n, q_{i}^{\alpha_{i}})>1$.}
\end{equation}Let $n\in \mathbb{Z}$ be such that $(n, q_{i}^{\alpha_{i}})>1$. We put
\[
m=nA_{i}+\sum_{\substack{1\leq j \leq r\\ j\neq i}}A_{j}.
\]We have
\[
\chi_{i}(n)=\chi(m).
\]The number $q_{i}$ divides $n$ and divides $A_{j}$ for any $1 \leq j \leq r$, $j\neq i$. Hence, the number $q_{i}$ divides $m$. Therefore $(m,q)>1$. We obtain $\chi(m)=0$ and, hence, $\chi_{i}(n)=0$. Thus, \eqref{L2.17_ZERO} is proved.

4) Let us show that
\begin{equation}\label{L2.17_ONE}
\chi_{i}(1)=1.
\end{equation}We put
\[
m=\sum_{j=1}^{r}A_{j}.
\]Then
\[
\chi_{i}(1)=\chi(m).
\]We see from \eqref{L2.17_Ai} that $m \equiv 1$ (mod $q_{i}^{\alpha_{i}}$), $i=1, \ldots, r$. The number $m-1$ is divisible by $q_{i}^{\alpha_{i}}$ for any $i=1,\ldots, r$. Hence, the number $m-1$ is divisible by $q$, i.\,e.
\[
m\equiv 1\ \text{(mod $q$)}.
\]We obtain $\chi(m)=\chi(1)=1$ and, hence, $\chi_{i}(1)=1$. Thus, \eqref{L2.17_ONE} is proved.

5) Let us show that
\begin{equation}\label{L2.17_MULTIPLE}
\chi_{i}(mn)=\chi_{i}(m)\chi_{i}(n)\ \text{for all $m, n\in \mathbb{Z}$}.
\end{equation}Let $m,n\in \mathbb{Z}$.

i) If $(m,q_{i})>1$, then (see \eqref{L2.17_ZERO})
\[
\chi_{i}(m)=0,\qquad \chi_{i}(mn)=0.
\]Hence, \eqref{L2.17_MULTIPLE} holds in this case.

ii) Similarly, \eqref{L2.17_MULTIPLE} holds if $(n,q_{i})>1$.

iii) Let $(m,q_{i})=(n,q_{i})=1$. We put
\begin{gather*}
l_{1}=mA_{i}+\sum_{\substack{1\leq j \leq r\\j\neq i}}A_{j},\quad l_{2}=nA_{i}+\sum_{\substack{1\leq j \leq r\\j\neq i}}A_{j},\quad
l_{3}=mnA_{i}+\sum_{\substack{1\leq j \leq r\\j\neq i}}A_{j}.
\end{gather*}We have
\begin{gather*}
\chi_{i}(m)\chi_{i}(n)=\chi(l_{1})\chi(l_{2})=\chi(l_{1}l_{2}),\\
\chi_{i}(mn)=\chi(l_{3}).
\end{gather*}Since
\begin{gather*}
l_{3}-l_{1}l_{2}=mnA_{i}+\sum_{\substack{1\leq j \leq r\\j\neq i}}A_{j}-
mnA_{i}^{2}-\\
-mA_{i}\sum_{\substack{1\leq j \leq r\\j\neq i}}A_{j} - nA_{i}\sum_{\substack{1\leq j \leq r\\j\neq i}}A_{j}
- \Bigl(\sum_{\substack{1\leq j \leq r\\j\neq i}}A_{j}\Bigr)^{2},
\end{gather*}we have (see \eqref{L2.17_Ai})
\[
l_{3}-l_{1}l_{2}\equiv 0\ \text{(mod $q_{j}^{\alpha_{j}}$)}\ \text{for any $1 \leq j \leq r$}.
\]Hence,  $l_{3} - l_{1}l_{2}$ is divisible by $q$, i.\,e.
\[
l_{3} \equiv l_{1}l_{2}\ \text{(mod $q$)}.
\]We obtain $\chi(l_{3})=\chi(l_{1}l_{2})$ and, hence, $\chi_{i}(mn)=\chi_{i}(m)\chi_{i}(n)$. Thus, \eqref{L2.17_MULTIPLE} is proved. We have proved that $\chi_{i}$ is a Dirichlet character modulo $q_{i}^{\alpha_{i}}$, $i=1,\ldots, r$.

Now we prove that \eqref{L2.17_Predstavlenie_char} holds. Let $n\in \mathbb{Z}$. We put
\[
n_{i}=nA_{i}+\sum_{\substack{1 \leq j \leq r\\ j\neq i}}A_{j},\quad i=1,\ldots, r.
\]We have
\[
\chi_{1}(n)\cdots\chi_{r}(n)=\chi(n_{1})\cdots\chi(n_{r})=\chi(n_{1}\cdots n_{r}).
\]From \eqref{L2.17_Ai} we obtain
\[
n_{1}\cdots n_{r} \equiv n\ \text{(mod $q_{s}^{\alpha_{s}}$)}\ \text{for any $1\leq s \leq r$}.
\]Hence,  $(n_{1}\cdots n_{r} - n)$ is divisible by $q$, i.\,e.
\[
n_{1}\cdots n_{r}\equiv n\ \text{(mod $q$)}.
\]Hence, $\chi(n_{1}\cdots n_{r})=\chi(n)$ and \eqref{L2.17_Predstavlenie_char} is proved.

Now we prove the uniqueness. Suppose that
\begin{equation}\label{L2.17_REPR2}
\chi(n)=\widetilde{\chi}_{1}(n)\cdots \widetilde{\chi}_{r}(n),
\end{equation}where $\widetilde{\chi}_{i}$ is a Dirichlet character modulo $q_{i}^{\alpha_{i}}$, $i= 1,\ldots, r$. Let $1 \leq i \leq r$ and $n\in \mathbb{Z}$. We have (see \eqref{L2.17_Ai})
\[
nA_{i}+\sum_{\substack{1\leq j\leq r\\ j\neq i}}A_{j}\equiv 1\ \text{(mod $q_{s}^{\alpha_{s}}$)}\ \text{for any $1\leq s \leq r$, $s\neq i$},
\]and
\[
nA_{i}+\sum_{\substack{1\leq j\leq r\\ j\neq i}}A_{j}\equiv n\ \text{(mod $q_{i}^{\alpha_{i}}$)}.
\]Hence,
\[
\widetilde{\chi}_{s}(nA_{i}+\sum_{\substack{1\leq j\leq r\\ j\neq i}}A_{j})=
1\ \text{for any $1 \leq s \leq r$, $s\neq i$ },
\]and
\[
\widetilde{\chi}_{i}(nA_{i}+\sum_{\substack{1\leq j\leq r\\ j\neq i}}A_{j})=
\widetilde{\chi}_{i}(n).
\]From \eqref{L2.17_REPR2} we obtain
\[
\chi(nA_{i}+\sum_{\substack{1\leq j\leq r\\ j\neq i}}A_{j})=
\widetilde{\chi}_{i}(n).
\]Hence (see \eqref{L2.17_BASIC}), $\widetilde{\chi}_{i}(n)= \chi_{i}(n)$. Since this equation holds for any $n\in \mathbb{Z}$, we have $\widetilde{\chi}_{i} = \chi_{i}$, $i=1,\ldots, r$. 

We see from \eqref{L2.17_BASIC} that if the character $\chi$ is real, then all characters $\chi_{i}$, $i=1,\ldots, r$, are real. We claim that if the character $\chi$ is primitive, then all characters $\chi_{i}$, $i=1,\ldots, r$, are primitive. Assume the converse: there is $1\leq i \leq r$ such that the character $\chi_{i}$ is imprimitive. Then $c(\chi_{i})<q_{i}^{\alpha_{i}}$. Since $c(\chi_{i})| q_{i}^{\alpha_{i}}$ (see Lemma \ref{L_c(chi)_delit_q}), we have
\[
c(\chi_{i})=q_{i}^{\beta},\quad \beta < \alpha_{i}.
\]We put
\[
\widetilde{q}=q_{i}^{\beta}\prod_{\substack{1\leq j \leq r\\ j\neq i}}q_{j}^{\alpha_{j}}.
\]Let us show that the character $\chi$ restricted by $(n,q)=1$ has period $\widetilde{q}$. Let $m$ and $n$ be integers such that $(m,q)=(n,q)=1$ and $m\equiv n$ (mod $\widetilde{q}$). Let $1\leq j \leq r$, $j\neq i$. Since
\[
m\equiv n\ \text{(mod $q_{j}^{\alpha_{j}}$)},
\]we have $\chi_{j}(m)=\chi_{j}(n)$. Since $(m,q_{i}^{\alpha_{i}})=(n, q_{i}^{\alpha_{i}})=1$,
\[
m \equiv n\ \text{(mod $q_{i}^{\beta}$)}
\]and $\chi_{i}$ restricted by $(n, q_{i}^{\alpha_{i}})=1$ has period $q_{i}^{\beta}$, we have
$\chi_{i}(m)=\chi_{i}(n)$. We obtain
\[
\chi(m)=\chi_{i}(m)\prod_{\substack{1\leq j \leq r\\ j \neq i}}\chi_{j}(m)=
\chi_{i}(n)\prod_{\substack{1\leq j \leq r\\ j \neq i}}\chi_{j}(n)=\chi(n).
\]We have proved that $\chi$ restricted by $(n,q)=1$ has period $\widetilde{q}$. But then $c(\chi)\leq \widetilde{q}<q$. This contradicts the fact that the character $\chi$ is primitive. Hence, all characters $\chi_{i}$, $i=1,\ldots, r$, are primitive. Lemma \ref{Lemma_char_razl} is proved.

\begin{lemma}\label{L2.18_PRIMITIVE}
   Let $q$ be a positive integer such that there exists a real primitive character $\chi$ modulo $q$. Then the number $q$ is of the form $2^{\alpha}k$, where $\alpha\in\{0,\ldots, 3\}$ and $k\geq 1$ is an odd square-free integer.
\end{lemma}
\textsc{Proof of Lemma \ref{L2.18_PRIMITIVE}.} Modulo $q=1$ there exists a real primitive character; namely, $\chi(n)=1$ for all $n\in \mathbb{Z}$. The number 1 is of the form $2^{\alpha}k$; namely, $\alpha = 0$ and $k=1$.

Let $q$ be an integer such that $q>1$ and there exists a real primitive character $\chi$ modulo $q$. Suppose that $q=p^{r}s$, where $p\geq 3$ is a prime number, $(p,s)=1$, $r\geq 2$. We put $\widetilde{q}=p^{r-1}s$. We claim that the character $\chi$ restricted by $(n,q)=1$ has period $\widetilde{q}$. Let $m$ and $n$ be integers such that $(m,q)=(n,q)=1$ and $m\equiv n\ \text{(mod $\widetilde{q})$}$. We have $m=n+ \widetilde{q}t$, $t\in \mathbb{Z}$, and
\begin{align}
m^{p^{r-1}}&=(n+ \widetilde{q}t)^{p^{r-1}}=n^{p^{r-1}}+\sum_{i=1}^{p^{r-1}}\binom{p^{r-1}}{i}(\widetilde{q}t)^{i}
n^{p^{r-1}-i}=\notag\\
&=n^{p^{r-1}}+\sum_{i=1}^{p^{r-1}}A_{i}t^{i}
n^{p^{r-1}-i},\ \text{where $A_{i} = \binom{p^{r-1}}{i}(\widetilde{q})^{i}$.}\label{L18_m_equiv_n_POWER}
\end{align}Let $2\leq i \leq p^{r-1}$. Then
\[
A_{i} = \binom{p^{r-1}}{i}(p^{r-1}s)^{i}=p^{r}s\binom{p^{r-1}}{i}p^{(i-1)r-i}s^{i-1}.
\]It is clear that $i-1\geq 1$. We claim that
\begin{equation}\label{L18_INEQ_i_r}
(i-1)r-i\geq 0
\end{equation}or, that is equivalent,
\[
i(r-1)\geq r.
\]In fact, since $i\geq 2$ and $r\geq 2$, we have
\[
i(r-1)\geq 2(r-1)\geq r.
\]The inequality \eqref{L18_INEQ_i_r} is proved. Hence, $A_{i}=p^{r}sN$, where $N\in \mathbb{N}$. Thus, for any $2 \leq i \leq p^{r-1}$ we have
\[
A_{i}\equiv 0\ \text{(mod $q$)}.
\]We have
\[
A_{1}=p^{r-1}(p^{r-1}s)=p^{r}s p^{r-2}.
\]Since $r\geq 2$, we obtain
\[
A_{1}\equiv 0\ \text{(mod $q$)}.
\]Hence (see \eqref{L18_m_equiv_n_POWER}),
\[
m^{p^{r-1}}\equiv n^{p^{r-1}}\ \text{(mod $q$)}.
\]Using properties of a character, we obtain
\[
(\chi(m))^{p^{r-1}}=(\chi(n))^{p^{r-1}}.
\]Since $(m,q)= (n,q) = 1$ and the character $\chi$ is real, by Lemma \ref{Lemma_B5} we have $\chi(m),\chi(n)\in\{-1,1\}$. Since $p\geq3$ is a prime number and $r\geq 2$ is an integer, we have $p^{r-1}$ is an odd positive integer. Therefore, if $\chi(m)=1$, then $\chi(n)=1$; if $\chi(m)=-1$, then $\chi(n)=-1$. Thus, $\chi(m)=\chi(n)$. We have proved that the character $\chi$ restricted by $(n,q)=1$ has period $\widetilde{q}$. We obtain
\[
c(\chi) \leq \widetilde{q} < q.
\]This contradicts the fact that $\chi$ is a primitive character. Hence, the number $q$ is of the form $2^{\alpha}k$, where $\alpha \geq 0$ is an integer and $k\geq 1$ is an odd square-free integer.

We claim that $\alpha\leq 3$. Assume the converse: $\alpha  \geq 4$. Let
\[
k=q_{1}\cdots q_{r},
\]where $q_{1}<\ldots < q_{r}$ are odd primes. By Lemma \ref{Lemma_char_razl}, we have
\begin{equation}\label{L18_chi_PROD_chi_r}
\chi(n)=\chi_{1}(n)\chi_{2}(n)\cdots\chi_{r+1}(n),
\end{equation}where $\chi_{1}$ is a real primitive character modulo $2^{\alpha}$, $\chi_{i}$ is a real primitive character modulo $q_{i-1}$, $i=2,\ldots, r+1$ (if $k=1$, then $\chi_{2},\ldots,\chi_{r+1}$ are omitted in \eqref{L18_chi_PROD_chi_r}). It is well-known (see, for example, \cite[Chapter 6]{Vinogradov}), if numbers $\nu$ and $\gamma$ run independently through sets $\{0,1\}$ and $\{0,\ldots, 2^{\alpha-2}-1\}$ respectively, then $(-1)^{\nu}5^{\gamma}$ runs through a reduced system of residues modulo $2^{\alpha}$. Hence, for any $n$ such that $(n,2)=1$ there are unique numbers $\nu(n)\in\{0,1\}$ and $\gamma(n)\in\{0,\ldots, 2^{\alpha-2}-1\}$ such that
\begin{equation}\label{L18_equa_mod_2}
n \equiv (-1)^{\nu(n)}5^{\gamma(n)}\ \text{(mod $2^{\alpha}$)}.
\end{equation}Since $(-1)^{2}=1$, we have $(\chi_{1}(-1))^{2}=1$. We obtain
\[
\chi_{1}(-1)= (-1)^{a},\ a\in\{0,1\}.
\]It is well-known (see, for example, \cite[Chapter 6]{Vinogradov}), the number $5$ belongs to $2^{\alpha-2}$ (mod $2^{\alpha}$) and, in particular,
\[
5^{2^{\alpha-2}}\equiv 1\ \text{(mod $2^{\alpha}$)}.
\]Hence,
\[
(\chi_{1}(5))^{2^{\alpha-2}}=1.
\]We obtain
\[
\chi_{1}(5) = \textup{exp}\Bigl(2\pi i\frac{b}{2^{\alpha-2}}\Bigr),\quad b\in\{0,\ldots, 2^{\alpha-2}-1\}.
\]We see from \eqref{L18_equa_mod_2} that if $n$ is such that $(n,2)=1$, then
\begin{equation}\label{L18_chi1}
\chi_{1}(n)=(-1)^{a\nu(n)}\textup{exp}\Bigl(2\pi i\frac{b\gamma(n)}{2^{\alpha-2}}\Bigr).
\end{equation} We claim that $(b,2)=1$. Assume the converse: $(b,2)>1$. We show that then $\chi_{1}$ restricted by $(n,2^{\alpha})=1$ has period $2^{\alpha-1}$. Let $m$ and $n$ be integers such that $(m,2^{\alpha})=(n,2^{\alpha})=1$ and $m\equiv n\ \text{(mod $2^{\alpha-1}$)}$. We have
\begin{gather*}
m \equiv (-1)^{\nu(m)}5^{\gamma(m)}\ \text{(mod $2^{\alpha}$)},\\
n \equiv (-1)^{\nu(n)}5^{\gamma(n)}\ \text{(mod $2^{\alpha}$)}.
\end{gather*}Since these congruences hold also modulo $2^{\alpha-1}$, we have
\begin{equation}\label{L18_index_equa}
(-1)^{\nu(m)}5^{\gamma(m)}\equiv
(-1)^{\nu(n)}5^{\gamma(n)}\ \text{(mod $2^{\alpha -1}$)}.
\end{equation}Since $\alpha\geq 4$, we obtain
\[
(-1)^{\nu(m)}5^{\gamma(m)}\equiv
(-1)^{\nu(n)}5^{\gamma(n)}\ \text{(mod $4$)}.
\]It is clear that
\begin{gather*}
(-1)^{\nu(m)}5^{\gamma(m)}\equiv (-1)^{\nu(m)}\ \text{(mod $4$)},\\
(-1)^{\nu(n)}5^{\gamma(n)}\equiv (-1)^{\nu(n)}\ \text{(mod $4$)}.
\end{gather*}Hence,
\[
(-1)^{\nu(m)}\equiv (-1)^{\nu(n)}\ \text{(mod $4$)}.
\]If $\nu(m)=0$, then $\nu(n)=0$; if $\nu(m)=1$, then $\nu(n)=1$. Thus,
\begin{equation}\label{L18_nu_m}
\nu(m)=\nu(n).
\end{equation}We obtain (see \eqref{L18_index_equa})
\[
5^{\gamma(m)}\equiv
5^{\gamma(n)}\ \text{(mod $2^{\alpha -1}$)}.
\]Suppose, for the sake of definiteness, that $\gamma(m) \geq \gamma(n)$. We have
\[
5^{\gamma(n)}(5^{\gamma(m)-\gamma(n)}-1)\equiv0\ \text{(mod $2^{\alpha -1}$)}.
\]Since $(5^{\gamma(n)}, 2^{\alpha-1})=1$, we obtain
\[
5^{\gamma(m)-\gamma(n)}-1\equiv0\ \text{(mod $2^{\alpha -1}$)}.
\]Hence,
\[
5^{\gamma(m)-\gamma(n)}\equiv 1\ \text{(mod $2^{\alpha -1}$)}.
\]Since $5$ belongs to $2^{\alpha-3}$ (mod $2^{\alpha-1}$), we have (see \cite[Chapter 6]{Vinogradov})
\[
\gamma(m)-\gamma(n)\equiv 0\ \text{(mod $2^{\alpha -3}$)}.
\]We obtain
\begin{equation}\label{L18_gammam}
\gamma(m)=\gamma(n)+ 2^{\alpha -3}t,
\end{equation}where $t\geq 0$ is an integer. Since $(b,2)>1$, we have
\begin{equation}\label{L18_b}
b=2\widetilde{b},
\end{equation} where $\widetilde{b}\geq 0$ is an integer. We obtain (see \eqref{L18_chi1}, \eqref{L18_nu_m}, \eqref{L18_gammam} and \eqref{L18_b})
\begin{align*}
&\chi_{1}(m)=(-1)^{a\nu(m)}\textup{exp}\Bigl(2\pi i\frac{\widetilde{b}\gamma(m)}{2^{\alpha-3}}\Bigr)=\\
&=(-1)^{a\nu(n)}\textup{exp}\Bigl(2\pi i\frac{\widetilde{b}(\gamma(n)+ 2^{\alpha -3}t)}{2^{\alpha-3}}\Bigr)=\\
&=(-1)^{a\nu(n)}\textup{exp}\Bigl(2\pi i\frac{\widetilde{b}\gamma(n)}{2^{\alpha-3}}\Bigr)
\textup{exp}(2\pi i \widetilde{b} t)=\\
&=(-1)^{a\nu(n)}\textup{exp}\Bigl(2\pi i\frac{\widetilde{b}\gamma(n)}{2^{\alpha-3}}\Bigr)=
\chi_{1}(n).
\end{align*}We have proved that $\chi_1$ restricted by $(n,2^{\alpha})=1$ has period $2^{\alpha-1}$. Hence,
\[
c(\chi_1)\leq 2^{\alpha-1}<2^{\alpha}.
\]This contradicts the fact that $\chi_1$ is a primitive character. Hence, $(b,2)=1$.

 For $n=5$ we have $\nu(5)=0$ and $\gamma(5)=1$. Hence (see \eqref{L18_chi1}),
 \[
 \chi_{1}(5)=\textup{exp}\Bigl(2\pi i\frac{b}{2^{\alpha-2}}\Bigr)=
 \textup{exp}\Bigl(\pi i\frac{b}{2^{\alpha-3}}\Bigr).
 \]Since $\alpha\geq 4$ and $(b,2)=1$, we have
 \[
 \textup{Im}(\chi_{1}(5))\neq 0.
 \] This contradicts the fact that $\chi_{1}$ is a real character. Hence, $0 \leq \alpha\leq 3$. Lemma \ref{L2.18_PRIMITIVE} is proved.

 \begin{lemma}\label{L_DIFFERENCE_PRIM_CHAR}
 Let $q_1$ and $q_2$ be positive integers with $q_{1} \neq q_{2}$, $\chi_{1}$ be a primitive character modulo $q_{1}$, $\chi_{2}$ be a primitive character modulo $q_{2}$. Then $\chi_{1} \neq \chi_{2}$.
 \end{lemma}
\textsc{Proof of Lemma \ref{L_DIFFERENCE_PRIM_CHAR}.} Assume the converse: $\chi_{1}=\chi_{2}$. Let $m$ and $n$ be integers such that $(m, q_{1})=(n, q_{1})=1$ and $m\equiv n$ (\textup{mod} $q_{2}$). Then
\[
\chi_{1}(m)=\chi_{2}(m)=\chi_{2}(n)=\chi_{1}(n).
\]Hence, $\chi_{1}$ restricted by $(n,q_{1})=1$ has period $q_{2}$. Hence, $c(\chi_{1}) \leq q_{2}$. Since $\chi_{1}$ is a primitive character modulo $q_{1}$, we have $c(\chi_{1})=q_{1}$. Hence, $q_{1}\leq q_{2}$. Similarly, it can be proved that $q_{2} \leq q_{1}$. Hence, $q_{1}=q_{2}$. We obtain a contradiction. Hence, $\chi_1 \neq \chi_2$. Lemma \ref{L_DIFFERENCE_PRIM_CHAR} is proved.

\section{Lemmas on $\psi(x, \chi)$}\label{S_Lemmas_Psi}

In this section we give some lemmas on $\psi(x,\chi)$. Most of these lemmas are well known.
 The proof of Lemma \ref{L4} is based on ideas of Maynard (see the proof of Theorem 3.2 in \cite{Maynard}). The proof of Lemma \ref{L5_Bombieri} follows a standard proof of the Bombieri--Vinogradov Theorem (see, for example, \cite[Chapter 28]{Davenport}).

\begin{lemma}\label{L_STAR_psi(u,Q,W)}
Let $Q$ and $W$ be integers with $Q\geq 2$ and $(W,Q) =1$.  Let $u$ be a real number with $u \geq 2$. Then
\[
\psi (u; Q, W) - \frac{u}{\varphi (Q)} = \frac{1}{\varphi (Q)} \sum_{\chi \in X_Q} \overline{\chi(W)} \psi^{\prime} (u,\chi).
\]Here the line denotes a complex conjugation.
\end{lemma}
\textsc{Proof of Lemma \ref{L_STAR_psi(u,Q,W)}.} We define
\[
I_{Q,W}(n)=\begin{cases}
1, &\text{if $n\equiv W$ (mod Q);}\\
0, &\text{otherwise}.
\end{cases}
\] Since (see, for example, \cite[Chapter 4]{Davenport})
\[
\frac{1}{\varphi (Q)}\sum_{\chi \in X_Q}\overline{\chi(W)}\chi(n)= I_{Q,W}(n),
\]we have
\begin{align*}
\psi(u; Q, W)&= \sum_{\substack{n\leq u\\
n\equiv W\,(\textup{mod }Q)}} \Lambda (n) = \sum_{n \leq u} \Lambda (n) I_{Q,W}(n)=\\
&=\sum_{n \leq u} \Lambda (n) \frac{1}{\varphi (Q)}\sum_{\chi \in X_Q}\overline{\chi(W)}\chi(n)=\\
&=\frac{1}{\varphi (Q)}\sum_{\chi \in X_Q}\overline{\chi(W)}\Bigl(\sum_{n\leq u}\Lambda(n) \chi(n)\Bigr)=\\
&= \frac{1}{\varphi (Q)}\sum_{\chi \in X_Q}\overline{\chi(W)}\psi (u,\chi).
\end{align*}Let $\chi_0$ be the principal character modulo $Q$. Since $(W,Q)=1$, we have $\chi_0 (W)=1$. We have
\[
\sum_{\chi \in X_Q} \overline{\chi(W)} E_{\chi_0}(\chi) u = \overline{\chi_{0} (W)} u = u.
\]Hence,
\begin{align*}
&\psi(u; Q, W) - \frac{u}{\varphi (Q)} = \frac{1}{\varphi (Q)}\sum_{\chi \in X_Q}
\overline{\chi(W)} \bigl(\psi (u,\chi) - E_{\chi_0}(\chi) u\bigr)=\\
&=\frac{1}{\varphi (Q)}\sum_{\chi \in X_Q}
\overline{\chi(W)} \psi^{\prime}(u,\chi).
\end{align*}Lemma \ref{L_STAR_psi(u,Q,W)} is proved.

\begin{lemma}[see, for example, \mbox{\cite[Chapter 14]{Davenport}}]\label{L1_Fact1}
 There is a positive absolute constant $a >0$ such that the following holds. If $\chi$ is a complex character modulo $q$, then $L(s,\chi)$ has no zero in the region defined by \textup{(}here $s=\sigma + i t$, $\sigma = \textup{Re}(s)$, $t=\textup{Im}(s)$\textup{)}
\[
\Omega_{q}:\qquad\quad\sigma >
\begin{cases}
1- a /\ln (q|t|), &\text{if $|t|\geq 1$,}\\
1- a /\ln q, &\text{if $|t|<1$.}
\end{cases}
\]If $\chi$ is a real nonprincipal character modulo $q$, the only possible zero of $L(s,\chi)$ in this region is a single \textup{(}simple\textup{)} real zero. Furthermore, for at most one of the real nonprincipal characters $\chi$ \textup{(}mod $q$\textup{)} can $L(s,\chi)$ have a zero in the region $\Omega_{q}$.
\end{lemma}

\textsc{Remark.} It is easy to see that the constant $a$ can be replaced by any constant $a^{*}$ such that $0< a^{*}< a$.
\begin{lemma}[see \mbox{\cite[Chapter 20]{Davenport}}]\label{L2_razlozh_psi}
Let $\chi$ be a nonprincipal character modulo $q$ and $2\leq T\leq u$. Then
\[
\psi (u,\chi)=-\frac{u^{\beta_1}}{\beta_1}+ R_{4}(u, T),
\]where
\[
|R_{4}(u, T)| \leq C\Bigl(u \ln^{2}(qu)\textup{exp}\bigl(-a\ln u/\ln (qT)\bigr)+uT^{-1}\ln^{2}(qu)+u^{1/4}\ln u\Bigr).
\]Here $C>0$ is an absolute constant, $a>0$ is the absolute constant in Lemma \textup{\ref{L1_Fact1}}. The term $-u^{\beta_1}/\beta_1$ is to be omitted unless $\chi$ is a real character for which $L(s, \chi)$ has a zero $\beta_1$ \textup{(}which is necessary unique, real and simple\textup{)} satisfying
\[
\beta_1 > 1- a/\ln q.
\]
\end{lemma}

\begin{lemma}[Theorem of Page (see, for example, \mbox{\cite[Chapter 14]{Davenport}})]\label{L3_Page}
There are absolute constants $a_1 >0$ and $a^{\prime}_{1}>0$ such that the following holds. Let $z$ be a real number with $z\geq 3$. Then there is at most one real primitive $\chi$ to a modulus $q_0$, $3\leq q_0 \leq z$, for which $L(s,\chi)$ has a real zero $\beta$ satisfying
\[
\beta> 1 - a_1/\ln z.
\]If such a character $\chi$ exists, then
\[
q_0 \geq \frac{a^{\prime}_{1} (\ln z)^{2}}{(\ln\ln z)^{4}}.
\]
\end{lemma} Such a modulus $q_0$ is said to be \emph{an exceptional modulus} in the interval $[3,z]$.

\begin{lemma}\label{L_q0_square_free}
 Let $z$ be a real number with $z\geq 3$. If an exceptional modulus $q_{0}$ in the interval $[3,z]$ exists, then the number $q_{0}$ is of the form $2^{\alpha}k$, where $\alpha\in\{0,\ldots, 3\}$ and $k\geq 1$ is an odd square-free integer.
\end{lemma}
\textsc{Proof of Lemma \ref{L_q0_square_free}.} Let an exceptional modulus $q_{0}$ in the interval $[3,z]$ exist. In particular, this means that there exists a real primitive character $\chi$ modulo $q_{0}$. By Lemma \ref{L2.18_PRIMITIVE}, the number $q_{0}$ is of the form $2^{\alpha}k$ where $\alpha\in\{0,\ldots, 3\}$ and $k\geq 1$ is an odd square-free integer. Lemma \ref{L_q0_square_free} is proved.

\begin{lemma}\label{L4}
There are positive absolute constants $c_0$, $c_1$, $\gamma_0$  and $C$ such that the following holds. Let $x$ be a real number with $x\geq c_0$, $q_0$ be an exceptional modulus in the interval $[3, \textup{exp}(2 c_1 \sqrt{\ln x})]$, $Q$ be an integer with $3\leq Q\leq \textup{exp}(2 c_1 \sqrt{\ln x})$ and  $Q\neq q_0$ \textup{(}the last inequality should be interpreted as follows: if $q_0$ exists, then $Q\neq q_0$; if $q_0$ does not exist, then $Q$ is any integer in the mentioned interval\textup{)}, $\chi$ be a primitive character modulo $Q$. Then
\[
\max_{2\leq u\leq x^{1+ \gamma_0/\sqrt{\ln x}}}|\psi(u, \chi)|\leq C x\,\textup{exp}(-3c_1\sqrt{\ln x}).
\]
\end{lemma}
\textsc{Proof of Lemma \ref{L4}.} We choose $c_1$ and $\gamma_0$ later. The number $c_0$ depends on $c_1$ and $\gamma_0$; the number $c_0$ is large enough and $x\geq c_{0}(c_1, \gamma_0)$. We put
\[
z=\textup{exp}(2 c_1 \sqrt{\ln x}).
\]We have
\[
z\geq 3,
\]if the number $c_0 (c_1, \gamma_0)$ is chosen large enough. By Lemma \ref{L3_Page}, there is at most one real primitive $\chi$ to a modulus $q_0$, $3\leq q_0 \leq z$, for which $L(s,\chi)$ has a real zero $\beta$ satisfying
  \begin{equation}\label{T6:Page_est}
 \beta > 1-a_1/\ln z= 1-a_1/ (2c_1 \sqrt{\ln x}).
 \end{equation} If such a character $\chi$ exists, then
 \begin{equation}\label{L4:q0_est}
 q_0 \geq \frac{a^{\prime}_{1}(\ln z)^{2}}{(\ln\ln z)^{4}}=
 \frac{a^{\prime}_{1}(2c_1\sqrt{\ln x})^{2}}{\bigl((1/2)\ln\ln x+ \ln (2c_1)\bigr)^{4}}\geq
 \frac{a^{\prime}_{1}c_{1}^{2}\ln x}{(\ln\ln x)^{4}},
 \end{equation}if $c_0 (c_1, \gamma_0)$ is chosen large enough. Let $Q$ be an integer with $3 \leq Q\leq \textup{exp}(2 c_1 \sqrt{\ln x})$ and $Q\neq q_0$, let $\chi$ be a primitive character modulo $Q$. Since $Q>1$, we see that $\chi$ is a nonprincipal character. By Lemma \ref{L2_razlozh_psi}, if $2 \leq T \leq u$, then
 \begin{equation}\label{T6:psi_formula}
\psi (u,\chi)=-\frac{u^{\beta_1}}{\beta_1}+ R_{4}(u, T),
\end{equation}where
\begin{align}
&|R_{4}(u, T)| \leq C\Bigl(u \ln^{2}(Qu)\textup{exp}\bigl(-a\ln u/\ln (QT)\bigr)+uT^{-1}\ln^{2}(Qu)+\notag\\
&+u^{1/4}\ln u\Bigr)= C(\Delta_1 + \Delta_2 +\Delta_3).\label{T6:R_4}
\end{align}The term $-u^{\beta_1}/\beta_1$ is to be omitted unless $\chi$ is a real character for which $L(s, \chi)$ has a zero $\beta_1$ \textup{(}which is necessary unique, real and simple\textup{)} satisfying
\[
\beta_1 > 1- a/\ln Q.
\] Let
\[
2 \leq u \leq x^{1+ \gamma_0/\sqrt{\ln x}}.
\]Let $u \geq c_{2}(c_1)$, where $c_2 (c_1)>0$ is a number depending only on $c_1$. We choose
\begin{equation}\label{T6:def_T}
T= \textup{exp}(4 c_1 \sqrt{\ln u}).
\end{equation}Then
\[
2 \leq T \leq u,
\]if $c_2(c_1)$ is chosen large enough.

I) Now we estimate
\[
\Delta_1 = u \ln^{2}(Qu)\textup{exp}\bigl(-a\ln u/\ln (QT)\bigr).
\] If $c_0 (c_1, \gamma_0)$ is chosen large enough, then
\begin{equation}\label{T6:<2}
1+\frac{\gamma_0}{\sqrt{\ln x}}\leq 2.
\end{equation}Hence,
\begin{gather}
\ln u \leq \Bigl(1+\frac{\gamma_0}{\sqrt{\ln x}}\Bigr)\ln x\leq 2 \ln x,\label{T6:ln<2}\\
QT \leq \textup{exp} (2 c_1 \sqrt{\ln x} +4 c_1 \sqrt{\ln u})\leq\textup{exp}(10 c_1 \sqrt{\ln x}),\notag\\
\ln (QT) \leq 10 c_1 \sqrt{\ln x},\notag\\
-\frac{a \ln u}{\ln (QT)}\leq -\frac{a\ln u}{10 c_1 \sqrt{\ln x}}.\notag
\end{gather}If $c_0 (c_1, \gamma_0)$ is chosen large enough, then
\[
\ln Q \leq 2 c_1 \sqrt{\ln x} \leq \ln x.
\]Hence,
\begin{equation}\label{T6:ln^2(Qu)}
\ln^{2} (Qu)\leq 2(\ln^2 Q + \ln^2 u)\leq 10 \ln^{2} x=10\,\textup{exp} (2\ln\ln x).
\end{equation}We have
\[
\Delta_1 \leq 10 u\,\textup{exp}\Bigl(-\frac{a\ln u}{10c_1\sqrt{\ln x}} + 2\ln\ln x\Bigr).
\]Let us consider two cases.

1) $x^{1/4} \leq u \leq x^{1+ \gamma_0 / \sqrt{\ln x}}$. Then
\[
\frac{\ln x}{4} \leq \ln u \leq  \Bigl(1+ \frac{\gamma_0}{\sqrt{\ln x}}\Bigr)\ln x\leq 2 \ln x.
\]If
\[
0< c_1 \leq \sqrt{\frac{a}{160}},
\]then
\[
-\frac{a}{40 c_1}\leq -4 c_1.
\]Hence,
\begin{gather*}
-\frac{a \ln u}{10 c_1 \sqrt{\ln x}} \leq -\frac{(a/4) \ln x}{10 c_1 \sqrt{\ln x}}=-\frac{a \sqrt{\ln x}}{40 c_1}
\leq - 4c_1 \sqrt{\ln x},\\
-\frac{a \ln u}{10 c_1 \sqrt{\ln x}} + 2\ln\ln x \leq -4 c_1 \sqrt{\ln x} + 2\ln\ln x\leq -3.5 c_1 \sqrt{\ln x},
\end{gather*}if $c_0 (c_1, \gamma_0)$ is chosen large enough. If
\[
0< \gamma_0 \leq 0.5 c_1,
\]then
\begin{align*}
&\Delta_1 \leq 10 x^{1+\gamma_0/\sqrt{\ln x}}\,\textup{exp}(-3.5 c_1 \sqrt{\ln x})=\\
&=10 x\,\textup{exp}
(-3.5 c_1 \sqrt{\ln x} + \gamma_0 \sqrt{\ln x})\leq 10 x\,\textup{exp}(-3 c_1\sqrt{\ln x}).
\end{align*}

2) $c_2 (c_1) \leq u < x^{1/4}$ (we may assume that $c_0 (c_1, \gamma_0)> \bigl(c_2 (c_1)\bigr)^{4}$ and $c_2 (c_1) \geq 10$). We have
\begin{align*}
\Delta_1 &\leq 10 u\,\textup{exp}\Bigl(-\frac{a\ln u}{10c_1\sqrt{\ln x}} + 2\ln\ln x\Bigr)\leq\\
&\leq 10 u\,\textup{exp}(2\ln\ln x)\leq 10 x^{1/4}\,\textup{exp}(2\ln\ln x)\leq\\
&\leq 10 x\,\textup{exp}(-3 c_1\sqrt{\ln x}),
\end{align*}if $c_0 (c_1, \gamma_0)$ is chosen large enough.

Thus, if $0< c_1 < \sqrt{a/160}$, $0 < \gamma_0 \leq 0.5 c_1$, $x \geq c_0 (c_1, \gamma_0)$, $c_{2}(c_1) \leq u \leq x^{1+\gamma_0/\sqrt{\ln x}}$, then
\[
\Delta_1 \leq 10 x\,\textup{exp}(-3 c_1\sqrt{\ln x}).
\]

II) Now we estimate
\[
\Delta_2 = u T^{-1} \ln^{2}(Qu).
\]From \eqref{T6:def_T} and \eqref{T6:ln^2(Qu)} we obtain
\[
\Delta_2 \leq 10 u\,\textup{exp}(-4 c_1\sqrt{\ln u} + 2\ln\ln x).
\]

Let us consider two cases.

1) $x^{0.9} \leq u \leq x^{1+ \gamma_0 / \sqrt{\ln x}}$. Then
\begin{gather*}
0.9 \ln x \leq \ln u \leq \Bigl(1+\frac{\gamma_0}{\sqrt{\ln x}}\Bigr)\ln x \leq 2 \ln x,\\
-4 c_1 \sqrt{\ln u}\leq -4 c_1 \sqrt{0.9 \ln x}< - 3.7 c_1 \sqrt{\ln x}.
\end{gather*}
Since $0< \gamma_0 \leq 0.5 c_1$, we have
\begin{align*}
\Delta_2&\leq 10 x^{1+ \gamma_0 / \sqrt{\ln x}}\,\textup{exp} (- 3.7 c_1 \sqrt{\ln x} + 2\ln\ln x)=\\
&= 10 x\,\textup{exp} (- 3.7 c_1 \sqrt{\ln x} + 2\ln\ln x+ \gamma_0 \sqrt{\ln x})\leq\\
&\leq 10 x\,\textup{exp} (- 3.2 c_1 \sqrt{\ln x} + 2\ln\ln x)\leq\\
&\leq 10 x\,\textup{exp} (- 3 c_1 \sqrt{\ln x}),
\end{align*}if $c_{0}(c_1, \gamma_0)$ is chosen large enough.

2) $c_2 (c_1) \leq u < x^{0.9}$. Then
\begin{align*}
\Delta_2 &\leq 10 u\,\textup{exp}(-4 c_1\sqrt{\ln u} + 2\ln\ln x)\leq\\
&\leq 10 u\,\textup{exp}( 2\ln\ln x)\leq 10 x^{0.9}\,\textup{exp}( 2\ln\ln x)\leq\\
&\leq 10 x\,\textup{exp}(-3 c_1 \sqrt{\ln x}),
\end{align*} if $c_0 (c_1, \gamma_0)$ is chosen large enough.

Thus, if $0< c_1 < \sqrt{a/160}$, $0 < \gamma_0 \leq 0.5 c_1$, $x \geq c_0 (c_1, \gamma_0)$, $c_{2}(c_1) \leq u \leq x^{1+\gamma_0/\sqrt{\ln x}}$, then
\[
\Delta_2 \leq 10 x\,\textup{exp}(-3 c_1\sqrt{\ln x}).
\]

III) Now we estimate
\[
\Delta_3 = u^{1/4}\ln u.
\]Since (see \eqref{T6:<2} and \eqref{T6:ln<2})
\begin{gather*}
\ln u\leq 2 \ln x,\\
u^{1/4}\leq x^{(1+\gamma_0/\sqrt{\ln x})/4}\leq x^{1/2},
\end{gather*}we have
\[
\Delta_3 \leq 2 x^{1/2}\ln x \leq x\,\textup{exp}(-3 c_1 \sqrt{\ln x}),
\]if $c_0 (c_1, \gamma_0)$ is chosen large enough.

Finally, we obtain (see \eqref{T6:R_4}): if $0< c_1 < \sqrt{a/160}$, $0 < \gamma_0 \leq 0.5 c_1$, $x \geq c_0 (c_1, \gamma_0)$, $c_{2}(c_1) \leq u \leq x^{1+\gamma_0/\sqrt{\ln x}}$, then
\begin{equation}\label{T6:R4_final_est}
|R_{4}(u, T)| \leq 21 C x\,\textup{exp}(-3 c_1\sqrt{\ln x}),
\end{equation}where $C>0$ is an absolute constant.

IV) Now we estimate (see \eqref{T6:psi_formula})
\[
\Delta_4 = \left|-\frac{u^{\beta_1}}{\beta_1}\right|.
\] If $\chi$ is not  a real character for which $L(s, \chi)$ has a zero $\beta_1$ \textup{(}which is necessary unique, real and simple\textup{)} satisfying
\[
\beta_1 > 1- a/\ln Q,
\] then the term $-u^{\beta_1}/\beta_1$ is to be omitted, and there is nothing to estimate. Let $\chi$ be such a character. Then $\chi$ is a real primitive character modulo $Q$. Since $Q\neq q_0$, we have (see Lemma \ref{L_DIFFERENCE_PRIM_CHAR} and \eqref{T6:Page_est})
\[
\beta_1 \leq 1-a_1/\ln z= 1-a_1/ (2c_1 \sqrt{\ln x}).
\]Hence,
\[
|u^{\beta_1}|= u^{\beta_1} \leq u^{1-a_1/ (2c_1 \sqrt{\ln x}) }=
u\,\textup{exp}\bigl(-(a_1 \ln u)/(2 c_1 \sqrt{\ln x})\bigr).
\]By Remark made below Lemma \ref{L1_Fact1}, we may assume that $0< a < 1/2$. Since $Q\geq 3$, we have
\[
\beta_1 > 1 - a/\ln Q> 1-\frac{1}{2\ln 3}> \frac{1}{2}.
\]Hence,
\[
0<\frac{1}{\beta_1}\leq 2.
\]We obtain
\begin{equation}\label{T6:Delta4_est}
\Delta_4 \leq 2 u\,\textup{exp}\bigl(-(a_1 \ln u)/(2 c_1 \sqrt{\ln x})\bigr).
\end{equation} Let us consider two cases.

1) $x^{1/2} \leq u \leq x^{1+\gamma_0/\sqrt{\ln x}}$. We have (see \eqref{T6:<2})
\[
\frac{\ln x}{2} \leq \ln u \leq \Bigl(1+\frac{\gamma_0}{\sqrt{\ln x}}\Bigr)\ln x \leq 2\ln x.
\]If we take
\[
0< c_1 < \sqrt{\frac{\min(a, a_1)}{160}},
\]then
\[
-\frac{a_1}{4c_1}\leq -3.5 c_1.
\]Hence,
\[
-\frac{a_1 \ln u}{2 c_1 \sqrt{\ln x}}\leq -\frac{(a_1/2)\ln x}{2 c_1\sqrt{\ln x}}= -\frac{a_1 \sqrt{\ln x}}{4 c_1}\leq
-3.5 c_1 \sqrt{\ln x}.
\]Since $0< \gamma_0 \leq 0.5 c_1$, we obtain (see \eqref{T6:Delta4_est})
\begin{align*}
\Delta_4&\leq 2 x^{1+\gamma_0/\sqrt{\ln x}}\,\textup{exp}(-3.5 c_1 \sqrt{\ln x})=\\
&=2x\,\textup{exp}(-3.5 c_1 \sqrt{\ln x}+\gamma_0\sqrt{\ln x})\leq\\
&\leq 2x\,\textup{exp}(-3 c_1 \sqrt{\ln x}).
\end{align*}

2) $c_2 (c_1) \leq u < x^{1/2}$. Then (see \eqref{T6:Delta4_est})
\[
\Delta_4 \leq 2 u \leq 2 x^{1/2}\leq 2 x\,\textup{exp}(-3 c_1 \sqrt{\ln x}),
\]if $c_0 (c_1, \gamma_0)$ is chosen large enough. Putting I--IV together, we obtain (see \eqref{T6:psi_formula}, \eqref{T6:R4_final_est}): if $0< c_1 < \sqrt{\min(a, a_1)/160}$, $0 < \gamma_0 \leq 0.5 c_1$, $x \geq c_0 (c_1, \gamma_0)$, $c_{2}(c_1) \leq u \leq x^{1+\gamma_0/\sqrt{\ln x}}$, then
\[
|\psi(u, \chi)| \leq (21 C +2) x\,\textup{exp}(-3 c_1\sqrt{\ln x}),
\]where $C>0$ is an absolute constant.

 There is a number $d(c_1)>0$, depending only on $c_1$, such that
\[
t\,\textup{exp}(-3 c_1 \sqrt{\ln t})\geq 1,\quad \text{if $t\geq d(c_1)$}.
\]We may assume that $c_0 (c_1, \gamma_0) > d(c_1)$. Hence, if $2 \leq u < c_2 (c_1)$, then (see \eqref{Notation:psi_est})
\begin{align*}
|\psi (u,\chi)|&=\Bigl|\sum_{n\leq u} \Lambda(n)\chi(n)\Bigr|\leq \sum_{n\leq u} \Lambda(n)=\psi(u)\leq\\
&\leq b_6 u\leq b_6 c_2 (c_1) \leq b_6 c_2 (c_1) x\,\textup{exp}(-3 c_1 \sqrt{\ln x}).
\end{align*}

Thus, if $0< c_1 < \sqrt{\min(a, a_1)/160}$, $0 < \gamma_0 \leq 0.5 c_1$, $x \geq c_0 (c_1, \gamma_0)$, then
\[
\max_{2\leq u \leq x^{1+\gamma_0 / \sqrt{\ln x}}}|\psi(u, \chi)| \leq \bigl(21 C +2+ b_6 c_2 (c_1)\bigr) x\,\textup{exp}(-3 c_1\sqrt{\ln x}),
\]where $C>0$ is an absolute constant. We take
\[
c_1=\frac{\sqrt{\min (a, a_1)}}{16},\quad \gamma_0 = 0.5 c_1= \frac{\sqrt{\min (a, a_1)}}{32}.
\] Since $a>0$ and $a_1 >0$ are absolute constants, we see that $c_1$, $\gamma_0$, $c_0 (c_1, \gamma_0)$ and $c_2 (c_1)$ are positive absolute constants. Lemma \ref{L4} is proved.

\begin{lemma}\label{L:chi_minus_chi1}
 Let $u$ be a real number with $u \geq 2$, $Q$ be an integer with $Q\geq 2$, $\chi \in X_Q$, $\chi_1$ be a primitive character modulo $q_1$ inducing $\chi$. Then
\begin{equation}\label{raznost_chi_minus_chi_1_est}
|\psi^{\prime}(u, \chi) - \psi^{\prime}(u, \chi_1)|\leq \ln^{2} (Qu).
\end{equation}
\end{lemma}
\textsc{Proof of Lemma \ref{L:chi_minus_chi1}.} From Lemma \ref{L_induced} and the definition of the inducing character, given below Lemma \ref{L_induced}, we have $q_1 = c(\chi)$ and, hence, $q_1| Q$ (see Lemma \ref{L_c(chi)_delit_q}). Let us prove \eqref{raznost_chi_minus_chi_1_est}. We consider two cases.

1) $\chi$ is a nonprincipal character modulo $Q$. We claim that then $\chi_1$ is a nonprincipal character modulo $q_{1}$. In fact, assume the converse. Then
\[
\chi_{1}(n)=\begin{cases}
1, &\text{if $(n, q_{1})=1;$}\\
0, &\text{if $(n, q_{1})>1$.}
\end{cases}
\]Since $q_1| Q$, we see that if $(n, Q) =1$, then $(n, q_1) =1$. We obtain (see \eqref{chi_induced_chi1_formula})
\[
\chi(n)=\begin{cases}
1, &\text{if $(n, Q)=1;$}\\
0, &\text{if $(n, Q)>1$.}
\end{cases}
\]Hence, $\chi$ is the principal character modulo $Q$. We obtain a contradiction. Thus, $\chi_1$ is a nonprincipal character modulo $q_1$. We have (see \eqref{chi_induced_chi1_formula})
\begin{align*}
&\psi^{\prime}(u, \chi)= \psi (u,\chi)=\sum_{n\leq u} \Lambda (n) \chi(n)=\\
&=\sum_{\substack{n\leq u\\
(n,Q)=1}} \Lambda (n) \chi(n)= \sum_{\substack{n\leq u\\
(n,Q)=1}} \Lambda (n) \chi_{1}(n);
\end{align*}
\begin{align*}
&\psi^{\prime}(u, \chi_{1})= \psi (u,\chi_{1})=\sum_{n\leq u} \Lambda (n) \chi_{1}(n)=\sum_{\substack{n\leq u\\
(n,q_{1})=1}} \Lambda (n) \chi_{1}(n)=\\
&=\sum_{\substack{n\leq u\\
(n,q_{1})=1\\ (n,Q)=1}} \Lambda (n) \chi_{1}(n)+ \sum_{\substack{n\leq u\\
(n,q_{1})=1\\ (n,Q)>1}} \Lambda (n) \chi_{1}(n)=\\
&= \sum_{\substack{n\leq u\\
 (n,Q)=1}} \Lambda (n) \chi_{1}(n)+ \sum_{\substack{n\leq u\\
(n,q_{1})=1\\ (n,Q)>1}} \Lambda (n) \chi_{1}(n).
\end{align*}Hence,
\begin{align*}
&\psi^{\prime}(u, \chi_{1}) - \psi^{\prime}(u, \chi) = \sum_{\substack{n\leq u\\
(n,q_{1})=1\\ (n,Q)>1}} \Lambda (n) \chi_{1}(n)=\sum_{\substack{p^{m}\leq u\\
p|Q\\
 (p,q_1)=1}}  \chi_{1}(p^{m})\ln p=\\
 &= \sum_{\substack{p|Q\\
 (p,q_1)=1}}\ln p \sum_{1 \leq m \leq \log_{p}u} \chi_{1}(p^{m}).
\end{align*}For any prime $p$ we have
\[
\Bigl|\sum_{1 \leq m \leq \log_{p}u} \chi_{1}(p^{m})\Bigr|\leq \log_{p}u = \frac{\ln u}{\ln p}\leq
\frac{\ln u}{\ln 2}\leq 2\ln u.
\]
 We obtain
\begin{align*}
&|\psi^{\prime}(u, \chi_{1}) - \psi^{\prime}(u, \chi)|\leq 2\ln u \sum_{\substack{p|Q\\
 (p,q_1)=1}}\ln p\leq 2\ln u\sum_{p|Q}\ln p\leq\\
 &\leq 2\ln u \ln Q\leq (\ln u + \ln Q)^{2}= \ln^{2} (Qu).
\end{align*}

2) $\chi$ is the principal character modulo $Q$. We have
\[
\psi^{\prime}(u, \chi) = \psi (u, \chi) - u;
\]\[
\chi(n)=\begin{cases}
1, &\text{if $(n, Q)=1$;}\\
0, &\text{if $(n, Q)>1$.}
\end{cases}
\]In this case the least period of $\chi$, restricted by $(n,Q)=1$, is $1$. Therefore $q_1 = 1$, and $\chi_1$ is a primitive character modulo $1$, i.\,e. $\chi_1 (n)=1$ for any integer $n$. Hence, $\chi_1$ is the principal character modulo $1$. We have
\begin{gather*}
\psi^{\prime}(u, \chi_1) = \psi(u, \chi_1) - u;\\
\psi^{\prime}(u, \chi_1) - \psi^{\prime}(u, \chi) = \psi(u, \chi_1) - \psi (u, \chi).
\end{gather*}Since
\begin{gather*}
\psi(u,\chi)=\sum_{n \leq u} \Lambda(n) \chi(n)=\sum_{\substack{n\leq u\\
(n, Q)=1}}\Lambda(n),\\
\psi(u,\chi_{1})=\sum_{n \leq u} \Lambda(n) \chi_{1}(n)=\sum_{n\leq u} \Lambda(n)=
\sum_{\substack{n\leq u\\
(n, Q)=1}}\Lambda(n) + \sum_{\substack{n\leq u\\
(n, Q)>1}}\Lambda(n),
\end{gather*}we have
\[
\psi^{\prime}(u, \chi_1) - \psi^{\prime}(u, \chi) = \sum_{\substack{n\leq u\\
(n, Q)>1}}\Lambda(n)=\sum_{\substack{p^{m} \leq u\\
p|Q}} \ln p=\sum_{p|Q}\ln p \sum_{1\leq m \leq \log_{p}u}1.
\]For any prime $p$ we have
\[
\sum_{1\leq m \leq \log_{p}u}1 \leq \log_{p}u = \frac{\ln u}{\ln p}\leq
\frac{\ln u}{\ln 2}\leq 2\ln u.
\]Hence,
\[
|\psi^{\prime}(u, \chi_1) - \psi^{\prime}(u, \chi)|\leq 2\ln u \sum_{p|Q}\ln p\leq 2\ln u\ln Q\leq
\ln^{2}(Qu).
\]Lemma \ref{L:chi_minus_chi1} is proved.

 \begin{lemma}[see \mbox{\cite[Chapter 28]{Davenport}}]\label{L:Dave_character}
 Let $Q_{1}$, $Q_{2}$ and $t$ be real numbers with $1\leq Q_1 < Q_2$, $t\geq 2$. Then
 \[
 \sum_{Q_1 < Q \leq Q_2}   \frac{1}{\varphi (Q)} \sum_{\chi\in X^{*}_{Q}} \max_{2\leq u \leq t}   |\psi (u,\chi)|
 \leq C\ln^{4}(t Q_{2})\Bigl(\frac{t}{Q_{1}}+ t^{5/6} \ln Q_{2}+ t^{1/2}Q_{2}\Bigr),
 \]where $C>0$ is an absolute constant.
 \end{lemma}

\begin{lemma}\label{L5_Bombieri}
 Let $\varepsilon$ and $\delta$ be real numbers with $0< \varepsilon < 1$ and $0< \delta < 1/2$. Then there is a number $c(\varepsilon, \delta)>0$, depending only on $\varepsilon$ and $\delta$, such that if $x\in \mathbb{R}$ and $q\in \mathbb{Z}$ are such that  $x\geq c(\varepsilon, \delta)$ and $1 \leq q \leq (\ln x)^{1-\varepsilon}$, then there is a positive integer $B$ such that
\[
1\leq B \leq \textup{exp} (c_1 \sqrt{\ln x}),\quad 1\leq \frac{B}{\varphi(B)}\leq 2,\quad (B,q)=1
\]and
\[
\sum_{\substack{1\leq Q\leq x^{1/2 - \delta}\\
(Q,B)=1}} \max_{2 \leq u \leq x^{1+ \gamma/\sqrt{\ln x}}}\max_{\substack{W\in \mathbb{Z}:\\ (W,Q)=1}}
\Bigl|\psi (u; Q, W)-\frac{u}{\varphi (Q)}\Bigr| \leq c_2 x\,\textup{exp}(-c_3\sqrt{\ln x}).
\]Here $c_1$, $\gamma$, $c_2$ and $c_3$ are positive absolute constants.
\end{lemma}
\textsc{Proof of Lemma \ref{L5_Bombieri}.} Let $c_0$, $c_1$, $\gamma_0$ and $C$ be the positive absolute constants in Lemma \ref{L4}. We choose $\gamma$ and $c(\varepsilon, \delta)=c(\varepsilon, \delta, \gamma)$ later; $\gamma$ is small enough, and $c(\varepsilon, \delta, \gamma)$ is large enough; $0< \gamma \leq \gamma_0$, $c(\varepsilon, \delta, \gamma) \geq c_0$ and $x \geq c(\varepsilon, \delta, \gamma)$. Let $q_0$ be the exceptional modulus in the interval $[3, \textup{exp} (2 c_1 \sqrt{\ln x})]$. If $q_0$ does not exist, then we take $B=1$. If $q_0$ exists, then (see \eqref{L4:q0_est})
\[
q_0 \geq \frac{a^{\prime}_{1}c_{1}^{2}\ln x}{(\ln\ln x)^{4}}=\frac{c_4\ln x}{(\ln\ln x)^{4}},
\]where $c_4 > 0$ is an absolute constant. We have $q_{0}\geq 24$, if $c(\varepsilon, \delta, \gamma)$ is chosen large enough. By Lemma \ref{L_q0_square_free}, the number $q_{0}$ is of the form $2^{\alpha}k$, where $\alpha\in\{0,\ldots,3\}$ and $k\geq 3$ is an odd square-free integer. We put
\[
M_{1}=\frac{q_{0}}{2^{\alpha}}\geq \frac{q_{0}}{8}\geq \frac{c_4\ln x}{8(\ln\ln x)^{4}}.
\] Let $\tau = (M_{1}, q)$, $M_{2}=M_{1}/\tau$. Then $(M_{2}, q)=1$. Since
\[
\tau \leq q \leq (\ln x)^{1-\varepsilon},
\]we have
\[
M_{2}=\frac{M_{1}}{\tau}\geq \frac{M_{1}}{(\ln x)^{1-\varepsilon}}\geq \frac{c_4\ln x}{8(\ln\ln x)^{4} (\ln x)^{1-\varepsilon}}=
\frac{c_4(\ln x)^{\varepsilon}}{8(\ln\ln x)^{4}}.
\]The number $M_{2} \geq 3$, if $c(\varepsilon, \delta, \gamma)$ is chosen large enough. Hence, $M_{2}\geq 3$ is an odd square-free integer. Furthermore, we have $(M_{2}, q)=1$ and $M_{2}$ divides $q_0$. Let $B$ to be the largest prime divisor of $M_{2}$. Hence, $B\geq 3$ is a prime number and $B$ divides $q_0$. We have (see Lemma \ref{L_about_Euler_func})
\[
\frac{B}{\varphi (B)} = \frac{B}{B(1-1/B)}=\frac{1}{1-1/B}\leq \frac{1}{1-1/3}=\frac{3}{2}.
\]Thus, $1 \leq B \leq \textup{exp} (2c_1 \sqrt{\ln x})$ is an integer, $(B,q)=1$,
\[
1\leq\frac{B}{\varphi (B)} \leq 2
\]and $B\geq 3$ is a prime divisor of $q_0$, if $q_0$ exists.

 Let $u$ be a real number with $2 \leq u \leq x^{1+\gamma/ \sqrt{\ln x}}$, let $Q$ and $W$ be integers with $2\leq Q \leq x^{1/2 - \delta}$, $(Q, B) = 1$ and $(W, Q) = 1$. By Lemma \ref{L_STAR_psi(u,Q,W)}, we have
\[
\psi (u; Q, W) - \frac{u}{\varphi (Q)} = \frac{1}{\varphi (Q)} \sum_{\chi \in X_Q} \overline{\chi(W)} \psi^{\prime} (u,\chi).
\]Hence,
\[
\Bigl|\psi (u; Q, W) - \frac{u}{\varphi (Q)}\Bigr| \leq \frac{1}{\varphi (Q)} \sum_{\chi \in X_Q} |\psi^{\prime} (u,\chi)|.
\]Since the right-hand side of this inequality does not depend on $W$, we have
\[
\max_{\substack{W\in \mathbb{Z}:\\
(W,Q)=1}}\Bigl|\psi (u; Q, W) - \frac{u}{\varphi (Q)}\Bigr| \leq \frac{1}{\varphi (Q)} \sum_{\chi \in X_Q} |\psi^{\prime} (u,\chi)|.
\]Let $\chi\in X_Q$, let $\chi_1$ be a primitive character modulo $q_1$ inducing $\chi$. From Lemma \ref{L_induced} and the definition of the inducing character, given below Lemma \ref{L_induced}, we have $q_1 = c(\chi)$, and hence $q_1| Q$ (see Lemma \ref{L_c(chi)_delit_q}). Applying Lemma \ref{L:chi_minus_chi1} we have
\[
|\psi^{\prime}(u, \chi)| \leq |\psi^{\prime}(u, \chi_1)| +\ln^{2} (Qu).
\]Since $\# X_Q = \varphi (Q)$, we obtain
\begin{align*}
\max_{\substack{W\in \mathbb{Z}:\\
(W,Q)=1}}\Bigl|\psi (u; Q, W) - \frac{u}{\varphi (Q)}\Bigr| &\leq \frac{1}{\varphi (Q)} \sum_{\chi \in X_Q} \bigl(|\psi^{\prime} (u,\chi_1)|+\ln^{2}(Qu)\bigr)=\\
&= \ln^{2}(Qu) + \frac{1}{\varphi (Q)} \sum_{\chi \in X_Q} |\psi^{\prime} (u,\chi_1)|.
\end{align*}We have
\begin{equation}\label{eq_<2}
1+ \frac{\gamma}{\sqrt{\ln x}} \leq 2,
\end{equation}if $c(\varepsilon, \delta, \gamma)$ is chosen large enough. Hence,
\begin{gather*}
0< \ln u \leq \Bigl(1+ \frac{\gamma}{\sqrt{\ln x}} \Bigr) \ln x \leq 2\ln x,\\
\ln^{2}u\leq 4\ln^{2}x,\\
0<\ln Q \leq \Bigl(\frac{1}{2} - \delta\Bigr)\ln x \leq \ln x,\\
\ln^{2}Q \leq \ln^{2}x,\\
\ln^{2}(Qu)\leq 2(\ln^{2}Q+\ln^{2}u)\leq 10 \ln^{2}x.
\end{gather*}We obtain
\begin{align*}
&\max_{2\leq u \leq x^{1+\gamma/\sqrt{\ln x}}}\max_{\substack{W\in \mathbb{Z}:\\
(W,Q)=1}}\Bigl|\psi (u; Q, W) - \frac{u}{\varphi (Q)}\Bigr| \leq \\
 &\leq 10\ln^{2}x +  \frac{1}{\varphi (Q)} \sum_{\chi \in X_Q} \max_{2\leq u \leq x^{1+\gamma/\sqrt{\ln x}}}   |\psi^{\prime} (u,\chi_1)| .
\end{align*}Hence,
\begin{align}
&S=\sum_{\substack{1\leq Q \leq x^{1/2 - \delta}\notag\\
(Q,B)=1}} A_Q = A_1 + \sum_{\substack{2\leq Q \leq x^{1/2 - \delta}\notag\\
(Q,B)=1}} A_Q \leq A_1 +\notag\\
&+ \sum_{\substack{2\leq Q \leq x^{1/2 - \delta}\notag\\
(Q,B)=1}}\Bigl(10\ln^{2}x +  \frac{1}{\varphi (Q)} \sum_{\chi \in X_Q} \max_{2\leq u \leq x^{1+\gamma/\sqrt{\ln x}}}   |\psi^{\prime} (u,\chi_1)|\Bigr)\leq\notag\\
&\leq 10 x^{1/2 - \delta}\ln^{2}x+ A_1+ \sum_{\substack{2\leq Q \leq x^{1/2 - \delta}\notag\\
(Q,B)=1}}  \sum_{\chi \in X_Q}\frac{1}{\varphi (Q)} \max_{2\leq u \leq x^{1+\gamma/\sqrt{\ln x}}}   |\psi^{\prime} (u,\chi_1)|=\notag\\
&= 10 x^{1/2 - \delta}\ln^{2}x+ A_1+ S^{\prime},\label{S_Super_Final}
\end{align}where
\[
A_Q:= \max_{2\leq u \leq x^{1+\gamma/\sqrt{\ln x}}}\max_{\substack{W\in \mathbb{Z}:\\
(W,Q)=1}}\Bigl|\psi (u; Q, W) - \frac{u}{\varphi (Q)}\Bigr|.
\]Let us estimate the sum $S^{\prime}$. Let $Q$ be an integer with $2 \leq Q \leq x^{1/2 - \delta}$ and $(Q,B)=1$, let $\chi\in X_Q$, let $\chi_1$ be the primitive character modulo $q_1$ inducing $\chi$. Since $q_1|Q$, we have $1\leq q_1 \leq x^{1/2 - \delta}$ and $(q_1,B)=1$. Hence,
\begin{align*}
&S^{\prime} = \sum_{\substack{2\leq Q \leq x^{1/2 - \delta}\\
(Q,B)=1}} \sum_{\chi \in X_Q}\frac{1}{\varphi (Q)}  \max_{2\leq u \leq x^{1+\gamma/\sqrt{\ln x}}}   |\psi^{\prime} (u,\chi_1)|\leq\\
&\leq \sum_{\substack{1\leq q_1 \leq x^{1/2 - \delta}\\
(q_1,B)=1}} \sum_{\chi_{1} \in X^{*}_{q_{1}}} \max_{2\leq u \leq x^{1+\gamma/\sqrt{\ln x}}}   |\psi^{\prime} (u,\chi_1)|
\Bigl(\sum_{1 \leq m \leq x^{1/2 - \delta}/q_1}\frac{1}{\varphi (m q_1)}\Bigr) .
\end{align*}Applying Lemmas \ref{L:Euler_func_ineq} and \ref{L:series_1/Euler_func} we have
\begin{align*}\
&\sum_{1 \leq m \leq x^{1/2 - \delta}/q_1}\frac{1}{\varphi (m q_1)}\leq \frac{1}{\varphi (q_1)}
\sum_{1 \leq m \leq x^{1/2 - \delta}/q_1}\frac{1}{\varphi (m)}\leq \\
&\leq \frac{1}{\varphi (q_1)}
\sum_{1 \leq m \leq x^{1/2}}\frac{1}{\varphi (m)}\leq \frac{1}{\varphi (q_1)} C\ln x,
\end{align*}where $C>0$ is an absolute constant. We obtain
\[
S^{\prime} \leq C\ln x \sum_{\substack{1\leq q_1 \leq x^{1/2 - \delta}\\
(q_1,B)=1}}\frac{1}{\varphi (q_1)} \sum_{\chi_{1} \in X^{*}_{q_{1}}} \max_{2\leq u \leq x^{1+\gamma/\sqrt{\ln x}}}   |\psi^{\prime} (u,\chi_1)| .
\]Replacing $q_1$ by $Q$ and $\chi_1$ by $\chi$, we have
\begin{align}
&S^{\prime} \leq C\ln x \sum_{\substack{1\leq Q \leq x^{1/2 - \delta}\notag\\
(Q,B)=1}}\frac{1}{\varphi (Q)} \sum_{\chi\in X^{*}_{Q}} \max_{2\leq u \leq x^{1+\gamma/\sqrt{\ln x}}}   |\psi^{\prime} (u,\chi)|=\notag\\
&=C\ln x \Bigl(\sum_{\substack{1\leq Q \leq \ln x\notag\\
(Q,B)=1}} R_{Q}  +  \sum_{\substack{\ln x < Q \leq \textup{exp} (c_1\sqrt{\ln x})\notag\\
(Q,B)=1}}R_{Q}  +   \sum_{\substack{\textup{exp} (c_1\sqrt{\ln x}) < Q \leq x^{1/2 - \delta}\notag\\
(Q,B)=1}} R_{Q}   \Bigr)=\notag\\
&=C\ln x\bigl(S^{\prime}_{1} + S^{\prime}_{2} +  S^{\prime}_{3} \bigr),\label{S_prime_final_est}
 \end{align}where
\[
R_{Q}:= \frac{1}{\varphi (Q)} \sum_{\chi\in X^{*}_{Q}} \max_{2\leq u \leq x^{1+\gamma/\sqrt{\ln x}}}   |\psi^{\prime} (u,\chi)|,
\]
$c_1>0$ is the absolute constant in Lemma \ref{L4}.

I) Now we estimate $S^{\prime}_{1}$. We have
\begin{equation}\label{S1_prime_estimate}
S^{\prime}_{1} = \sum_{\substack{1\leq Q \leq \ln x\\
(Q,B)=1}} R_{Q} \leq R_1 + \sum_{2\leq Q \leq \ln x} R_{Q} = R_1 + S^{\prime}_{4}.
\end{equation}

1) Now we estimate $R_1$. Since $\# X^{*}_{1} = 1$, we have
\[
R_1 = \max_{2\leq u \leq x^{1+\gamma/\sqrt{\ln x}}}   |\psi^{\prime} (u,\chi)|,
\]where $\chi\in X^{*}_{1}$, i.\,e. $\chi (n)=1$ for any $n\in \mathbb{Z}$. Since $\chi$ is the principal character modulo 1, we have
\[
\psi^{\prime}(u, \chi) = \psi (u,\chi) - u .
\]We have
\begin{gather*}
\psi (u,\chi) = \sum_{n\leq u} \Lambda (n)\chi(n)=\sum_{n\leq u}\Lambda(n)=\psi (u),\\
\psi^{\prime}(u, \chi) = \psi (u) - u .
\end{gather*}It is well-known that (see, for example, \mbox{\cite[Chapter 18]{Davenport}})
\begin{equation}\label{Psi_assymptotic_formula}
|\psi (u) - u| \leq C u\, \textup{exp} (-c\sqrt{\ln u}),\quad u \geq 2,
\end{equation}where $C>0$ and $c>0$ are absolute constants.

 Let us consider two cases.\\
 i) $x^{1/4} \leq u \leq x^{1+ \gamma/\sqrt{\ln x}}$ (we may assume that $c(\varepsilon, \delta, \gamma)>16$). We have (see \eqref{eq_<2})
 \begin{gather*}
 \frac{1}{4} \ln x \leq \ln u \leq \Bigl(1+ \frac{\gamma}{\sqrt{\ln x}}\Bigr)\ln x \leq 2 \ln x,\\
 -c\sqrt{\ln u} \leq -\frac{c}{2} \sqrt{\ln x}.
 \end{gather*}Hence,
 \begin{align*}
 &|\psi^{\prime}(u,\chi)|\leq C u\, \textup{exp} (-c\sqrt{\ln u})\leq C x^{1+\gamma/\sqrt{\ln x}}\, \textup{exp} \bigl(-(c/2)\sqrt{\ln x}\bigr)=\\
 &= Cx\,\textup{exp} \bigl((\gamma-c/2)\sqrt{\ln x}\bigr) \leq Cx\,\textup{exp} \bigl((-c/4)\sqrt{\ln x}\bigr),
 \end{align*}if $0< \gamma \leq c/4$.\\
 ii) $2 \leq u < x^{1/4}$. Then
 \begin{align*}
 |\psi^{\prime}(u,\chi)|&\leq C u\, \textup{exp} (-c\sqrt{\ln u})\leq Cu\leq \\
  & \leq C x^{1/4}\leq Cx\,\textup{exp} \bigl((-c/4)\sqrt{\ln x}\bigr),
 \end{align*} if $c(\varepsilon, \delta, \gamma)$ is chosen large enough.

 We obtain
 \begin{equation}\label{R1_est}
 R_1 = \max_{2\leq u \leq x^{1+\gamma/\sqrt{\ln x}}}   |\psi^{\prime} (u,\chi)|\leq Cx\,\textup{exp} \bigl((-c/4)\sqrt{\ln x}\bigr) .
 \end{equation}

2) Now we estimate
\begin{equation}\label{def_S_prime_4}
S^{\prime}_{4} = \sum_{2\leq Q \leq \ln x} \frac{1}{\varphi (Q)} \sum_{\chi\in X^{*}_{Q}} \max_{2\leq u \leq x^{1+\gamma/\sqrt{\ln x}}}   |\psi^{\prime} (u,\chi)|.
\end{equation}Let $Q$ be an integer with $2 \leq Q \leq \ln x$, and let $\chi \in X^{*}_{Q}$. Then $\chi$ is a nonprincipal character modulo $Q$, and hence
\[
\psi^{\prime}(u, \chi) = \psi(u,\chi).
\]Let us consider two cases.\\
i) $x^{1/4} \leq u \leq x^{1+ \gamma/\sqrt{\ln x}}$. Then (see \eqref{eq_<2})
\[
\frac{1}{4} \ln x  \leq \ln u\leq \Bigl(1+\frac{\gamma}{\sqrt{\ln x}}\Bigr)\ln x \leq 2 \ln x.
\]We may assume that $c(\varepsilon, \delta, \gamma) \geq e^{16}$. Hence, $\ln u \geq (\ln x)/4\geq 4$. We have
\[
2 \leq Q \leq \ln x \leq 4 \ln u \leq \ln^{2} u.
\]Therefore (see, for example, \mbox{\cite[Chapter 22]{Davenport}})
\[
|\psi (u,\chi)| \leq C u\,\textup{exp}(-c(2) \sqrt{\ln u}),
\]where $C>0$ and $c(2)>0$ are absolute constants. We have
\[
-c(2)\sqrt{\ln u} \leq -\frac{c(2)}{2}\sqrt{\ln x},
\]
\begin{align*}
|\psi(u,\chi)| &\leq C u\,\textup{exp}\bigl(-(c(2)/2) \sqrt{\ln x}\bigr) \leq\\
 &\leq C x^{1+\gamma/\sqrt{\ln x}}\,\textup{exp}\bigl(-(c(2)/2) \sqrt{\ln x}\bigr)=\\
 &= C x\,\textup{exp}\bigl((\gamma-c(2)/2) \sqrt{\ln x}\bigr)\leq\\
 & \leq Cx\,\textup{exp}\bigl(-(c(2)/4) \sqrt{\ln x}\bigr),
\end{align*}if $0< \gamma \leq c(2)/4$.\\
ii) $2 \leq u < x^{1/4}$. We have (see \eqref{Notation:psi_est})
\[
\psi(u,\chi) = \sum_{n \leq u} \Lambda(n) \chi(n),
\]
\begin{align*}
|\psi(u,\chi)|&\leq \sum_{n \leq u} \Lambda(n) = \psi(u)\leq b_6 u \leq b_6 x^{1/4}\leq\\
&\leq Cx\,\textup{exp}\bigl(- (c(2)/4)\sqrt{\ln x}\bigr),
\end{align*}if $c(\varepsilon, \delta, \gamma)$ is chosen large enough.

Hence,
\[
\max_{2 \leq u \leq x^{1+\gamma/\sqrt{\ln x}}} |\psi (u,\chi)| \leq Cx\,\textup{exp}\bigl(- (c(2)/4)\sqrt{\ln x}\bigr).
\] Substituting this estimate into \eqref{def_S_prime_4} and using the fact that $\#X^{*}_{Q}\leq \#X_{Q}=\varphi (Q)$, we obtain
\begin{align}
S^{\prime}_{4} &\leq Cx\,\textup{exp}\bigl(- (c(2)/4)\sqrt{\ln x}\bigr)\ln x=\notag\\
&= Cx\,\textup{exp}\bigl(- (c(2)/4)\sqrt{\ln x} + \ln\ln x\bigr)\leq\notag\\
&\leq Cx\,\textup{exp}\bigl(- (c(2)/8)\sqrt{\ln x}\bigr),\label{S4_prime_est}
\end{align}if $c(\varepsilon, \delta, \gamma)$ is chosen large enough.

Substituting \eqref{R1_est} and \eqref{S4_prime_est} into \eqref{S1_prime_estimate}, we obtain
\begin{equation}\label{S_1_prime}
S^{\prime}_{1} \leq Cx\,\textup{exp}(-c\sqrt{\ln x}),
\end{equation}where $C>0$ and $c>0$ are absolute constants.

II) Now we estimate
\[
S^{\prime}_{3} = \sum_{\substack{\textup{exp} (c_1\sqrt{\ln x}) < Q \leq x^{1/2 - \delta}\\
(Q,B)=1}}   \frac{1}{\varphi (Q)} \sum_{\chi\in X^{*}_{Q}} \max_{2\leq u \leq x^{1+\gamma/\sqrt{\ln x}}}   |\psi^{\prime} (u,\chi)| .
 \]Let $Q$ be an integer with $\textup{exp} (c_1\sqrt{\ln x}) < Q \leq x^{1/2 - \delta}$ and $(Q,B)=1$, and let $\chi\in X^{*}_{Q}$. Since $Q>1$, we see that $\chi$ is a nonprincipal character modulo $Q$. Hence,
 \[
 \psi^{\prime}(u,\chi) = \psi(u,\chi).
 \]We obtain
 \begin{align*}
 &S^{\prime}_{3} = \sum_{\substack{\textup{exp} (c_1\sqrt{\ln x}) < Q \leq x^{1/2 - \delta}\\
(Q,B)=1}}   \frac{1}{\varphi (Q)} \sum_{\chi\in X^{*}_{Q}} \max_{2\leq u \leq x^{1+\gamma/\sqrt{\ln x}}}   |\psi (u,\chi)|\leq\\
&\leq \sum_{\textup{exp} (c_1\sqrt{\ln x}) < Q \leq x^{1/2 - \delta}}   \frac{1}{\varphi (Q)} \sum_{\chi\in X^{*}_{Q}} \max_{2\leq u \leq x^{1+\gamma/\sqrt{\ln x}}}   |\psi (u,\chi)| .
 \end{align*}Applying Lemma \ref{L:Dave_character} with $Q_1 = \textup{exp} (c_1\sqrt{\ln x})$, $Q_2= x^{1/2 - \delta}$ and $t=x^{1+\gamma/\sqrt{\ln x}}$, we obtain
 \begin{align*}
 &S^{\prime}_{3} \leq C\ln^{4}\bigl(x^{3/2 - \delta +\gamma/\sqrt{\ln x}}\bigr)\Bigl(x\,\textup{exp}\bigl((\gamma - c_1)\sqrt{\ln x}\bigr)+\\
 &+ x^{(5/6)(1+\gamma/\sqrt{\ln x})} \ln (x^{1/2 - \delta})+ x^{1-\delta+ \gamma/(2\sqrt{\ln x})}\Bigr) .
 \end{align*}We have
 \[
 \frac{\gamma}{\sqrt{\ln x}} \leq \delta,\quad \frac{5}{6}\Bigl(1+\frac{\gamma}{\sqrt{\ln x}}\Bigr)\leq 0.9,
 \]if $c(\varepsilon, \delta, \gamma)$ is chosen large enough. Replacing $C$, we have
  \[
 S^{\prime}_{3} \leq C\ln^{4} x\Bigl(x\,\textup{exp}\bigl((\gamma - c_1)\sqrt{\ln x}\bigr)+ x^{0.9} \ln x+ x^{1-\delta / 2}\Bigr).
 \]We have
 \[
 (\gamma - c_1)\sqrt{\ln x} \leq -\frac{c_1}{2}\sqrt{\ln x},
 \]if $0< \gamma\leq c_1 / 2$. We obtain
 \begin{gather*}
 x\,\textup{exp}\bigl((\gamma - c_1)\sqrt{\ln x}\bigr)\ln^{4}x \leq x\,\textup{exp}\bigl(- (c_1/2)\sqrt{\ln x}+ 4\ln\ln x\bigr)\leq\\
 \leq x\,\textup{exp}\bigl(- (c_{1}/4)\sqrt{\ln x}\bigr),\\
 x^{0.9}\ln^{5}x\leq x\,\textup{exp}\bigl(- (c_{1}/4)\sqrt{\ln x}\bigr),\\
 x^{1-\delta / 2}\ln^{4} x \leq x\,\textup{exp}\bigl(- (c_{1}/4)\sqrt{\ln x}\bigr),
  \end{gather*}if $c(\varepsilon, \delta, \gamma)$ is chosen large enough. Replacing $3C$ by $C$ and $c_{1}/4$ by $c$, we obtain
  \begin{equation}\label{S_3_prime}
  S^{\prime}_{3} \leq Cx\,\textup{exp}(-c\sqrt{\ln x}),
  \end{equation}where $C>0$ and $c>0$ are absolute constants.

  III) Now we estimate
  \[
  S^{\prime}_{2}= \sum_{\substack{\ln x < Q \leq \textup{exp} (c_1\sqrt{\ln x})\\
(Q,B)=1}}\frac{1}{\varphi (Q)} \sum_{\chi\in X^{*}_{Q}} \max_{2\leq u \leq x^{1+\gamma/\sqrt{\ln x}}}   |\psi^{\prime} (u,\chi)| .
  \]Let $Q$ be an integer with $\ln x < Q \leq \textup{exp} (c_1\sqrt{\ln x})$ and $(Q,B)=1$, and let $\chi\in X^{*}_{Q}$. Since $Q>1$, we see that $\chi$ is a nonprincipal character modulo $Q$, and hence $\psi^{\prime} (u,\chi)= \psi(u,\chi)$. We recall that if the exceptional modulus $q_0$ in the interval $[3, \textup{exp}(2c_{1}\sqrt{\ln x})]$ does not exist, then $B=1$; if $q_0$ exists, then $B \geq3$ is a prime divisor of $q_0$, and hence $Q\neq q_0$. Since $0<\gamma \leq \gamma_0$ and $c(\varepsilon, \delta, \gamma) \geq c_0$, we see from Lemma \ref{L4} that
  \[
  \max_{2 \leq u \leq x^{1+\gamma/\sqrt{\ln x}}} |\psi (u,\chi)| \leq Cx\,\textup{exp}(-3c_{1}\sqrt{\ln x}).
  \]Since $\#X^{*}_{Q}\leq \#X_{Q}=\varphi (Q)$, we obtain
  \begin{align}
  S^{\prime}_{2} &\leq \sum_{\substack{\ln x < Q \leq \textup{exp} (c_1\sqrt{\ln x})\notag\\
(Q,B)=1}}  Cx\,\textup{exp}(-3c_{1}\sqrt{\ln x})\leq\notag\\
&\leq Cx\,\textup{exp}(-3c_{1}\sqrt{\ln x}) \textup{exp} (c_1\sqrt{\ln x})=\notag\\
&= Cx\,\textup{exp}(-2c_{1}\sqrt{\ln x}).\label{S2_prime}
  \end{align}

  From \eqref{S_1_prime}, \eqref{S_3_prime} and \eqref{S2_prime} we obtain
  \begin{equation}\label{S1+S2+S3_prime}
  S^{\prime}_{1} +  S^{\prime}_{2}+  S^{\prime}_{3} \leq  \widetilde{C}x\,\textup{exp}(-\widetilde{c}\sqrt{\ln x}),
  \end{equation}where $\widetilde{C}>0$ and $\widetilde{c}> 0$ are absolute constants. Substituting \eqref{S1+S2+S3_prime} into \eqref{S_prime_final_est} we obtain
  \[
  S^{\prime} \leq C^{\prime}x\,\textup{exp}(-\widetilde{c}\sqrt{\ln x} + \ln\ln x)\leq
  C^{\prime}x\,\textup{exp}\bigl(-(\widetilde{c}/2)\sqrt{\ln x} \bigr),
  \]if $c(\varepsilon, \delta, \gamma)$ is chosen large enough. Replacing $C^{\prime}$ by $C$ and $\widetilde{c}/2$ by $c$, we obtain
  \begin{equation}\label{S_prime_absolute_fin}
  S^{\prime} \leq Cx\,\textup{exp}(-c\sqrt{\ln x}),
  \end{equation}where $C>0$ and $c>0$ are absolute constants.

  IV) We have
  \begin{equation}\label{x_delta_absolute}
  x^{1/2 - \delta}\ln^{2}x\leq x^{1/2}\ln^{2}x \leq x\,\textup{exp}(-c\sqrt{\ln x}),
  \end{equation}if $c(\varepsilon, \delta, \gamma)$ is chosen large enough (here $c>0$ is the absolute constant in \eqref{S_prime_absolute_fin}).

  V) Now we estimate
  \[
  A_1= \max_{2\leq u \leq x^{1+\gamma/\sqrt{\ln x}}}\max_{W\in \mathbb{Z}}\bigl|\psi (u; 1, W) - u\bigr|.
    \] Let $W\in \mathbb{Z}$. We have
   \[
    \psi (u; 1, W)=\sum_{\substack{n\leq u\\ n\equiv W\,(\textup{mod} 1)}}\Lambda (n) =
    \sum_{n\leq u} \Lambda(n) = \psi (u).
    \]Hence,
    \[
    A_1 = \max_{2\leq u \leq x^{1+\gamma/\sqrt{\ln x}}} | \psi (u) - u|.
    \]Using \eqref{Psi_assymptotic_formula} and arguing as in I, 1, i) and I, 1, ii), we obtain
    \begin{equation}\label{A1_absolute}
    A_1 \leq Cx\,\textup{exp}(-c\sqrt{\ln x}),
    \end{equation}where $C>0$ and $c>0$  are absolute constants.

    Substituting \eqref{S_prime_absolute_fin}, \eqref{x_delta_absolute} and \eqref{A1_absolute} into \eqref{S_Super_Final}, we obtain
    \[
    S\leq Cx\,\textup{exp}(-c\sqrt{\ln x}),
    \]where $C>0$ and $c>0$ are absolute constants. Thus, if $\gamma$ is sufficiently small positive absolute constant, $x$ is a real number with $x \geq c(\varepsilon, \delta, \gamma)$, $q$ is an integer with $1\leq q \leq (\ln x)^{1-\varepsilon}$, then there is an integer $B$ such that
    \[
    1\leq B \leq \textup{exp}(2 c_1\sqrt{\ln x}),\quad 1\leq \frac{B}{\varphi(B)}\leq 2,\quad (B,q)=1
    \]and
    \[
    \sum_{\substack{1\leq Q\leq x^{1/2 - \delta}\\
     (Q,B)=1}} \max_{2 \leq u \leq x^{1+ \gamma/\sqrt{\ln x}}}\max_{\substack{W\in \mathbb{Z}:\\ (W,Q)=1}}
    \Bigl|\psi (u; Q, W)-\frac{u}{\varphi (Q)}\Bigr| \leq C x\,\textup{exp}(-c\sqrt{\ln x}),
    \]where $c_1$, $C$ and $c$ are positive absolute constants. Let us denote $2c_1$ by $c_1$, $C$ by $c_2$, $c$ by $c_3$. Since $\gamma$ is an absolute constant, we see that the positive number $c(\varepsilon, \delta, \gamma)= c(\varepsilon, \delta)$ depends only on $\varepsilon$ and $\delta$. Lemma \ref{L5_Bombieri} is proved.

    \begin{lemma}\label{L6_Bombieri_Pi}
    Let $\varepsilon$ and $\delta$ be real numbers with $0< \varepsilon < 1$ and $0< \delta < 1/2$. Then there is a number $c(\varepsilon, \delta)>0$, depending only on $\varepsilon$ and $\delta$, such that if $x\in \mathbb{R}$ and $q\in \mathbb{Z}$ are such that $x\geq c(\varepsilon, \delta)$ and $1 \leq q \leq (\ln x)^{1-\varepsilon}$, then there is a positive integer $B$ such that
  \[
1\leq B \leq \textup{exp} (c_1 \sqrt{\ln x}),\quad 1\leq \frac{B}{\varphi(B)}\leq 2,\quad (B,q)=1
\]and
\[
\sum_{\substack{1\leq Q\leq x^{1/2 - \delta}\\
(Q,B)=1}} \max_{2 \leq u \leq x^{1+ \gamma/\sqrt{\ln x}}}\max_{\substack{W\in \mathbb{Z}:\\ (W,Q)=1}}
\Bigl|\pi (u; Q, W)-\frac{li(u)}{\varphi (Q)}\Bigr| \leq c_2 x\,\textup{exp}(-c_3\sqrt{\ln x}).
\]Here $c_1$, $\gamma$, $c_2$ and $c_3$ are positive absolute constants.
    \end{lemma}
    \textsc{Proof of Lemma \ref{L6_Bombieri_Pi}.} We choose $\widetilde{c}(\varepsilon, \delta)$ later; this number is large enough. Let $\widetilde{c}(\varepsilon, \delta) \geq c(\varepsilon, \delta)$, where $c(\varepsilon, \delta)$ is the number in Lemma \ref{L5_Bombieri}. Let $x\in\mathbb{R}$ and $q\in\mathbb{Z}$ be such that $x \geq \widetilde{c}(\varepsilon, \delta)$ and $1 \leq q \leq (\ln x)^{1-\varepsilon}$. Then, by Lemma \ref{L5_Bombieri}, there is a positive integer $B$ such that
    \begin{equation}\label{B_for_L6}
1\leq B \leq \textup{exp} (c_1 \sqrt{\ln x}),\quad 1\leq \frac{B}{\varphi(B)}\leq 2,\quad (B,q)=1
\end{equation}and
\begin{equation}\label{psi_c2_c3_ineq_for_L6}
\sum_{\substack{1\leq Q\leq x^{1/2 - \delta}\\
(Q,B)=1}} \max_{2 \leq u \leq x^{1+ \gamma/\sqrt{\ln x}}}\max_{\substack{W\in \mathbb{Z}:\\ (W,Q)=1}}
|R (u; Q, W)| \leq c_2 x\,\textup{exp}(-c_{3}\sqrt{\ln x}),
\end{equation}where
\[
R (u; Q, W):= \psi(u; Q, W) - \frac{u}{\varphi(Q)}
\]and $c_1$, $\gamma$, $c_2$ and $c_3$ are positive absolute constants.

We put
\begin{equation}\label{def_R1_for_L6}
R_{1}(u; Q, W) := \pi (u; Q, W) - \frac{li (u)}{\varphi(Q)}.
\end{equation}Let $Q\in \mathbb{Z}$, $W\in \mathbb{Z}$ and $u\in \mathbb{Z}$ be such that $1 \leq Q \leq x^{1/2 - \delta}$, $(Q,B)=1$, $(W,Q)=1$, $3\leq u \leq x^{1+\gamma/\sqrt{\ln x}}$. We claim that
\begin{equation}\label{R1_Super_Estim}
|R_{1}(u; Q, W)|\leq C_{1} u^{1/2} + |R(u; Q, W)|+ \sum_{2 \leq n \leq u-1} \frac{|R(n; Q, W)|}{n\ln^{2}n},
\end{equation}where $C_{1}>0$ is an absolute constant. We define
\[
\alpha (n) = \begin{cases}
1, &\text{if $n\equiv W$ (mod $Q$),}\\
0, &\text{otherwise};
\end{cases}
\]
\[
\pi_{1}(u; Q, W) = \sum_{n\leq u}\frac{\Lambda (n) \alpha (n)}{\ln n}.
\]Let us show that
\begin{equation}\label{Pi_minus_P_1_eq}
\pi (u; Q, W) = \pi_{1}(u; Q, W) + \widetilde{R}(u; Q, W),\quad |\widetilde{R}(u; Q, W)| \leq C u^{1/2},
\end{equation}where $C>0$ is an absolute constant. Let $u \geq 8$. Then
\begin{align*}
&\pi_{1}(u; Q, W) = \sum_{p^{m} \leq u}\frac{ \alpha (p^{m})\ln p}{m\ln p}=\sum_{1 \leq m \leq \ln u/ \ln 2} \sum_{p \leq u^{1/m}}\frac{\alpha (p^m)}{m}=\\
&=\sum_{p \leq u} \alpha (p) + \sum_{2 \leq m \leq \ln u/ \ln 2}\frac{1}{m} \sum_{p \leq u^{1/m}}\alpha (p^m)=
S_1 + S_2.
\end{align*}We have
\[
S_1 = \sum_{\substack{p\leq u\\ p\equiv W\, \text{(mod $Q$)}}} 1 = \pi(u; Q, W),
\]
\begin{align*}
&S_2\leq  \sum_{2 \leq m \leq \ln u/ \ln 2}\frac{u^{1/m}}{m} = \frac{1}{2} u^{1/2} + \sum_{3 \leq m \leq \ln u/ \ln 2}\frac{u^{1/m}}{m}\leq\\
&\leq \frac{1}{2} u^{1/2} + \frac{1}{3} u^{1/3}\frac{\ln u}{\ln 2}\leq u^{1/2} + u^{1/3} \ln u \leq C^{\prime} u^{1/2},
\end{align*}where $C^{\prime}>0$ is an absolute constant. If $3 \leq u <8$, then
\[
\Bigl|\sum_{2 \leq m \leq \ln u/ \ln 2}\frac{1}{m} \sum_{p \leq u^{1/m}}\alpha (p^m)\Bigr|
\leq \frac{1}{2} \sum_{p \leq 8^{1/2}}1 + \frac{1}{3} \sum_{p \leq 8^{1/3}}1 = C^{\prime \prime} \leq C^{\prime \prime} u^{1/2}.
\]Thus, \eqref{Pi_minus_P_1_eq} is proved.

Since
\[
\psi(x; Q, W) = \sum_{m\leq x}\Lambda (m)\alpha(m),
\] we have
\begin{align*}
&\pi_{1}(u; Q, W) = \sum_{2 \leq n \leq u} \frac{\psi (n; Q, W) - \psi(n-1; Q, W)}{\ln n}=\\
&=\sum_{2 \leq n \leq u-1}\psi (n; Q, W)\Bigl(\frac{1}{\ln n} - \frac{1}{\ln (n+1)}\Bigr)
+\frac{\psi(u; Q, W)}{\ln u}=\\
&= \sum_{2 \leq n \leq u-1}\Bigl(\frac{n}{\varphi(Q)}+ R(n; Q, W)\Bigr)\Bigl(\frac{1}{\ln n} - \frac{1}{\ln (n+1)}\Bigr)
+\frac{u}{\varphi(Q)\ln u}+\\
&+ \frac{R(u; Q, W)}{\ln u}.
\end{align*}We have
\begin{align*}
&\sum_{2 \leq n \leq u-1}\frac{n}{\varphi (Q)}\Bigl(\frac{1}{\ln n} - \frac{1}{\ln (n+1)}\Bigr)=
\sum_{2 \leq n \leq u-1}\frac{n}{\varphi (Q)}\int_{n}^{n+1}\frac{dt}{t\ln^{2}t}=\\
&=\frac{1}{\varphi (Q)}\sum_{2 \leq n \leq u-1}\int_{n}^{n+1}\frac{t - \{t\} }{t\ln^{2}t}\,dt=
\frac{1}{\varphi (Q)}\Bigl(\int_{2}^{u}\frac{dt}{\ln^{2}t} - \int_{2}^{u}\frac{\{t\}dt}{t\ln^{2}t} \Bigr).
\end{align*}Since
\begin{align*}
&\int_{2}^{u}\frac{dt}{\ln^{2}t}=\int_{2}^{u}t\, d\Bigl(-\frac{1}{\ln t}\Bigr)=\left. -\frac{t}{\ln t}\right|_{2}^{u}+
\int_{2}^{u}\frac{dt}{\ln t}=\\
&=-\frac{u}{\ln u} + \frac{2}{\ln 2} + li (u),
\end{align*}we obtain
\begin{align*}
&\pi_{1}(u; Q, W) = \frac{u}{\varphi(Q)\ln u}+ \frac{R(u; Q, W)}{\ln u} -\frac{u}{\varphi(Q)\ln u} + \frac{2}{\varphi(Q)\ln 2} +\\
&+ \frac{li (u)}{\varphi(Q)}- \frac{1}{\varphi (Q)}\int_{2}^{u}\frac{\{t\}}{t\ln^{2}t}\,dt+
\sum_{2 \leq n \leq u-1}R(n; Q, W)\Bigl(\frac{1}{\ln n} - \frac{1}{\ln (n+1)}\Bigr).
\end{align*} We have (see \eqref{Pi_minus_P_1_eq})
\[
\pi(u; Q, W) = \frac{li(u)}{\varphi(Q)}+ R_{1}(u; Q, W),
\]where
\begin{align*}
&R_{1}(u; Q, W) = \frac{2}{\varphi(Q)\ln 2} - \frac{1}{\varphi (Q)}\int_{2}^{u}\frac{\{t\}}{t\ln^{2}t}\,dt + \widetilde{R}(u; Q, W) +\\
&+  \frac{R(u; Q, W)}{\ln u} +\sum_{2 \leq n \leq u-1}R(n; Q, W)\Bigl(\frac{1}{\ln n} - \frac{1}{\ln (n+1)}\Bigr).
\end{align*}We obtain
\begin{align}
&|R_{1}(u; Q, W)| \leq \frac{2}{\ln 2}+ \Bigl|\int_{2}^{u}\frac{\{t\}}{t\ln^{2}t}\,dt\Bigr|+ |\widetilde{R}(u; Q, W)| +\notag\\
&+  \frac{|R(u; Q, W)|}{\ln u} +\sum_{2 \leq n \leq u-1}|R(n; Q, W)|\Bigl(\frac{1}{\ln n} - \frac{1}{\ln (n+1)}\Bigr).\label{R1_General_Estimate}
\end{align}Since $u \geq 3$, we have
\begin{equation}\label{R/ln_u_trivial_est}
\frac{|R(u; Q, W)|}{\ln u} \leq |R(u; Q, W)|.
\end{equation}Since
\[
\Bigl|\int_{2}^{u}\frac{\{t\}}{t\ln^{2}t}\,dt\Bigr|\leq \int_{2}^{u}\frac{dt}{t\ln^{2}t}=\left. -\frac{1}{\ln t}\right|_{2}^{u}=\frac{1}{\ln 2} - \frac{1}{\ln u}\leq \frac{1}{\ln 2},
\]we have (see \eqref{Pi_minus_P_1_eq})
\begin{equation}\label{R_with_wave_est<<u_1/2}
\frac{2}{\ln 2}+ \Bigl|\int_{2}^{u}\frac{\{t\}}{t\ln^{2}t}\,dt\Bigr|+ |\widetilde{R}(u; Q, W)|\leq
\frac{3}{\ln 2} + C u^{1/2}\leq \Bigl(C+\frac{3}{\ln 2}\Bigr)u^{1/2}.
\end{equation}Let $f(x)=-\ln^{-1}x$ and let $n$ be an integer with $n \geq 2$. By the mean value theorem, there is $\xi\in (n, n+1)$ such that
\begin{equation}\label{1/ln_n_Lagrange}
\frac{1}{\ln n} - \frac{1}{\ln (n+1)} = f(n+1) - f(n) = f^{\prime}(\xi)=\frac{1}{\xi \ln^{2}\xi}\leq \frac{1}{n\ln^{2}n}.
\end{equation}Substituting \eqref{R/ln_u_trivial_est}, \eqref{R_with_wave_est<<u_1/2} and \eqref{1/ln_n_Lagrange} into \eqref{R1_General_Estimate}, we obtain \eqref{R1_Super_Estim}. Hence,
\begin{align}
&\sum_{\substack{1 \leq Q \leq x^{1/2 - \delta}\notag\\ (Q,B)=1}}\max_{\substack{3 \leq u \leq x^{1+\gamma/\sqrt{\ln x}}\notag\\ u\in \mathbb{Z}}}\max_{\substack{W\in \mathbb{Z}\notag\\ (W,Q)=1}}|R_{1}(u; Q, W)|\leq\notag\\
&\leq \sum_{\substack{1 \leq Q \leq x^{1/2 - \delta}\notag\\ (Q,B)=1}}\max_{\substack{3 \leq u \leq x^{1+\gamma/\sqrt{\ln x}}\notag\\ u\in \mathbb{Z}}}\max_{\substack{W\in \mathbb{Z}\notag\\ (W,Q)=1}} |R(u; Q, W)|  +\notag\\
 &+\sum_{\substack{1 \leq Q \leq x^{1/2 - \delta}\notag\\ (Q,B)=1}}\max_{\substack{3 \leq u \leq x^{1+\gamma/\sqrt{\ln x}}\notag\\ u\in \mathbb{Z}}}\max_{\substack{W\in \mathbb{Z}\notag\\ (W,Q)=1}} |C_{1} u^{1/2}| +\notag\\
  &+ \sum_{\substack{1 \leq Q \leq x^{1/2 - \delta}\notag\\ (Q,B)=1}}\max_{\substack{3 \leq u \leq x^{1+\gamma/\sqrt{\ln x}}\notag\\ u\in \mathbb{Z}}}\max_{\substack{W\in \mathbb{Z}\notag\\ (W,Q)=1}} \sum_{2 \leq n \leq u-1} \frac{|R(n; Q, W)|}{n\ln^{2}n}=\notag\\
  &= S_1 + S_2 + S_3.\label{S1+S2+S3_for_L6}
\end{align}\newline
I) Now we estimate $S_1$. We have (see \eqref{psi_c2_c3_ineq_for_L6})
\begin{equation}\label{S1_for_L6}
S_1 \leq \sum_{\substack{1\leq Q\leq x^{1/2 - \delta}\\
(Q,B)=1}} \max_{2 \leq u \leq x^{1+ \gamma/\sqrt{\ln x}}}\max_{\substack{W\in \mathbb{Z}\\ (W,Q)=1}}
|R (u; Q, W)| \leq c_2 x\,\textup{exp}(-c_{3}\sqrt{\ln x}).
\end{equation}\\ II) Now we estimate $S_2$. We have
\[
\frac{\gamma}{\sqrt{\ln x}} \leq \delta,
\]if $\widetilde{c}(\varepsilon, \delta)$ is chosen large enough. We have
\begin{equation}\label{S2_for_L6}
S_2 \leq C_1 x^{1- \delta + \gamma/(2\sqrt{\ln x})}\leq C_1 x^{1- \delta / 2} \leq x\,\textup{exp}(-c_{3}\sqrt{\ln x}),
\end{equation}if $\widetilde{c}(\varepsilon, \delta)$ is chosen large enough.\\
III) Now we estimate $S_3$. Let $Q$, $W$, $u$ and $n$ be integers such that $1 \leq Q \leq x^{1/2 - \delta}$, $(Q,B)=1$, $(W,Q)=1$, $3 \leq u \leq x^{1+ \gamma/\sqrt{\ln x}}$  and $2\leq n \leq u-1$. Then
\[
|R(n; Q, W)| \leq \max_{2 \leq m \leq x^{1+\gamma/\sqrt{\ln x}}}\max_{\substack{V\in \mathbb{Z}\\ (V,Q)=1}}|R(m; Q, V)| .
\]Hence,
\begin{align*}
&\sum_{2 \leq n \leq u-1} \frac{|R(n; Q, W)|}{n\ln^{2}n} \leq \max_{2 \leq m \leq x^{1+\gamma/\sqrt{\ln x}}}\max_{\substack{V\in \mathbb{Z}\\ (V,Q)=1}}|R(m; Q, V)|\cdot\\
&\cdot\sum_{2 \leq n \leq u-1} \frac{1}{n\ln^{2}n} \leq c_{0} \max_{2 \leq m \leq x^{1+\gamma/\sqrt{\ln x}}}
\max_{\substack{V\in \mathbb{Z}\\ (V,Q)=1}}|R(m; Q, V)|,
\end{align*} where
\[
c_{0}:= \sum_{n=2}^{\infty}\frac{1}{n\ln^{2}n}<+\infty.
\]We have
\begin{align*}
&\max_{\substack{3 \leq u \leq x^{1+\gamma/\sqrt{\ln x}}\\ u\in \mathbb{Z}}} \max_{\substack{W\in \mathbb{Z}\\ (W, Q)=1}}
\sum_{2 \leq n \leq u-1} \frac{|R(n; Q, W)|}{n\ln^{2}n} \leq\\
 &\leq c_{0} \max_{2 \leq m \leq x^{1+\gamma/\sqrt{\ln x}}}\max_{\substack{V\in \mathbb{Z}\\ (V,Q)=1}}|R(m; Q, V)|.
\end{align*}We obtain (see \eqref{psi_c2_c3_ineq_for_L6})
\begin{align}
&S_3 \leq c_0 \sum_{\substack{1 \leq Q \leq x^{1/2 - \delta}\notag\\ (Q,B)=1}}\max_{2 \leq m \leq x^{1+\gamma/\sqrt{\ln x}}}\max_{\substack{V\in \mathbb{Z}\\ (V,Q)=1}}|R(m; Q, V)|\leq\notag\\
&\leq c_{0}c_{2}x\,\textup{exp}(-c_3 \sqrt{\ln x}).\label{S3_for_L6}
\end{align}Substituting \eqref{S1_for_L6}, \eqref{S2_for_L6} and \eqref{S3_for_L6} into \eqref{S1+S2+S3_for_L6}, we obtain (see \eqref{def_R1_for_L6})
\begin{align}
&\sum_{\substack{1 \leq Q \leq x^{1/2 - \delta}\notag\\ (Q,B)=1}}\max_{\substack{3 \leq u \leq x^{1+\gamma/\sqrt{\ln x}}\notag\\ u\in \mathbb{Z}}}\max_{\substack{W\in \mathbb{Z}\notag\\ (W,Q)=1}}\Bigl|\pi(u; Q, W)- \frac{li (u)}{\varphi (Q)}\Bigr| \leq\\
&\leq c_{4}x\,\textup{exp}(-c_3 \sqrt{\ln x}),\label{Discrete_Ineq_for_L6}
\end{align}where $c_4 = c_2+1+ c_{0} c_{2}> 0$ is an absolute constant.

 Let $Q$ and $W$ be integers such that $1 \leq Q \leq x^{1/2 - \delta}$, $(Q,B)=1$ and $(W,Q)=1$, and let $u$ be a real number with $2 \leq u \leq x^{1+\gamma/\sqrt{\ln x}}$. Let us consider two cases.\\
1) $2 \leq u \leq 3$. Then
\begin{gather}
|\pi(u; Q, W)| \leq \pi(u) \leq 2,\notag\\
\left|\frac{li (u)}{\varphi(Q)}\right|\leq li (u) \leq li (3),\notag\\
\left| \pi(u; Q, W) - \frac{li (u)}{\varphi(Q)} \right| \leq
|\pi(u; Q, W)| + \left|\frac{li (u)}{\varphi(Q)}\right| \leq 2+ li (3).\label{pi-li_1_for_L6}
\end{gather}\\
2) $3< u \leq x^{1+\gamma/\sqrt{\ln x}}$. Then
\[
\left|\frac{li(u) -li([u])}{\varphi (Q)} \right| \leq \int_{[u]}^{[u]+1}\frac{dt}{\ln t}\leq
\int_{2}^{3}\frac{dt}{\ln t} =li(3).
\]Hence,
\begin{align}
&\left| \pi(u; Q, W) - \frac{li (u)}{\varphi (Q)} \right|=
\left| \pi([u]; Q, W) - \frac{li (u)}{\varphi (Q)} - \frac{li ([u])}{\varphi (Q)} + \frac{li ([u])}{\varphi (Q)}  \right|\leq\notag\\
&\leq \left| \pi([u]; Q, W) - \frac{li ([u])}{\varphi (Q)} \right| + \left|\frac{li(u) -li([u])}{\varphi (Q)} \right|\leq
li(3)+\notag\\
&+ \left| \pi([u]; Q, W) - \frac{li ([u])}{\varphi (Q)} \right|.\label{pi-li_2_for_L6}
\end{align}From \eqref{pi-li_1_for_L6} and \eqref{pi-li_2_for_L6} we obtain
\begin{align}
&\max_{2 \leq u \leq x^{1+\gamma/\sqrt{\ln x}}} \max_{\substack{W\in \mathbb{Z}\notag\\ (W,Q)=1}}\left| \pi(u; Q, W) - \frac{li (u)}{\varphi (Q)} \right|\leq\notag\\
&\leq \max_{\substack{3 \leq u \leq x^{1+\gamma/\sqrt{\ln x}}\\ u\in \mathbb{Z}}} \max_{\substack{W\in \mathbb{Z}\\ (W,Q)=1}}\left| \pi(u; Q, W) - \frac{li (u)}{\varphi (Q)} \right| + 2li(3)+ 2.\label{f<discrete_f_for_L6}
\end{align}We have
\begin{equation}\label{x_1/2<exp_for_L6}
x^{1/2} \leq x\,\textup{exp}(-c_{3}\sqrt{\ln x}),
\end{equation}if $\widetilde{c}(\varepsilon, \delta)$ is chosen large enough. From \eqref{Discrete_Ineq_for_L6}, \eqref{f<discrete_f_for_L6} and \eqref{x_1/2<exp_for_L6} we obtain
\begin{align}
&\sum_{\substack{1 \leq Q \leq x^{1/2 - \delta}\\ (Q,B)=1}}\max_{2 \leq u \leq x^{1+\gamma/\sqrt{\ln x}}} \max_{\substack{W\in \mathbb{Z}\\ (W,Q)=1}}\left| \pi(u; Q, W) - \frac{li (u)}{\varphi (Q)} \right|\leq\notag\\
&\leq \sum_{\substack{1 \leq Q \leq x^{1/2 - \delta}\\ (Q,B)=1}}\max_{\substack{3 \leq u \leq x^{1+\gamma/\sqrt{\ln x}}\\ u\in \mathbb{Z}}} \max_{\substack{W\in \mathbb{Z}\\ (W,Q)=1}}\left| \pi(u; Q, W) - \frac{li (u)}{\varphi (Q)} \right|+\notag\\
&+ (2li(3)+2) x^{1/2}\leq (c_{4} + 2li(3)+2)x\,\textup{exp}(-c_{3}\sqrt{\ln x}).\label{Sum_Pi_for_L6}
\end{align}Thus, if $x$ is a real number with $x \geq \widetilde{c}(\varepsilon, \delta)$ and $q$ is an integer with $1 \leq q \leq (\ln x)^{1-\varepsilon}$, then there is a positive integer $B$ such that \eqref{B_for_L6} and \eqref{Sum_Pi_for_L6} hold. Let us denote $\widetilde{c}(\varepsilon, \delta)$ by $c(\varepsilon, \delta)$ and $c_{4} + 2li(3)+2$ by $c_2$. Lemma \ref{L6_Bombieri_Pi} is proved.

\section{Proof of Theorem \ref{T6} and Corollary \ref{C5}}\label{S_T6_C5}

Let us introduce some notation. Let $\mathcal{A}$ be a set of integers, $\mathcal{P}$ be a set of primes, $L(n)=l_{1} n + l_{2}$ be a linear function with coefficients in the integers. We define
\begin{gather*}
\mathcal{A}(x)=\{n\in \mathcal{A}:\ x\leq n< 2 x\},\\
\mathcal{A}(x; q, a)=\{n\in \mathcal{A}(x):\ n\equiv a \text{ (mod $q$)}\},\\
L(\mathcal{A})=\{L(n),\ n\in \mathcal{A}\},\\
\mathcal{P}_{L, \mathcal{A}}(x)= L(\mathcal{A}(x))\cap \mathcal{P},\\
\mathcal{P}_{L,\mathcal{A}}(x; q, a)= L(\mathcal{A}(x; q, a))\cap \mathcal{P},\\
\varphi_{L}(q)=\varphi(|l_1|q)/\varphi(|l_1|).
\end{gather*}

 Let $\mathcal{L}=\{L_1,\ldots, L_k\}$ be a set of distinct linear functions $L_i (n) = a_i n+ b_i$, $i=1,\ldots,k$, with coefficients in the positive integers. We say such a set is \emph{admissible} if for every prime $p$ there is an integer $n_p$ such that $(\prod_{i=1}^{k} L_{i}(n_p), p)=1$.

We focus on sets which satisfy the following hypothesis, which is given in terms of $(\mathcal{A}, \mathcal{L}, \mathcal{P}, B, x, \theta)$ for $\mathcal{L}$ an admissible set of linear functions, $B\in \mathbb{N}$, $x$ a large real number, and $0<\theta <1$.

\begin{hypothesis}\label{Hypothesis_1}
$(\mathcal{A}, \mathcal{L}, \mathcal{P}, B, x, \theta)$. Let $k=\# \mathcal{L}$.\par
\textup{(1)} $\mathcal{A}$ is well-distributed in arithmetic progressions: we have
\[
\sum_{1\leq q \leq x^{\theta}}\max_{a\in \mathbb{Z}}\Bigl|\#\mathcal{A}(x; q, a) - \frac{\# \mathcal{A}(x)}{q}\Bigr|\ll \frac{\#\mathcal{A}(x)}{(\ln x)^{100 k^{2}}}.
\]

\textup{(2)} Primes in $L(\mathcal{A})\cap \mathcal{P}$ are well-distributed in most arithmetic progressions: for any $L\in \mathcal{L}$ we have
\[
\sum_{\substack{1\leq q\leq x^{\theta}\\ (q, B)=1}} \max_{\substack{a\in \mathbb{Z}\\(L(a), q)=1}}\Bigl| \#\mathcal{P}_{L,\mathcal{A}}(x; q, a) - \frac{\#\mathcal{P}_{L,\mathcal{A}}(x)}{\varphi_{L}(q)}\Bigr|\ll \frac{\#\mathcal{P}_{L,\mathcal{A}}(x)}{(\ln x)^{100 k^{2}}}.
\]

\textup{(3)} $\mathcal{A}$ is not too concentrated in any arithmetic progression: for any $1\leq q< x^{\theta}$ we have
\[
\max_{a\in\mathbb{Z}}\#\mathcal{A}(x; q, a) \ll \frac{\#\mathcal{A}(x)}{q}.
\]
\end{hypothesis}
Maynard proved the following result (see \cite[Proposition 6.1]{Maynard}).

\begin{proposition}\label{Proposition_Maynard}
 Let $\alpha$ and $\theta$ be real numbers with $\alpha >0$ and $0< \theta <1$. Let $\mathcal{A}$ be a set of integers, $\mathcal{P}$ be a set of primes, $\mathcal{L}=\{L_{1},\ldots, L_{k}\}$ be an admissible set of $k$ linear functions, and let $B$, $x$ be integers. Let the coefficients $L_{i}(n)= a_{i}n+ b_{i}\in \mathcal{L}$ satisfy $1 \leq a_{i}, b_{i}\leq x^{\alpha}$ for all $1 \leq i \leq k$, and let $k \leq (\ln x)^{1/5}$ and $1 \leq B \leq x^{\alpha}$. Let $x^{\theta/10}\leq R \leq x^{\theta / 3}$. Let $\rho$, $\xi$ satisfy $k(\ln \ln x)^{2}/\ln x \leq \rho, \xi \leq \theta / 10$, and define
\[
\mathcal{S}(\xi; D) = \{n\in \mathbb{N}:\ p|n\Rightarrow\ (p> x^{\xi}\text{ or }p|D)\}.
\]Then there is a number $C>0$ depending only on $\alpha$ and $\theta$ such that the following holds. If $k \geq C$ and $(\mathcal{A}, \mathcal{L}, \mathcal{P}, B, x, \theta)$ satisfy Hypothesis \ref{Hypothesis_1}, then there is a choice of nonnegative weights $w_{n}=w_{n}(\mathcal{L})$ satisfying
\begin{equation}\label{Maynard1}
w_{n} \ll (\ln R)^{2 k}\prod_{i=1}^{k}\prod_{\substack{p|L_{i}(n)\\ p\nmid B}}4
\end{equation} such that the following statements hold.

\textup{(1)} We have
\begin{equation}\label{Maynard2}
\sum_{n\in \mathcal{A}(x)}w_{n} = \Bigl(1+ O\Bigl(\frac{1}{(\ln x)^{1/10}}\Bigr)\Bigr)\frac{B^{k}}{\varphi(B)^{k}}\mathfrak{S}_{B}(\mathcal{L})\#\mathcal{A}(x)(\ln R)^{k}I_{k}.
\end{equation}

\textup{(2)} For $L(n)=a_{L}n+ b_{L}\in \mathcal{L}$ we have
\begin{align}
&\sum_{n\in \mathcal{A}(x)}\textup{\textbf{1}}_{\mathcal{P}}(L(n))w_{n}\geq
\Bigl(1+ O\Bigl(\frac{1}{(\ln x)^{1/10}}\Bigr)\Bigr)\frac{B^{k-1}}{\varphi(B)^{k-1}}\mathfrak{S}_{B}(\mathcal{L})
\frac{\varphi (a_{L})}{a_{L}}\cdot\notag\\
&\cdot \#\mathcal{P}_{L, \mathcal{A}}(x)(\ln R)^{k+1}J_{k}+ O\Bigl(\frac{B^{k}}{\varphi(B)^{k}} \mathfrak{S}_{B}(\mathcal{L}) \#\mathcal{A}(x)(\ln R)^{k-1} I_{k}\Bigr).\label{Maynard3}
\end{align}

\textup{(3)} For $L(n)=a_{0}n+ b_{0}\notin \mathcal{L}$ and $D \leq x^{\alpha}$, if $\Delta_{L}\neq 0$ we have
\begin{align}
&\sum_{n\in \mathcal{A}(x)} \textup{\textbf{1}}_{\mathcal{S}(\xi; D)}(L(n))w_{n} \ll \xi^{-1} \frac{\Delta_{L}}{\varphi (\Delta_{L})}\frac{D}{\varphi(D)}\frac{B^{k}}{\varphi(B)^{k}}\mathfrak{S}_{B}(\mathcal{L})\#\mathcal{A}(x)\cdot\notag\\
&\cdot(\ln R)^{k-1}\label{Maynard4}
 I_{k},
\end{align}where
\[
\Delta_{L}=|a_{0}|\prod_{i=1}^{k}|a_{0}b_{i} - b_{0}a_{i}|.
\]

\textup{(4)} For $L\in \mathcal{L}$ we have
\begin{equation}\label{Maynard5}
\sum_{n\in \mathcal{A}(x)} \Bigl(\sum_{\substack{p|L(n)\\ p< x^{\rho}\\ p\nmid B}} 1\Bigr)w_{n}\ll \rho^{2} k^{4}(\ln k)^{2}
\frac{B^{k}}{\varphi(B)^{k}}\mathfrak{S}_{B}(\mathcal{L})\#\mathcal{A}(x) (\ln R)^{k}I_{k}.
\end{equation}

 Here $I_{k}$, $J_{k}$ are quantities depending only on $k$, and $\mathfrak{S}_{B}(\mathcal{L})$ is a quantity depending only on $\mathcal{L}$, and these satisfy
\begin{gather}
\mathfrak{S}_{B}(\mathcal{L}) = \prod_{p\nmid B}\Bigl(1 - \frac{\#\{1\leq n \leq p:\ p|\prod_{i=1}^{k} L_{i}(n)\}}{p}\Bigr)
\Bigl(1-\frac{1}{p}\Bigr)^{-k} \geq\notag\\
\geq \textup{exp}(-ck),\label{exp_Ineq_for_sigma_B_L}\\
I_{k}= \int_{0}^{\infty}\cdots\int_{0}^{\infty} F^{2}(t_{1},\ldots, t_{k})\,dt_{1}\ldots dt_{k}\gg (2k\ln k)^{-k},\label{Maynard_I_k_est}\\
J_{k}= \int_{0}^{\infty}\cdots\int_{0}^{\infty}\Bigl(\int_{0}^{\infty} F(t_{1},\ldots, t_{k})\, dt_{k} \Bigr)^{2}\,
dt_{1}\ldots dt_{k-1} \gg \frac{\ln k}{k} I_{k},\label{Maynard_Jk_est}
\end{gather}for a smooth function $F=F_{k}: \mathbb{R}^{k}\to \mathbb{R}$ depending only on $k$.

Here the implied constants depend only on $\alpha$, $\theta$, and the implied constants from Hypothesis \ref{Hypothesis_1}. The constant $c$ in the inequality \eqref{exp_Ineq_for_sigma_B_L} is a positive absolute constant.
\end{proposition}

\textsc{Proof of Theorem \ref{T6}.} First we prove the following

\begin{lemma}\label{L:admissible_set}
Let $k$ be a positive integer. Let $a$, $q$, $b_1,\ldots, b_k$ be positive integers with $b_{1} <\ldots < b_k$ and  $(a, q) =1$. Let
\[
L_{i}(n) = q n+ a+ q b_{i},\quad i=1,\ldots, k.
\]Then $\mathcal{L}=\{L_{1},\ldots, L_{k}\}$ is an admissible set iff for any prime $p$ such that $p\nmid q$ there is an integer $m_{p}$ such that $m_{p}\not\equiv b_{i}$ \textup{(mod $p$)} for all $1 \leq i\leq k$.
\end{lemma}
\textsc{Proof of Lemma \ref{L:admissible_set}.} 1) Let $\mathcal{L}=\{L_{1},\ldots, L_{k}\}$ be an admissible set. Let $p$ be a prime such that $p \nmid q$. Since $\mathcal{L}$ is an admissible set, there is an integer $n_{p}$ such that $(\prod_{i=1}^{k}L_{i}(n_{p}),p)=1$. Since $(q,p) = 1$, there is an integer $q^{\prime}$ such that $q q^{\prime} \equiv 1$ (mod $p$). We put $m_{p} = -(n_{p}+q^{\prime}a)$. Let $i$ be an integer with $1 \leq i \leq k$. Since  $(q^{\prime}, p)=1$ and $(L_{i}(n_{p}),p)=1$, we have $(q^{\prime}L_{i}(n_{p}), p)=1$. We have
\[
q^{\prime}L_{i}(n_{p}) \equiv -m_{p} + b_{i}\quad (\text{mod $p$}).
\]Hence, $m_{p} \not\equiv b_{i}$ (mod $p$).

2) Suppose that for any prime $p$ such that $p\nmid q$ there is an integer $m_{p}$ such that $m_{p}\not\equiv b_{i}$ \textup{(mod $p$)} for all $1 \leq i\leq k$. Let us show that then $\mathcal{L}$ is an admissible set. First we observe that $\mathcal{L}=\{L_{1},\ldots, L_{k}\}$ is the set of distinct linear functions $L_i (n) = q n+ l_i$, $i=1,\ldots,k$, with coefficients in the positive integers. Thus, we must prove that for any prime $p$ there is an integer $n_p$ such that $(\prod_{i=1}^{k} L_{i}(n_p), p)=1$. Let $p$ be a prime number. Let us consider two cases.\\
i) $p|q$. Since $(a,q)=1$, we have $(a,p)=1$. Let $i$ be an integer with $1 \leq i \leq k$. For any integer $n$ we have
\[
L_{i}(n)\equiv a\quad (\text{mod $p$}),
\]and hence $L_{i}(n) \not\equiv 0$ (mod $p$). Therefore $(\prod_{i=1}^{k} L_{i}(n), p)=1$. Therefore in this case we may take as $n_{p}$ any integer.\\
ii) $p\nmid q$. Then $(q,p)=1$, and hence there is an integer $c$ such that
\begin{equation}\label{qc_equiv_a_for_L}
qc\equiv a\quad (\text{mod $p$}).
\end{equation}By assumption, there is an integer $m_{p}$ such that $m_{p}\not\equiv b_{i}$ (mod $p$) for all $1 \leq i \leq k$. We put $n_{p}=-m_{p} - c$. Let $i$ be an integer with $1 \leq i \leq k$. We have
\[
n_{p}+c+ b_{i}\not\equiv 0\quad (\text{mod $p$}).
\]Since $(q,p)=1$, we obtain
\[
qn_{p}+ qc+ qb_{i}\not\equiv 0\quad (\text{mod $p$}).
\]Using \eqref{qc_equiv_a_for_L}, we obtain $L_{i}(n_{p})\not\equiv 0$ (mod $p$). Hence, $(L_{i}(n_{p}), p) =1$. Since this holds for all $1 \leq i \leq k$, we have $(\prod_{i=1}^{k}L_{i}(n_{p}), p) =1$. Lemma \ref{L:admissible_set} is proved.

The proof of the following lemma is based on ideas of Maynard (see the proof of Lemma 8.1 in \cite{Maynard}).

\begin{lemma}\label{L:sum_Delta_L}
There are positive absolute constants $c$ and $C$ such that the following holds. Let $x$ and $\eta$ be real numbers with $x\geq c$ and $(\ln x)^{-9/10}\leq \eta \leq 1$. Let $k$ and $a$ be positive integers. Let $b_{1},\ldots, b_{k}$ be integers with $1\leq b_{i} \leq \ln x$, $i=1,\ldots, k$.  Let $\mathcal{L}=\{L_{1},\ldots, L_{k}\}$ be the set of $k$ linear functions, where $$L_{i}(n)=an+b_{i},\quad i=1,\ldots, k.$$
For $L(n)=an+b$, $b\in \mathbb{Z}$, we define
\[
\Delta_{L} = a^{k+1}\prod_{i=1}^{k}|b_{i}-b|.
\]Then
\[
\sum_{\substack{1\leq b \leq \eta \ln x\\ L=an+b\notin \mathcal{L}}} \frac{\Delta_{L}}{\varphi(\Delta_{L})}\leq
C\ln\ln (a+2)\ln (k+1)\eta \ln x.
\]
\end{lemma}
\textsc{Proof of Lemma \ref{L:sum_Delta_L}.} Let us consider two cases.

1) Let $k> \ln\ln x$. We have
\[
\ln\ln x \geq 100,
\]if $c$ is chosen large enough. Therefore $k\geq 100$. Let $b$ to be an integer such that $1 \leq b \leq \eta\ln x$ and $L=an+b\notin \mathcal{L}$. Then $\Delta_{L}\in\mathbb{N}$. Applying Lemma \ref{L:low_est_Euler}, we have
\begin{equation}\label{Euler_Delta_1}
\frac{\Delta_{L}}{\varphi(\Delta_{L})}\leq c_{0}\ln\ln (\Delta_{L}+2),
\end{equation}where $c_{0}>0$ is an absolute constant. We have
\[
\ln \Delta_{L} = (k+1)\ln a + \sum_{i=1}^{k}\ln |b_{i} - b|.
\]For any $1 \leq i \leq k$ we have
\[
|b_{i} - b| \leq \ln x.
\] Hence,
\[
\ln \Delta_{L} \leq (k+1)\ln a  + k \ln\ln x \leq 2 k \ln a + k^{2}.
\]Since
\begin{gather*}
2k\ln a\leq k^{2}\ln (a+2),\\ k^{2}\leq k^{2}\ln (a+2),
\end{gather*}we have
\[
\ln \Delta_{L} \leq 2 k^{2}\ln(a+2).
\]We observe that if $u$ and $v$ are real numbers with $u\geq 2$ and $v\geq 2$, then
\begin{equation}\label{L:u+v<uv}
u+v\leq uv.
\end{equation} Applying \eqref{L:u+v<uv}, we obtain
\begin{align*}
&\ln (\Delta_{L}+2)\leq \ln(3\Delta_{L})=\ln \Delta_{L} + \ln 3\leq\\
&\leq 2 k^{2}\ln(a+2) + 3\leq 6k^{2}\ln(a+2).
\end{align*}Applying \eqref{L:u+v<uv} again, we have
\begin{align*}
&\ln\ln (\Delta_{L}+2)\leq \ln 6 + 2\ln k+\ln\ln(a+2)\leq\\
&\leq 2+2\ln k+25\ln\ln(a+2)\leq 4\ln k+ 25\ln\ln (a+2)\leq\\
&\leq 100 \ln k \ln\ln (a+2)\leq 100 \ln (k+1)\ln\ln (a+2).
\end{align*}Substituting this estimate into \eqref{Euler_Delta_1}, we obtain
\[
\frac{\Delta_{L}}{\varphi(\Delta_{L})} \leq 100 c_{0}\ln\ln (a+2) \ln (k+1)=
c_{1} \ln\ln (a+2) \ln (k+1),
\]where $c_{1} = 100 c_{0}>0$ is an absolute constant. We have
\begin{align}
\sum_{\substack{1\leq b \leq \eta \ln x\\ L=an+b\notin \mathcal{L}}} \frac{\Delta_{L}}{\varphi(\Delta_{L})}&\leq
c_{1}\ln\ln (a+2)\ln (k+1)\sum_{\substack{1\leq b \leq \eta \ln x\\ L=an+b\notin \mathcal{L}}}1\leq\notag\\
&\leq c_{1}\ln\ln (a+2)\ln (k+1) [\eta \ln x]\leq\notag\\
&\leq c_{1}\ln\ln (a+2)\ln (k+1) \eta \ln x.\label{L6.2_FINAL_I}
\end{align}

2) Let $1 \leq k \leq \ln\ln x$. For an integer $b$ we define
\[
\Delta(b):=\prod_{i=1}^{k}|b-b_{i}|.
\]Let $b$ to be an integer such that $1\leq b \leq \eta\ln x$ and $L=an+b\notin \mathcal{L}$. Applying Lemmas \ref{L:Euler_func_ineq} and \ref{L_about_Euler_func}, we obtain
\[
\frac{\Delta_{L}}{\varphi(\Delta_{L})} = \frac{a^{k+1}\Delta(b)}{\varphi(a^{k+1}\Delta(b))}\leq
\frac{a^{k+1}}{\varphi(a^{k+1})}
\frac{\Delta(b)}{\varphi(\Delta(b))}=
\frac{a}{\varphi(a)}
\frac{\Delta(b)}{\varphi(\Delta(b))}.
\]Hence,
\begin{equation}\label{L6.2_S_FINAL}
S=\sum_{\substack{1\leq b \leq \eta \ln x\\ L=an+b\notin \mathcal{L}}} \frac{\Delta_{L}}{\varphi(\Delta_{L})}\leq
\frac{a}{\varphi (a)}\sum_{\substack{1\leq b \leq \eta \ln x\\ L=an+b\notin \mathcal{L}}} \frac{\Delta(b)}{\varphi(\Delta(b))}
= \frac{a}{\varphi(a)}\widetilde{S}.
\end{equation}Applying Lemma \ref{L:Euler_func_representation}, we have
\begin{align}
&\widetilde{S}=\sum_{\substack{1\leq b \leq \eta \ln x\\ L=an+b\notin \mathcal{L}}} \frac{\Delta(b)}{\varphi(\Delta(b))}=
\sum_{\substack{1\leq b \leq \eta \ln x\\ L=an+b\notin \mathcal{L}}}\sum_{d|\Delta(b)}\frac{\mu^{2}(d)}{\varphi(d)}=\notag\\
&=\sum_{\substack{1\leq b \leq \eta \ln x\\ L=an+b\notin \mathcal{L}}}\Bigl(\sum_{\substack{1\leq d\leq \eta\ln x\\ d|\Delta(b)}}\frac{\mu^{2}(d)}{\varphi(d)}+ \sum_{\substack{d> \eta\ln x\\ d|\Delta(b)}}\frac{\mu^{2}(d)}{\varphi(d)}\Bigr)= S_{1}+S_{2}.\label{L6.2_S_wave_FINAL}
\end{align}

First we estimate the sum $S_{2}$. Let $b$ and $d$ be positive integers such that $1\leq b \leq \eta\ln x$, $L=an+b\notin \mathcal{L}$, $d> \eta \ln x$ and $d|\Delta(b)$. We claim that
\begin{equation}\label{L:Euler_Super_Ineq}
\frac{\mu^{2}(d)}{\varphi(d)}\leq \frac{\mu^{2}(d)\sum_{p|d}\ln p}{\varphi(d)\ln(\eta\ln x)}.
\end{equation}We have
\[
d>\eta \ln x\geq (\ln x)^{1/10}\geq 100,
\]if $c$ is chosen large enough. If $\mu^{2}(d)=0$, then the inequality \eqref{L:Euler_Super_Ineq} holds. Let $\mu^{2}(d)\neq 0$. Then $d$ is square-free. Therefore
\[
\sum_{p|d}\ln p = \ln d.
\]The inequality \eqref{L:Euler_Super_Ineq} is equivalent to the inequality
\[
\ln(\eta\ln x)\leq \sum_{p|d}\ln p = \ln d,
\]which, obviously, holds. Thus, \eqref{L:Euler_Super_Ineq} is proved. We have
\begin{align*}
&S_{2}= \sum_{\substack{1\leq b \leq \eta \ln x\\ L=an+b\notin \mathcal{L}}}\sum_{\substack{d> \eta\ln x\\ d|\Delta(b)}}\frac{\mu^{2}(d)}{\varphi(d)}\leq
\sum_{\substack{1\leq b \leq \eta \ln x\\ L=an+b\notin \mathcal{L}}}\sum_{\substack{d> \eta\ln x\\ d|\Delta(b)}}
\frac{\mu^{2}(d)\sum_{p|d}\ln p}{\varphi(d)\ln(\eta\ln x)}=\\
&=\sum_{\substack{1\leq b \leq \eta \ln x\\ L=an+b\notin \mathcal{L}}}
\sum_{p|\Delta(b)}\frac{\ln p}{\ln(\eta\ln x)}
\sum_{\substack{d> \eta\ln x\\ \text{$d$ is a multiple of $p$}\\ d|\Delta(b)}}\frac{\mu^{2}(d)}{\varphi(d)}.
\end{align*}Let $b\in \mathbb{N}$, $d\in \mathbb{N}$ and $p\in \mathbb{P}$ be such that $1 \leq b \leq \eta\ln x$, $L=an+b\notin \mathcal{L}$, $p|\Delta(b)$, $d>\eta\ln x$, $d$ is a multiple of $p$ and $d|\Delta (b)$. Then $d=pt$, where $t\in \mathbb{N}$, $t> (\eta\ln x)/p$ and $t|\Delta(b)$. We have (see Lemmas \ref{L:Euler_func_ineq} and \ref{L_about_Euler_func})
\[
\varphi(d)=\varphi(pt)\geq \varphi(p)\varphi(t)= (p-1)\varphi(t)\geq \frac{p}{2}\varphi(t).
\]Hence,
\[
\frac{\mu^{2}(d)}{\varphi(d)}= \frac{\mu^{2}(pt)}{\varphi(pt)}\leq
\frac{2 \mu^{2}(pt)}{p\varphi(t)}\leq \frac{2 \mu^{2}(t)}{p\varphi(t)}.
\]We obtain (see Lemma \ref{L:Euler_func_representation})
\[
\sum_{\substack{d> \eta\ln x\\\text{$d$ is a multiple of $p$}\\ d|\Delta(b)}}\frac{\mu^{2}(d)}{\varphi(d)}\leq
\frac{2}{p}\sum_{\substack{t> (\eta\ln x)/p\\ t|\Delta(b)}}\frac{\mu^{2}(t)}{\varphi(t)}\leq
\frac{2}{p}\sum_{t|\Delta(b)}\frac{\mu^{2}(t)}{\varphi(t)}=\frac{2}{p}\frac{\Delta(b)}{\varphi(\Delta(b))}.
\]Hence,
\begin{align*}
S_{2}&\leq \sum_{\substack{1\leq b \leq \eta \ln x\\ L=an+b\notin \mathcal{L}}}
\sum_{p|\Delta(b)}\frac{\ln p}{\ln(\eta\ln x)}\frac{2}{p}\frac{\Delta(b)}{\varphi(\Delta(b))}=\\
&=\frac{2}{\ln(\eta\ln x)}\sum_{\substack{1\leq b \leq \eta \ln x\\ L=an+b\notin \mathcal{L}}}
\frac{\Delta(b)}{\varphi(\Delta(b))}\sum_{p|\Delta(b)}\frac{\ln p}{p}.
\end{align*}Since
\[
\eta\geq (\ln x)^{-9/10},
\]we have
\[
\frac{2}{\ln(\eta\ln x)}\leq \frac{2}{\ln((\ln x)^{1/10})}=\frac{20}{\ln\ln x}.
\]We obtain
\begin{equation}\label{L6.2_S2_INEQ}
S_{2}\leq \frac{20}{\ln\ln x}\sum_{\substack{1\leq b \leq \eta \ln x\\ L=an+b\notin \mathcal{L}}}
\frac{\Delta(b)}{\varphi(\Delta(b))}\sum_{p|\Delta(b)}\frac{\ln p}{p}.
\end{equation}Let $b$ to be an integer such that $1 \leq b \leq \eta\ln x$ and $L=an+b\notin\mathcal{L}$. Applying Lemmas \ref{L:low_est_Euler} and \ref{Lemma_sum_ln_p}, we have
\begin{gather}
\frac{\Delta(b)}{\varphi(\Delta(b))}\leq c_{2}\ln\ln (\Delta(b)+2)\leq
c_{2}\ln\ln (3\Delta(b)),\label{L6.2_Delta(b)_Ineq}\\
\sum_{p|\Delta(b)}\frac{\ln p}{p}\leq c_{3}\ln\ln (3\Delta(b)),\label{L6.2_sum_ln_p}
\end{gather}where $c_{2}>0$ and $c_{3}>0$ are absolute constants. We have
\[
\ln \Delta(b)=\sum_{i=1}^{k}\ln |b_{i}-b|\leq k \ln\ln x\leq (\ln\ln x)^{2}.
\]Hence,
\begin{align}
&\ln\ln(3\Delta(b))=\ln(\ln 3 + \ln \Delta(b))\leq\notag\\
&\leq \ln(\ln 3 + (\ln\ln x)^{2})\leq 3\ln\ln\ln x,\label{L6.2_lnln_Delta(b)_EST}
\end{align}if $c$ is chosen large enough. From \eqref{L6.2_Delta(b)_Ineq}, \eqref{L6.2_sum_ln_p} and
\eqref{L6.2_lnln_Delta(b)_EST} we obtain
\[
\frac{\Delta(b)}{\varphi(\Delta(b))}\sum_{p|\Delta(b)}\frac{\ln p}{p}\leq c_{2}c_{3}9 (\ln\ln\ln x)^{2}=
c_{4}(\ln\ln\ln x)^{2},
\]where $c_{4}=9c_{2}c_{3}>0$ is an absolute constant. Substituting this estimate into \eqref{L6.2_S2_INEQ}, we obtain \[
S_{2} \leq \frac{20 c_{4} (\ln\ln\ln x)^{2}}{\ln\ln x} \eta\ln x.
\]We have
\[
\frac{20 c_{4} (\ln\ln\ln x)^{2}}{\ln\ln x} \leq 1,
\]if $c$ is chosen large enough.  Hence,
\begin{equation}\label{L6.2_S2_FINAL}
S_{2} \leq \eta\ln x\leq \frac{1}{\ln 2}\ln (k+1) \eta \ln x\leq 2 \ln (k+1) \eta \ln x.
\end{equation}

 Now we estimate $S_{1}$. We have
\begin{align}
&S_{1} = \sum_{\substack{1\leq b \leq \eta\ln x\\ L=an+b\notin\mathcal{L}}}\sum_{\substack{1\leq d\leq\eta\ln x\\
d|\Delta(b)}}\frac{\mu^{2}(d)}{\varphi(d)}=
\sum_{1\leq d\leq\eta\ln x}
\sum_{\substack{1\leq b \leq \eta\ln x\\ L=an+b\notin\mathcal{L}\\d|\Delta(b)}}\frac{\mu^{2}(d)}{\varphi(d)}=\notag\\
&=\sum_{1\leq d\leq\eta\ln x}\frac{\mu^{2}(d)}{\varphi(d)}
\sum_{\substack{1\leq b \leq \eta\ln x\\ L=an+b\notin\mathcal{L}\\d|\Delta(b)}}1=
\sum_{1\leq d\leq\eta\ln x}\frac{\mu^{2}(d)}{\varphi(d)}N_{0}(d)=\notag\\
&=\sum_{\substack{1\leq d\leq\eta\ln x\\d\in \mathcal{M}}}\frac{1}{\varphi(d)}N_{0}(d).\label{L6.2_S1_UP}
\end{align}

Let $d$ be an integer such that $1\leq d \leq \eta\ln x$ and $d\in \mathcal{M}$. We claim that
\begin{equation}\label{L6.2_N0(d)_INEQ}
N_{0}(d)\leq \frac{2\eta\ln x}{d}\prod_{p|d}\min(p,k).
\end{equation}If $d=1$, then the inequality is obvious. Let $d>1$. We define
\[
R(b)=(b-b_{1})\cdots(b-b_{k}).
\]Then $\Delta(b)=|R(b)|$. We have
\[
N_{0}(d)=\sum_{\substack{1\leq b \leq \eta\ln x\\ L=an+b\notin\mathcal{L}\\d|\Delta(b)}}1=
\sum_{\substack{1\leq b \leq \eta\ln x\\ L=an+b\notin\mathcal{L}\\R(b)\equiv 0\, (\text{mod $d$})}}1.
\]Let $d$ be expressed in standard form
\[
d=q_{1}\cdots q_{r},
\]where $q_{1}<\ldots <q_{r}$ are prime numbers. It is well-known (see, for example, \cite[Chapter 4]{Vinogradov}), the congruence
\[
R(b)\equiv 0\ (\text{mod $d$})
\]is equivalent to the system of congruences
\begin{equation}\label{L6.2_II_SYSTEM}
\begin{cases}
R(b)&\equiv 0\ (\text{mod $q_{1}$}),\\
&\ \vdots\\
R(b)&\equiv 0\ (\text{mod $q_{r}$}).
\end{cases}
\end{equation}Let $1\leq j \leq r$. Let $\Omega_{j}$ be the set of numbers of a complete system of residues modulo $q_{j}$ satisfying the congruence
\[
R(b)\equiv 0\ (\text{mod $q_{j}$}).
\]Since $R(b_{1})=0$, we see that $\Omega_{j}\neq\varnothing$. Since the leading coefficient of the polynomial $R(b)$ is $1$ and the degree of the polynomial $R(b)$ is $k$, we have $\#\Omega_{j}\leq k$ (see, for example, \cite[Chapter 4]{Vinogradov})). It is clear that $\#\Omega_{j} \leq q_{j}$. We obtain
\[
\#\Omega_{j} \leq \min (q_{j},k).
\]The system \eqref{L6.2_II_SYSTEM} is equivalent to the union of
\[
T=\#\Omega_{1}\cdots\#\Omega_{r}
\]systems
\begin{equation}\label{L6.2_SYSTEM_III}
\begin{cases}
b&\equiv \tau_{1}\ (\text{mod $q_{1}$}),\\
&\ \vdots\\
b&\equiv \tau_{r}\ (\text{mod $q_{r}$}),
\end{cases}
\end{equation}where $\tau_{1}\in\Omega_{1},\ldots, \tau_{r}\in\Omega_{r}$. It is well-known (see, for example, \cite[Chapter 4]{Vinogradov}), the system of congruences \eqref{L6.2_SYSTEM_III} is equivalent to the congruence
\[
b\equiv x_{0}\ (\text{mod $d$}),
\]where $x_{0}=x_{0}(\tau_{1},\ldots, \tau_{r})$. It is also well-known that the numbers $x_{0}(\tau_{1},\ldots, \tau_{r})$, $\tau_{1}\in \Omega_{1}, \ldots, \tau_{r}\in \Omega_{r},$ are incongruent modulo $d$. Thus,
\[
\{b\in \mathbb{Z}:\ R(b)\equiv 0\ (\text{mod $d$})\}=\bigsqcup_{\tau_{1}\in\Omega_{1},\ldots, \tau_{r}\in \Omega_{r}}
\{x_{0}(\tau_{1},\ldots,\tau_{r})+ dt,\ t\in\mathbb{Z}\}.
\]Let $\tau_{1}\in \Omega_{1},\ldots, \tau_{r}\in\Omega_{r}, x_{0}=x_{0}(\tau_{1},\ldots, \tau_{r})$. We have
\begin{align*}
&\#\{t\in \mathbb{Z}:\ 1\leq x_{0}+dt\leq \eta\ln x\}=\left[\frac{\eta\ln x-x_{0}}{d}\right]-
\left\lceil\frac{1-x_{0}}{d}\right\rceil +1=\\
&=\frac{\eta\ln x-x_{0}}{d} - \theta_{1}- \left(\frac{1-x_{0}}{d}+\theta_{2}\right)+1=
\frac{\eta\ln x}{d} + 1 -\theta_{1}-\theta_{2} - \frac{1}{d},
\end{align*}where $\theta_{1}$ and $\theta_{2}$ are real numbers with $0 \leq \theta_{1} <1$ and $0 \leq \theta_{2} <1$. Since $1 \leq d \leq \eta\ln x$, we obtain
\[
\#\{t\in \mathbb{Z}:\ 1\leq x_{0}+dt\leq \eta\ln x\}\leq
\frac{\eta\ln x}{d} + 1\leq 2\frac{\eta\ln x}{d}.
\]We obtain
\[
N_{0}(d)\leq \frac{2\eta\ln x}{d}T\leq \frac{2\eta\ln x}{d}\prod_{p|d}\min(p,k).
\]The inequality \eqref{L6.2_N0(d)_INEQ} is proved.

 Substituting \eqref{L6.2_N0(d)_INEQ} into \eqref{L6.2_S1_UP}, we obtain
\begin{align}
S_{1}&\leq \sum_{\substack{1\leq d\leq\eta\ln x\\d\in \mathcal{M}}}
\frac{1}{\varphi(d)}\frac{2\eta\ln x}{d}\prod_{p|d}\min(p,k)=\notag\\
&=2\eta\ln x\sum_{\substack{1\leq d\leq\eta\ln x\\d\in \mathcal{M}}}
\frac{\prod_{p|d}\min(p,k)}{d\varphi(d)}= 2\eta\ln x S_{3}.\label{L6.2_S1_S3_INEQUA}
\end{align}Let $d$ be an integer with $1\leq d \leq \eta\ln x$ and $d\in \mathcal{M}$. We have (see Lemmas \ref{L:Euler_func_multipl} and \ref{L_about_Euler_func})
\begin{gather*}
d=\prod_{p|d}p,\\
\varphi(d)=\prod_{p|d}\varphi(p)=\prod_{p|d}(p-1),\\
\frac{\prod_{p|d}\min(p,k)}{d\varphi(d)}=
\frac{\prod_{p|d}\min(p,k)}{\prod_{p|d}p(p-1)}=\prod_{\substack{p|d\\ p\leq k}}\frac{1}{p-1}
\prod_{\substack{p|d\\ p> k}}\frac{k}{p(p-1)}.
\end{gather*}Hence,
\begin{align}
&S_{3}=\sum_{\substack{1\leq d\leq\eta\ln x\\d\in \mathcal{M}}}
\prod_{\substack{p|d\\ p\leq k}}\frac{1}{p-1}
\prod_{\substack{p|d\\ p> k}}\frac{k}{p(p-1)}\leq\notag\\
&\leq \prod_{p\leq k} \Bigl(1+\frac{1}{p-1}\Bigr)\prod_{p>k}\Bigl(1+\frac{k}{p(p-1)}\Bigr)=AB.\label{L6.2_S3_AB}
\end{align}

We have (see Lemma \ref{L_prod(1-1/p)})
\begin{align}
A&=\prod_{p\leq k}\Bigl(1+\frac{1}{p-1}\Bigr)\leq
\prod_{p\leq k+1}\Bigl(1+\frac{1}{p-1}\Bigr)=\notag\\
 &=\prod_{p\leq k+1}\Bigl(1-\frac{1}{p}\Bigr)^{-1}
\leq c_{5}\ln (k+1),\label{L6.2_A}
\end{align}where $c_{5}>0$ is an absolute constant.

 Now we estimate $B$. Since $\ln(1+u)\leq u$, $u> 0$, we have
 \begin{align*}
\ln B&= \sum_{p> k} \ln\Bigl(1+\frac{k}{p(p-1)}\Bigr)\leq
\sum_{p> k}\frac{k}{p(p-1)}=\\
&=k\sum_{p\geq k+1}\frac{1}{p(p-1)}\leq
k\sum_{n\geq k+1}\frac{1}{n(n-1)}.
\end{align*}We define
\begin{align*}
s_{m}&=\sum_{n=k+1}^{m}\frac{1}{n(n-1)}=\sum_{n=k+1}^{m}\Bigl(\frac{1}{n-1}-\frac{1}{n}\Bigr)=\\
&=\frac{1}{k} - \frac{1}{m},\ m\geq k+1.
\end{align*}Hence,
\[
\sum_{n\geq k+1}\frac{1}{n(n-1)}=\lim_{m\to +\infty} s_{m}=\frac{1}{k}.
\]We obtain $\ln B\leq 1$, i.\,e.
\begin{equation}\label{L6.2_B}
B\leq e<3.
\end{equation} We have (see \eqref{L6.2_S3_AB}, \eqref{L6.2_A} and \eqref{L6.2_B})
\[
S_{3} \leq c_{6}\ln (k+1),
\]where $c_{6}>0$ is an absolute constant. Substituting this estimate into \eqref{L6.2_S1_S3_INEQUA}, we obtain
\begin{equation}\label{L6.2_S1_FINAL}
S_{1} \leq c_{7}\ln (k+1)\eta\ln x,
\end{equation}where $c_{7}>0$ is an absolute constant.

 We obtain (see \eqref{L6.2_S_wave_FINAL}, \eqref{L6.2_S2_FINAL} and \eqref{L6.2_S1_FINAL})
\[
\widetilde{S} \leq (c_{7}+2)\ln (k+1)\eta\ln x = c_{8}\ln (k+1)\eta\ln x,
\]where $c_{8}=c_{7}+2>0$ is an absolute constant. We obtain (see \eqref{L6.2_S_FINAL} and Lemma \ref{L:low_est_Euler})
\begin{equation}\label{L6.2_FINAL_II}
S\leq c_{8}\frac{a}{\varphi(a)}\ln (k+1)\eta\ln x\leq c_{9}\ln\ln(a+2)\ln (k+1)\eta\ln x,
\end{equation}where $c_{9}>0$ is an absolute constant. We put
\[
C=c_{1}+c_{9},
\]where $c_{1}$ is the constant in \eqref{L6.2_FINAL_I}. Then $C>0$ is an absolute constant and in both cases, $1\leq k \leq \ln\ln x$ and $k>\ln\ln x$, we have
\[
\sum_{\substack{1\leq b \leq \eta \ln x\\ L=an+b\notin \mathcal{L}}} \frac{\Delta_{L}}{\varphi(\Delta_{L})}
\leq C \ln\ln(a+2)\ln (k+1)\eta\ln x.
\]Lemma \ref{L:sum_Delta_L} is proved.

\begin{lemma}\label{L_wn_for_T6}
Let  $\mathcal{A}=\mathbb{N}$, $\mathcal{P}=\mathbb{P}$, $\alpha = 1/5$, $\theta = 1/3$, let $C_{0}=C(1/5, 1/3)>0$ be the absolute constant in Proposition \textup{\ref{Proposition_Maynard}}. Let $\varepsilon$ be a real number with $0< \varepsilon < 1$. Then there is a number $c_{0}(\varepsilon)>0$ such that the following holds. Let $x\in \mathbb{N}$, $y\in\mathbb{R}$, $q\in\mathbb{N}$ be such that $x\geq c_{0}(\varepsilon)$, $1 \leq y \leq \ln x$, $1 \leq q \leq y^{1-\varepsilon}$. Then there is a positive integer $B$ such that
\begin{equation}\label{B_for_T6}
1\leq B \leq \textup{exp}(\vartheta\sqrt{\ln x}),\quad 1\leq\frac{B}{\varphi(B)}\leq 2,\quad (B,q)=1.
\end{equation}Furthermore, let
\begin{gather}
k\in \mathbb{N},\ \rho\in \mathbb{R},\ \xi\in \mathbb{R},\ R\in\mathbb{R},\ \eta\in\mathbb{R},\ a\in\mathbb{Z},  \notag\\
C_{0}\leq k \leq (\ln x)^{1/5},\label{L_w_n_for_T6_1}\\
\frac{k(\ln\ln x)^{2}}{\ln x} \leq \rho \leq \frac{1}{30},\quad \xi = \rho,\label{L_w_n_for_T6_2}\\
R=x^{1/9},\quad 0<\eta \leq \frac{1}{2},\label{L_w_n_for_T6_3}\\
1\leq a \leq q,\quad (a,q)=1.\label{L_w_n_for_T6_4}
\end{gather}Let $\mathcal{L}=\{L_{1},\ldots, L_{k}\}$ be an admissible set of $k$ linear functions, where $L_{i}(n) = q n + a+ q b_{i}$, $i=1,\ldots, k$, $b_1, \ldots, b_k$ are positive integers with $b_{1}<\ldots < b_{k}$ and $q b_{k} \leq \eta y$. Then the assumption of Proposition \ref{Proposition_Maynard} holds and there are nonnegative weights $w_{n}= w_{n}(\mathcal{L})$ satisfying the statement of Proposition \ref{Proposition_Maynard}; the implied constants in \eqref{Maynard1} -- \eqref{Maynard5} are positive and absolute.

Here $\vartheta>0$ is an absolute constant.
\end{lemma}
\textsc{Proof of Lemma \ref{L_wn_for_T6}.} We choose $c_{0}(\varepsilon)$ later; this number is large enough. We take $\delta =0.1$ and let $c_{0}(\varepsilon)\geq c(\varepsilon, \delta)=c(\varepsilon, 0.1)$, where  $c(\varepsilon, \delta)$ is the quantity in Lemma \ref{L6_Bombieri_Pi}. Let $x\in \mathbb{N}$, $y\in\mathbb{R}$, $q\in\mathbb{N}$ be such that $x\geq c_{0}(\varepsilon)$, $1 \leq y \leq \ln x$, $1 \leq q \leq y^{1-\varepsilon}$. By Lemma \ref{L6_Bombieri_Pi}, there is a positive integer $B$ such that
\[
1\leq B \leq \textup{exp}(c_{1}\sqrt{\ln x}),\quad 1\leq\frac{B}{\varphi(B)}\leq 2,\quad (B,q)=1
\] and
\begin{equation}\label{Sum_Pi_for_T6_Lemma}
\sum_{\substack{1\leq Q\leq x^{0.4}\\
(Q,B)=1}} \max_{2 \leq u \leq x^{1+ \gamma/\sqrt{\ln x}}}\max_{\substack{W\in \mathbb{Z}:\\ (W,Q)=1}}
\Bigl|\pi (u; Q, W)-\frac{li(u)}{\varphi (Q)}\Bigr| \leq c_2 x\,\textup{exp}(-c_3\sqrt{\ln x}),
\end{equation}where $c_1$, $\gamma$, $c_2$ and $c_3$ are positive absolute constants. Let \eqref{L_w_n_for_T6_1} -- \eqref{L_w_n_for_T6_4} hold. Let $\mathcal{L}=\{L_{1},\ldots, L_{k}\}$ be an admissible set of $k$ linear functions, where $L_{i}(n) = q n + a+ q b_{i}$, $i=1,\ldots, k$, $b_1, \ldots, b_k$ are positive integers with $b_{1}<\ldots < b_{k}$ and $q b_{k} \leq \eta y$. Let us show that the assumption of Proposition \ref{Proposition_Maynard} holds. First we show that $(\mathcal{A}, \mathcal{L}, \mathcal{P}, B, x, 1/3)$ satisfy Hypothesis \ref{Hypothesis_1}.

I) Let us show that part (2) of Hypothesis \ref{Hypothesis_1} holds. Let $L(n)=l_{1} n+ l_{2}\in \mathcal{L}$. It is clear that
\begin{equation}\label{l_1<ln_x}
1\leq l_{1} \leq \ln x,\qquad 1\leq l_{2} \leq \ln x.
\end{equation}Let us show that
\begin{equation}\label{Basic_I_for_T6}
S:=\sum_{\substack{1 \leq r \leq x^{1/3}\\ (r, B)=1}}\max_{\substack{b\in \mathbb{Z}\\ (L(b), r)=1}}
\Bigl|\#\mathcal{P}_{L, \mathcal{A}}(x; r, b) - \frac{\#\mathcal{P}_{L, \mathcal{A}}(x)}{\varphi_{L}(r)}\Bigr|
\leq \frac{\#\mathcal{P}_{L, \mathcal{A}}(x)}{(\ln x)^{100 k^{2}}}.
\end{equation}It is not hard to see that
\begin{gather*}
\mathcal{P}_{L, \mathcal{A}}(x)=\bigl\{ l_{1}x+ l_{2} \leq p < 2l_{1}x+ l_{2}:\ p\equiv l_{2}\ \text{(mod $l_{1}$)}\bigr\},\\
\mathcal{P}_{L, \mathcal{A}}(x; r, b)=\bigl\{l_{1}x + l_{2}\leq p < 2l_{1}x+l_{2}:\ p\equiv l_{1}b+l_{2}\ \text{(mod $l_{1}r$)}\bigr\},
\end{gather*}and hence
\begin{gather}
\#\mathcal{P}_{L, \mathcal{A}}(x)= \pi(2l_{1}x+l_{2}-1; l_1, l_2) -
\pi(l_{1}x+l_{2}-1; l_1, l_2),\label{Number_P_L_A_for_T6}\\
\#\mathcal{P}_{L, \mathcal{A}}(x; r, b)=\pi(2l_{1}x+l_{2}-1; l_{1}r, L(b))-
\pi(l_{1}x+l_{2}-1; l_{1}r, L(b)).\notag
\end{gather}We obtain
\begin{gather}
S=\sum_{\substack{1 \leq r \leq x^{1/3}\\ (r, B)=1}}\max_{\substack{b\in \mathbb{Z}\\ (L(b), r)=1}}
\Bigl|\pi(2l_{1}x+l_{2}-1; l_{1}r, L(b))-\notag\\
-\pi(l_{1}x+l_{2}-1; l_{1}r, L(b)) -\notag\\
- \frac{\pi(2l_{1}x+l_{2}-1; l_1, l_2) -
\pi(l_{1}x+l_{2}-1; l_1, l_2)}{\varphi(l_{1}r)/\varphi(l_{1})}\Bigr|\leq\notag\\
\leq S_{1}+S_{2}+S_{3}+S_{4},\label{S_view_for_T6}
\end{gather}where
\begin{align*}
S_{1} &= \sum_{\substack{1 \leq r \leq x^{1/3}\\ (r, B)=1}}\max_{\substack{b\in \mathbb{Z}\\ (L(b), r)=1}}
\Bigl| \pi(l_{1}x+l_{2}-1; l_{1}r, L(b)) - \frac{li(l_{1}x+l_{2}-1)}{\varphi(l_{1}r)}\Bigr|,\\
S_{2} &= \sum_{\substack{1 \leq r \leq x^{1/3}\\ (r, B)=1}}
\Bigl|\frac{\pi(l_{1}x+l_{2}-1; l_{1}, l_{2})}{\varphi(l_{1}r)/\varphi(l_{1})} -
\frac{li(l_{1}x+l_{2}-1)}{\varphi(l_{1}r)} \Bigr|,\\
S_{3} &= \sum_{\substack{1 \leq r \leq x^{1/3}\\ (r, B)=1}}\max_{\substack{b\in \mathbb{Z}\\ (L(b), r)=1}}
\Bigl| \pi(2l_{1}x+l_{2}-1; l_{1}r, L(b)) - \frac{li(2l_{1}x+l_{2}-1)}{\varphi(l_{1}r)}\Bigr|,\\
S_{4} &= \sum_{\substack{1 \leq r \leq x^{1/3}\\ (r, B)=1}}
\Bigl|\frac{\pi(2l_{1}x+l_{2}-1; l_{1}, l_{2})}{\varphi(l_{1}r)/\varphi(l_{1})} -
\frac{li(2l_{1}x+l_{2}-1)}{\varphi(l_{1}r)} \Bigr|.
\end{align*}

Let us show that
\begin{equation}\label{L(b)_coprime_l1}
(L(b), l_{1})=1
\end{equation}for any $b\in \mathbb{Z}$. Assume the converse: there is an integer $b$ such that $(L(b), l_{1})>1$. Then there is a prime $p$ such that $p|l_{1}$ and $p|L(b)$. Hence $p|l_{2}$, and we see that $p|L(n)$ for any integer $n$. Since $L\in \mathcal{L}$, we see that $p|L_{1}(n)\cdots L_{k}(n)$ for any integer $n$. But this contradicts the fact that $\mathcal{L}=\{L_{1},\ldots, L_{k}\}$ is an admissible set. Thus, \eqref{L(b)_coprime_l1} is proved.

 We observe that since $(B,q)=1$ and $l_{1}=q$, we have
\begin{equation}\label{B_coprime_l1}
(B,l_{1})=1.
\end{equation}Let $r$ be an integer with $1\leq r \leq x^{1/3}$ and $(r,B)=1$. Applying \eqref{l_1<ln_x}, we have \begin{gather*}
l_{1} r\leq  x^{1/3} \ln x \leq x^{0.4},\\
l_{1}x+l_{2}-1\geq l_{1}x\geq x\geq 2,\\
l_{1}x+ l_{2} -1 \leq 2 x\ln x\leq x^{1+\gamma/\sqrt{\ln x}},
\end{gather*}if $c_{0}(\varepsilon)$ is chosen large enough. Hence, we obtain (see \eqref{L(b)_coprime_l1}, \eqref{B_coprime_l1} and \eqref{Sum_Pi_for_T6_Lemma})
\begin{align}
&S_{1} = \sum_{\substack{r:\\ l_{1} \leq l_{1}r \leq l_{1}x^{1/3}\\ (l_{1}r, B)=1}}\max_{\substack{b\in \mathbb{Z}\\ (L(b), l_{1}r)=1}}
\Bigl| \pi(l_{1}x+l_{2}-1; l_{1}r, L(b)) - \frac{li(l_{1}x+l_{2}-1)}{\varphi(l_{1}r)}\Bigr|\leq\notag\\
&\leq \sum_{\substack{1\leq Q\leq x^{0.4}\\
(Q,B)=1}} \max_{2 \leq u \leq x^{1+ \gamma/\sqrt{\ln x}}}\max_{\substack{W\in \mathbb{Z}:\\ (W,Q)=1}}
\Bigl|\pi (u; Q, W)-\frac{li(u)}{\varphi (Q)}\Bigr| \leq c_2 x\,\textup{exp}(-c_3\sqrt{\ln x}).\label{S1_for_T6_est}
\end{align} Applying Lemmas \ref{L:Euler_func_ineq} and \ref{L:series_1/Euler_func}, we obtain
\begin{align*}
S_{2}&=\varphi(l_{1})\Bigl|\pi(l_{1}x+l_{2}-1; l_{1}, l_{2})- \frac{li(l_{1}x+l_{2}-1)}{\varphi(l_{1})}\Bigr|
\sum_{\substack{1 \leq r \leq x^{1/3}\\ (r, B)=1}}\frac{1}{\varphi (l_{1}r)}\leq\\
&\leq \Bigl|\pi(l_{1}x+l_{2}-1; l_{1}, l_{2})- \frac{li(l_{1}x+l_{2}-1)}{\varphi(l_{1})}\Bigr|
\sum_{1 \leq r \leq x^{1/3}}\frac{1}{\varphi (r)}\leq\\
&\leq \widetilde{c}\ln x \Bigl|\pi(l_{1}x+l_{2}-1; l_{1}, l_{2})- \frac{li(l_{1}x+l_{2}-1)}{\varphi(l_{1})}\Bigr|,
\end{align*}where $\widetilde{c}>0$ is an absolute constant. Since (see \eqref{l_1<ln_x})
\[
l_{1} x+ l_{2}-1\geq l_{1}x\geq x,
\] we obtain
\[
1\leq l_{1} \leq \ln x \leq \ln (l_{1}x+l_{2} -1).
\]Hence, (see, for example, \cite[Chapter 22]{Davenport})
\begin{gather}
\Bigl|\pi (l_{1}x+l_{2}-1; l_{1}, l_{2}) - \frac{li(l_{1}x+l_{2}-1)}{\varphi (l_{1})}\Bigr| \leq\notag\\
\leq C(l_{1}x+l_{2}-1)\textup{exp}\bigl(-c\sqrt{\ln(l_{1}x+l_{2}-1)}\bigr),\label{Pi_formula_for_T6}
\end{gather}where $C$ and $c$ are positive absolute constants. We have
\begin{gather}
\textup{exp}\bigl(-c\sqrt{\ln(l_{1}x+l_{2}-1)}\bigr) \leq \textup{exp}(-c\sqrt{\ln x}),\label{Formula_1_for_T6}\\
l_{1}x+l_{2}-1\leq x\ln x + \ln x \leq 2 x\ln x.\label{Formula_2_for_T6}
\end{gather}We have
\begin{equation}\label{Formula_3_for_T6}
-c\sqrt{\ln x} + 2\ln\ln x \leq -\frac{c}{2}\sqrt{\ln x},
\end{equation}if $c_{0}(\varepsilon)$ is chosen large enough. Hence,
\begin{align}
S_{2} &\leq \widetilde{C} x\,\textup{exp}(-c\sqrt{\ln x} + 2\ln\ln x)\leq \notag\\
&\leq \widetilde{C} x\,\textup{exp}(-(c/2)\sqrt{\ln x}),\label{S2_for_T6_est}
\end{align}where $\widetilde{C}=2\widetilde{c}C$ is a positive absolute constant. Similarly, it can be shown that
\begin{equation}\label{S3_for_T6_est}
S_{3} \leq C x\,\textup{exp}(-c\sqrt{\ln x}),\qquad S_{4} \leq C x\,\textup{exp}(-c\sqrt{\ln x}),
\end{equation}where $C$ and $c$ are positive absolute constants. Substituting \eqref{S1_for_T6_est}, \eqref{S2_for_T6_est} and \eqref{S3_for_T6_est} into \eqref{S_view_for_T6}, we obtain
\begin{equation}\label{Basic_2_1_for_T6}
\sum_{\substack{1 \leq r \leq x^{1/3}\\ (r, B)=1}}\max_{\substack{b\in \mathbb{Z}\\ (L(b), r)=1}}
\Bigl|\#\mathcal{P}_{L, \mathcal{A}}(x; r, b) - \frac{\#\mathcal{P}_{L, \mathcal{A}}(x)}{\varphi_{L}(r)}\Bigr|
\leq c_{4} x\,\textup{exp}(-c_{5}\sqrt{\ln x}),
\end{equation}where $c_{4}$ and $c_{5}$ are positive absolute constants.

 Applying \eqref{Pi_formula_for_T6} -- \eqref{Formula_3_for_T6}, we have
\begin{gather*}
\pi (l_{1}x+l_{2}-1; l_{1}, l_{2}) = \frac{li(l_{1}x+l_{2}-1)}{\varphi (l_{1})} + R_{1},\\
|R_{1}| \leq C x\,\textup{exp}(-c\sqrt{\ln x}),
\end{gather*}where $C$ and $c$ are positive absolute constants. Similarly, it can be shown that
\begin{gather*}
\pi (2l_{1}x+l_{2}-1; l_{1}, l_{2}) = \frac{li(2l_{1}x+l_{2}-1)}{\varphi (l_{1})} + R_{2},\\
|R_{2}| \leq C x\,\textup{exp}(-c\sqrt{\ln x}),
\end{gather*}where $C$ and $c$ are positive absolute constants. Hence, we have (see \eqref{Number_P_L_A_for_T6})
\begin{gather}
\#\mathcal{P}_{L, \mathcal{A}}(x) = \frac{li(2l_{1}x+l_{2}-1) - li(l_{1}x+l_{2}-1)}{\varphi (l_{1})}+ R, \label{P_L_A_est_for_T6} \\
|R|\leq c_{6} x\,\textup{exp}(-c_{7}\sqrt{\ln x}),\label{|R|<_for_T6}
\end{gather}where $c_{6}$ and $c_{7}$ are positive absolute constants. We have
\begin{gather*}
2l_{1}x+l_{2} - 1\leq 2x\ln x + \ln x \leq 3 x \ln x,\\
\ln (2l_{1} x+ l_{2} -1) \leq \ln x + \ln\ln x + \ln 3\leq 2 \ln x,
\end{gather*}if $c_{0}(\varepsilon)$ is chosen large enough. Hence,
\begin{align}
&\frac{li(2l_{1}x+l_{2}-1)-li(l_{1}x+l_{2}-1)}{\varphi(l_{1})}=\frac{1}{\varphi(l_{1})}\int_{l_{1}x+l_{2}-1}^{2l_{1}x+l_{2}-1}
\frac{dt}{\ln t}\geq \notag\\
&\geq \frac{l_{1} x}{\varphi(l_{1})\ln (2l_{1}x+l_{2}-1)} \geq
\frac{l_{1} x}{2\varphi(l_{1})\ln x} .\label{li_est_for_T6}
\end{align}Let us show that
\begin{equation}\label{|R|<x/ln_x_for_T6}
|R| \leq \frac{l_{1} x}{4\varphi(l_{1})\ln x}.
\end{equation}Since $l_{1}/\varphi (l_{1})\geq 1$, we see from \eqref{|R|<_for_T6} that it is sufficient to show that   \[
c_{6} x\,\textup{exp}(-c_{7}\sqrt{\ln x}) \leq \frac{ x}{4\ln x}.
\]This inequality holds, if $c_{0}(\varepsilon)$ is chosen large enough. Thus, \eqref{|R|<x/ln_x_for_T6} is proved. From \eqref{P_L_A_est_for_T6}, \eqref{li_est_for_T6} and \eqref{|R|<x/ln_x_for_T6} we obtain
\begin{equation}\label{Basic_2_2_for_T6}
\#\mathcal{P}_{L, \mathcal{A}}(x) \geq \frac{l_{1} x}{4\varphi(l_{1})\ln x}.
\end{equation}

Now we prove \eqref{Basic_I_for_T6}. Since $l_{1}/\varphi (l_{1})\geq 1$, we see from \eqref{Basic_2_1_for_T6} and \eqref{Basic_2_2_for_T6} that it suffices to show that
\begin{equation}\label{exp<P_L_A_for_T6}
c_{4} x\,\textup{exp}(-c_{5}\sqrt{\ln x})\leq \frac{x}{4 (\ln x)^{100 k^{2}+1}}.
\end{equation}Taking logarithms, we obtain
\[
\ln (c_{4}) +\ln x - c_{5}\sqrt{\ln x} \leq \ln x - \ln 4 - 100k^{2}\ln\ln x - \ln\ln x
\]or, that is equivalent,
\[
100 k^{2}\ln\ln x \leq c_{5}\sqrt{\ln x} - \ln\ln x - \ln (4 c_{4}).
\]Since $k \leq (\ln x)^{1/5}$, we have
\[
100 k^{2}\ln\ln x \leq 100 (\ln x)^{2/5}\ln\ln x.
\]The inequality
\[
100 (\ln x)^{2/5}\ln\ln x \leq c_{5}\sqrt{\ln x} - \ln\ln x - \ln (4 c_{4})
\]holds, if $c_{0}(\varepsilon)$ is chosen large enough. The inequality \eqref{exp<P_L_A_for_T6} is proved. Thus, \eqref{Basic_I_for_T6} is proved.

II) Let us show that part (1) of Hypothesis \ref{Hypothesis_1} holds. We show that
\begin{equation}\label{Basic_II_for_T6}
S:= \sum_{1 \leq r \leq x^{1/3}}\max_{b\in \mathbb{Z}}\Bigl|\#\mathcal{A}(x; r, b) - \frac{\# \mathcal{A}(x)}{r}\Bigr|\leq \frac{\#\mathcal{A}(x)}{(\ln x)^{100 k^{2}}}.
\end{equation}Let $1 \leq r \leq x^{1/3}$, $b\in \mathbb{Z}$. We have
\begin{gather*}
\mathcal{A}(x)=\{x\leq n < 2x\},\\
\mathcal{A}(x; r, b) = \{x \leq n < 2x:\ n \equiv b\ \text{(mod r)}\}.
\end{gather*}Hence,
\begin{align}
\#\mathcal{A}(x)=x;\label{A(x)_view}\\
\#\mathcal{A}(x; r, b)=\frac{x}{r}+\rho,\quad &|\rho|\leq 1.\label{A(x,r,b)_view}
\end{align}We obtain
\[
\Bigl|\#\mathcal{A}(x; r, b) - \frac{\# \mathcal{A}(x)}{r}\Bigr|=
|\rho| \leq 1.
\]Hence,
\[
S \leq  x^{1/3}.
\]Thus, to prove \eqref{Basic_II_for_T6} it suffices to show that
\[
x^{1/3}\leq \frac{x}{(\ln x)^{100 k^{2}}}
\]or, that is equivalent,
\[
(\ln x)^{100 k^{2}} \leq x^{2/3}.
\]Taking logarithms, we obtain
\[
100 k^{2}\ln\ln x  \leq \frac{2}{3}\ln x.
\]Since $k \leq (\ln x)^{1/5}$, we have
\[
100 k^{2}\ln\ln x  \leq 100 (\ln x)^{2/5}\ln\ln x .
\]The inequality
\[
100 (\ln x)^{2/5}\ln\ln x \leq \frac{2}{3}\ln x
\]holds, if $c_{0}(\varepsilon)$ is chosen large enough. Thus, \eqref{Basic_II_for_T6} is proved.

III) Let us show that part (3) of Hypothesis \ref{Hypothesis_1} holds. Let us show that for any integer $r$ with $1\leq r < x^{1/3}$ we have
\begin{equation}\label{Basic_III_for_T6}
\max_{b\in\mathbb{Z}}\#\mathcal{A}(x; r, b) \leq 2\, \frac{\#\mathcal{A}(x)}{r}.
\end{equation}Let $1\leq r < x^{1/3}$ and $b\in \mathbb{Z}$. We may assume that $c_{0}(\varepsilon)\geq 2$. Hence,
\[
r \leq x^{1/3} \leq x.
\]Applying \eqref{A(x)_view} and \eqref{A(x,r,b)_view}, we obtain
\[
\#\mathcal{A}(x; r, b)\leq \frac{x}{r}+1\leq 2\,\frac{x}{r}= 2\,\frac{\#\mathcal{A}(x)}{r},
\]and \eqref{Basic_III_for_T6} is proved. Thus, $(\mathcal{A}, \mathcal{L}, \mathcal{P}, B, x, 1/3)$ satisfy Hypothesis \ref{Hypothesis_1}.

We have
\[
\textup{exp}(c_{1}\sqrt{\ln x}) \leq x^{1/5},\qquad \ln x \leq x^{1/5},
\]if $c_{0}(\varepsilon)$ is chosen large enough. Since $1\leq B \leq \textup{exp}(c_{1}\sqrt{\ln x})$, we obtain $1 \leq B \leq x^{1/5}$. Let $L=l_{1}n+ l_{2}\in \mathcal{L}$. Applying \eqref{l_1<ln_x}, we have $1\leq l_{1} \leq x^{1/5}$, $1\leq l_{2} \leq x^{1/5}$. Thus, the assumption of Proposition \ref{Proposition_Maynard} holds and there are nonnegative weights $w_{n}= w_{n}(\mathcal{L})$ satisfying the statement of Proposition \ref{Proposition_Maynard}. Since in Proposition \ref{Proposition_Maynard} the implied constants in \eqref{Maynard1} -- \eqref{Maynard5} depend only on $\alpha$, $\theta$ and on the implied constants from Hypothesis \textup{\ref{Hypothesis_1}}, and in our case these constants are absolute ($\alpha=1/5$, $\theta=1/3$ and see \eqref{Basic_I_for_T6}, \eqref{Basic_II_for_T6} and \eqref{Basic_III_for_T6}), we see that in our case the implied constants in \eqref{Maynard1} -- \eqref{Maynard5} are positive and absolute. Finally, let us denote $c_{1}$ by $\vartheta$. Lemma \ref{L_wn_for_T6} is proved.

\begin{lemma}\label{Basic_Lemma_for_T6}
There are positive absolute constants $c$ and $C$ such that the following holds. Let $\varepsilon$ be a real number with $0< \varepsilon < 1$. Then there is a number $c_{0}(\varepsilon)>0$, depending only on $\varepsilon$, such that if  $x\in \mathbb{N}$, $y\in\mathbb{R}$, $m\in \mathbb{Z}$, $q\in \mathbb{Z}$, $a\in\mathbb{Z}$ are such that $c_{0}(\varepsilon) \leq y \leq \ln x$, $1 \leq m \leq c\cdot\varepsilon \ln y$, $1 \leq q \leq y^{1-\varepsilon}$ and $(a,q)=1$, then
\begin{align*}
&\#\{qx< p_{n}\leq 2qx - 5q:\ p_{n}\equiv\ldots\equiv p_{n+m}\equiv a\ (\textup{\text{mod }}q),\ p_{n+m} - p_{n} \leq y\}\\
&\geq \pi(2qx)\left(\frac{y}{2q\ln x}\right)^{\textup{exp} (Cm)}.
\end{align*}
\end{lemma}
\textsc{Proof of Lemma \ref{Basic_Lemma_for_T6}.} Let $\mathcal{A}=\mathbb{N}$, $\mathcal{P}=\mathbb{P}$, $\alpha = 1/5$, $\theta = 1/3$, let $C_{0}=C(1/5, 1/3)>0$ be the absolute constant in Proposition \textup{\ref{Proposition_Maynard}}. Let $c_{0}(\varepsilon)$ be the quantity in Lemma \ref{L_wn_for_T6}. We choose $c(\varepsilon)$ later; this number is large enough. Let $c(\varepsilon) \geq c_{0}(\varepsilon)$. Let $x\in \mathbb{N}$, $y\in \mathbb{R}$, $q\in\mathbb{Z}$ are such that
\begin{equation}\label{y_zona_T6}
c(\varepsilon) \leq y \leq \ln x,
\end{equation}
\begin{equation}\label{q_zona_T6}
1\leq q \leq y^{1-\varepsilon}.
\end{equation}By Lemma \ref{L_wn_for_T6}, there is a positive integer $B$ such that \eqref{B_for_T6} holds. We assume that
\begin{equation}\label{k_zona_T6}
\widetilde{C}_{0} \leq k \leq y^{\varepsilon/14},
\end{equation} where $\widetilde{C}_{0}>0$ is an absolute constant. We choose $\widetilde{C}_{0}$ later; this number is large enough. We may assume that $\widetilde{C}_{0} \geq C_{0}$. It follows from \eqref{y_zona_T6} and \eqref{k_zona_T6} that $k\leq (\ln x)^{1/5}$. Thus, \eqref{L_w_n_for_T6_1} holds. Let \eqref{L_w_n_for_T6_2} -- \eqref{L_w_n_for_T6_4} hold. Let $\mathcal{L}=\{L_{1},\ldots, L_{k}\}$ be an admissible set of $k$ linear functions, where $L_{i}(n) = q n + a+ q b_{i}$, $i=1,\ldots, k$, $b_{1},\ldots, b_{k}$ are positive integers with $b_{1}<\ldots < b_{k}$ and $q b_{k} \leq \eta y$.  Then (see Lemma \ref{L_wn_for_T6}) the assumption of Proposition \textup{\ref{Proposition_Maynard}} holds and there are nonnegative weights $w_{n}= w_{n}(\mathcal{L})$ satisfying the statement of Proposition \textup{\ref{Proposition_Maynard}}; the implied constants in \eqref{Maynard1} -- \eqref{Maynard5} are positive and absolute. We write $\mathcal{L}=\mathcal{L}(\textbf{b})$ for such a set given by $b_{1},\ldots, b_{k}$.

 Let $m$ be a positive integer. We consider
\begin{align}
S=&\sum_{\substack{1\leq b_{1}<\ldots< b_{k}\\ qb_{k}\leq \eta y\\ \mathcal{L}=\mathcal{L}(\textbf{b})\text{ admissible}}} \sum_{n\in \mathcal{A}(x)}\Bigl(\sum_{i=1}^{k}\textbf{1}_{\mathcal{P}}(L_{i}(n)) - m -\notag\\
&-k\sum_{i=1}^{k}
\sum_{\substack{p|L_{i}(n)\\ p< x^{\rho}\\ p\nmid B}} 1 -
k\sum_{\substack{1\leq b \leq 2\eta y\\ L=qt+b\notin \mathcal{L}}} \textbf{1}_{\mathcal{S}(\rho; B)}(L(n))\Bigr)w_{n}(\mathcal{L})=\notag\\
&=\sum_{\substack{1\leq b_{1}<\ldots< b_{k}\\ qb_{k}\leq \eta y\\ \mathcal{L}=\mathcal{L}(\textbf{b})\text{ admissible}}} \sum_{n\in \mathcal{A}(x)}A_{n}(\mathcal{L})w_{n}(\mathcal{L}).\label{S_Basic_T6}
\end{align}

Let $\mathcal{L}=\{L_{1},\ldots, L_{k}\}$ and $n$ be in the range of summation of $S$ such that $A_{n}(\mathcal{L})>0$.
Then the following statements hold.\newline
1) The number of primes among $L_{1}(n),\ldots, L_{k}(n)$ is at least $m+1$.\newline
2) For any $1\leq i \leq k$,  $L_{i}(n)$ has not a prime factor $p$ such that $p< x^{\rho}$ and $p\nmid B$.\newline
3) For any linear function $L=qt+b\notin\mathcal{L}$, where $b$ is an integer with $1\leq b \leq 2\eta y$,  $L(n)$ has a prime factor $p$ such that $p<x^{\rho}$ and $p\nmid B$ (we choose $\rho$ so that $x^{\rho}$ is not an integer; therefore the statements $p \leq x^{\rho}$ and $p< x^{\rho}$ are equivalent). Since $$L(n)> n\geq x>x^{\rho},$$ we see that $L(n)$ is not a prime number.

As a consequence we obtain the following statements.\newline
 i) None of $n\in \mathcal{A}(x)$ can make a positive contribution to $S$ from two different admissible sets (since if $n$ makes a positive contribution for some admissible set $\mathcal{L}=\{L_{1},\ldots, L_{k}\}$, then the numbers $L_{1}(n),\ldots, L_{k}(n)$ are uniquely determined as the integers in $[qn+1, qn+2\eta y]$ with no prime factors $p$ such that $p< x^{\rho}$ and $p\nmid B$).\newline
 ii) If $\mathcal{L}=\{L_{1},\ldots, L_{k}\}$ and $n$ are in the range of summation of $S$ such that $A_{n}(\mathcal{L})>0$, then there can be no primes in the interval $[qn+1, qn+ 2\eta y]$ apart from possibly $L_{1}(n),\ldots, L_{k}(n)$, and so the primes counted in this way must be consecutive.

  Let $\mathcal{L}=\{L_{1},\ldots, L_{k}\}$ and $n$ be in the range of summation of $S$ such that $A_{n}(\mathcal{L})>0$. Let $1 \leq i \leq k$. If $p| L_{i}(n)$ and $p\nmid B$, then $p \geq x^{\rho}$. Let
 \[
 \Omega = \{p:\ p|L_{i}(n)\text{ and }p\nmid B\}.
 \]We have
 \[
 x^{\rho\#\Omega}\leq \prod_{p\in \Omega} p \leq L_{i}(n).
 \]Since
 \begin{gather*}
 q\leq y^{1-\varepsilon} \leq y\leq \ln x,\\
 a+qb_{i}\leq 2 \eta y\leq \ln x,
 \end{gather*}we have
\[
L_{i}(n) =qn + a+qb_{i}\leq n\ln x+ \ln x \leq 2x\ln x+\ln x \leq x^{2},
\]if $c(\varepsilon)$ is chosen large enough. Hence, $\rho\#\Omega \leq 2$, i.\,e.
\[
\#\Omega \leq \frac{2}{\rho}.
\]We have
\begin{gather*}
\prod_{\substack{p|L_{i}(n)\\ p\nmid B}} 4 = \prod_{p\in \Omega} 4 = 4^{\#\Omega}\leq 4^{2/\rho}=
e^{(2/\rho)\ln 4}\leq e^{4/\rho},\\
\prod_{i=1}^{k} \prod_{\substack{p|L_{i}(n)\\ p\nmid B}} 4 \leq e^{(4k)/\rho}.
\end{gather*}Thus, if $\mathcal{L}=\{L_{1},\ldots, L_{k}\}$ and $n$ are in the range of summation of $S$ such that $A_{n}(\mathcal{L})>0$, then (see \eqref{Maynard1})
\begin{equation}\label{w_n_UP_T6}
w_{n}(\mathcal{L}) \leq C (\ln R)^{2k} e^{(4k)/\rho},
\end{equation}where $C>0$ is an absolute constant.

Let $\mathcal{L}=\{L_{1},\ldots, L_{k}\}$ be in the range of summation of $S$. We consider
\begin{align*}
\widetilde{S}(\mathcal{L})&= \sum_{n\in \mathcal{A}(x)}\Bigl(\sum_{i=1}^{k}\textbf{1}_{\mathcal{P}}(L_{i}(n)) - m -k\sum_{i=1}^{k}
\sum_{\substack{p|L_{i}(n)\\ p< x^{\rho}\\ p\nmid B}} 1 -\\
&-k\sum_{\substack{1\leq b \leq 2\eta y\\ L=qt+b\notin \mathcal{L}}} \textbf{1}_{\mathcal{S}(\rho; B)}(L(n))\Bigr)w_{n}(\mathcal{L})= S_{1} - S_{2} - S_{3} - S_{4}.
\end{align*}We are going to obtain a lower bound for $\widetilde{S}(\mathcal{L})$. We write $w_{n}$ instead of $w_{n}(\mathcal{L})$ for the brevity. Let $1 \leq i \leq k$. Since $\#\mathcal{A}(x)= x$, we have (see \eqref{Maynard3})
\begin{align*}
&\sum_{n\in \mathcal{A}(x)}\textup{\textbf{1}}_{\mathcal{P}}(L_{i}(n))w_{n}\geq
\bigl(1+ o(1)\bigr)\frac{B^{k-1}}{\varphi(B)^{k-1}}\mathfrak{S}_{B}(\mathcal{L})
\frac{\varphi (q)}{q}\cdot\\
&\cdot\#\mathcal{P}_{L_{i}, \mathcal{A}}(x)(\ln R)^{k+1}J_{k}+O\Bigl(\frac{B^{k}}{\varphi(B)^{k}} \mathfrak{S}_{B}(\mathcal{L}) x(\ln R)^{k-1} I_{k}\Bigr).
\end{align*}Hence,
\begin{align*}
&S_{1} = \sum_{n\in \mathcal{A}(x)}\sum_{i=1}^{k} \textbf{1}_{\mathcal{P}}(L_{i}(n))w_{n}=
 \sum_{i=1}^{k}  \sum_{n\in \mathcal{A}(x)} \textbf{1}_{\mathcal{P}}(L_{i}(n))w_{n}\geq \\
 &\geq \bigl(1+ o(1)\bigr)\frac{B^{k-1}}{\varphi(B)^{k-1}}\mathfrak{S}_{B}(\mathcal{L})
\frac{\varphi (q)}{q}(\ln R)^{k+1}J_{k}\cdot\\
&\cdot\sum_{i=1}^{k}\#\mathcal{P}_{L_{i}, \mathcal{A}}(x)+O\Bigl(k\frac{B^{k}}{\varphi(B)^{k}} \mathfrak{S}_{B}(\mathcal{L}) x(\ln R)^{k-1} I_{k}\Bigr)=\\
&=\bigl(1+ o(1)\bigr)\frac{B^{k}}{\varphi(B)^{k}}\mathfrak{S}_{B}(\mathcal{L})
(\ln R)^{k+1}J_{k}\frac{\varphi(B)}{B}\frac{\varphi (q)}{q}\sum_{i=1}^{k}\#\mathcal{P}_{L_{i}, \mathcal{A}}(x)+\\
&+o\Bigl(\frac{B^{k}}{\varphi(B)^{k}} \mathfrak{S}_{B}(\mathcal{L}) x(\ln R)^{k} I_{k}\Bigr)=S_{1}^{\prime}+S_{1}^{\prime\prime},
\end{align*}since
\[
0<\frac{k}{\ln R}\leq \frac{(\ln x)^{1/5}}{(1/9)\ln x}\to 0\ \text{ as $x\to +\infty$}.
\]We have shown (see \eqref{Basic_2_2_for_T6}), if $x \geq c_{0}$, where $c_{0}> 0$ is an absolute constant, then for any $L\in \mathcal{L}$
\[
\#\mathcal{P}_{L, \mathcal{A}}(x) \geq \frac{q x}{4\varphi(q)\ln x}.
\]We may assume that $c(\varepsilon) \geq c_{0}$. Since $\varphi(B)/B \geq 1/2$ (see \eqref{B_for_T6}), we obtain
\[
\frac{\varphi(B)}{B}\frac{\varphi (q)}{q}\sum_{i=1}^{k}\#\mathcal{P}_{L_{i}, \mathcal{A}}(x)\geq
\frac{kx}{8\ln x} = \frac{kx}{72\ln R}.
\]We have $|o(1)|\leq 1/2$ in $S_{1}^{\prime}$, if $x \geq c^{\prime}$, where $c^{\prime}>0$ is an absolute constant. We may assume that $c(\varepsilon)\geq c^{\prime}$. Since (see \eqref{Maynard_Jk_est})
\[
J_{k}\geq c^{\prime\prime}\frac{\ln k}{k}I_{k},
\]where $c^{\prime\prime}>0$ is an absolute constant, we obtain
\[
S_{1}^{\prime}\geq \frac{c^{\prime\prime}}{144}\frac{B^{k}}{\varphi(B)^{k}}\mathfrak{S}_{B}(\mathcal{L})x
(\ln R)^{k}I_{k}\ln k.
\]We have
\[
|S_{1}^{\prime\prime}|\leq \frac{c^{\prime\prime}}{288}\frac{B^{k}}{\varphi(B)^{k}} \mathfrak{S}_{B}(\mathcal{L}) x(\ln R)^{k} I_{k}\leq
\frac{c^{\prime\prime}}{288}\frac{B^{k}}{\varphi(B)^{k}} \mathfrak{S}_{B}(\mathcal{L}) x(\ln R)^{k} I_{k}\ln k,
\]if $c(\varepsilon)$ is chosen large enough. We obtain
\begin{align}
S_{1} &\geq \frac{c^{\prime\prime}}{288}\frac{B^{k}}{\varphi(B)^{k}} \mathfrak{S}_{B}(\mathcal{L}) x(\ln R)^{k} I_{k}\ln k=\notag\\
&= c\,\frac{B^{k}}{\varphi(B)^{k}} \mathfrak{S}_{B}(\mathcal{L}) x(\ln R)^{k} I_{k}\ln k,\label{S1_LOW_for_T6}
\end{align}where $c>0$ is an absolute constant.

We have (see \eqref{Maynard2})
\begin{align}
S_{2} &= m\sum_{n\in \mathcal{A}(x)} w_{n} = m (1+ o(1))\frac{B^{k}}{\varphi(B)^{k}}\mathfrak{S}_{B}(\mathcal{L})x(\ln R)^{k}I_{k}\geq \notag\\
& \geq \frac{m}{2} \frac{B^{k}}{\varphi(B)^{k}}\mathfrak{S}_{B}(\mathcal{L})x(\ln R)^{k}I_{k},\label{S2_low_est_for_T6}
\end{align}if $c(\varepsilon)$ is chosen large enough. Applying \eqref{Maynard5}, we have
\begin{align*}
S_{3} &= k \sum_{n\in \mathcal{A}(x)} \sum_{i=1}^{k}\sum_{\substack{p|L_{i}(n)\\ p< x^{\rho}\\ p\nmid B}} w_{n}=
k\sum_{i=1}^{k}\Bigl(\sum_{n\in \mathcal{A}(x)} \bigl(\sum_{\substack{p|L_{i}(n)\\ p< x^{\rho}\\ p\nmid B}} 1\bigr) w_{n}\Bigr)\leq\\
&\leq c_{2} \rho^{2} k^{6}(\ln k)^{2}
\frac{B^{k}}{\varphi(B)^{k}}\mathfrak{S}_{B}(\mathcal{L})x (\ln R)^{k}I_{k},
\end{align*}where $c_{2}>0$ is an absolute constant. Let $c_{3}>0$ be an absolute constant such that
\[
c_{2}c_{3}^{2}\leq \frac{1}{12}\quad\text{and}\quad \frac{c_{3}}{j^{3}\ln j}\leq \frac{1}{30} \text{ for any }j\geq 2.
\]We choose a number $\rho$ in the interval
\begin{equation}\label{Interval_for_rho_T6}
\left[\frac{c_{3}}{2k^{3}\ln k}, \frac{c_{3}}{k^{3}\ln k}\right]
\end{equation}so that $x^{\rho}$ is not an integer. It is clear that $\rho \leq 1/30$. Let us show that the first inequality in \eqref{L_w_n_for_T6_2} holds. It suffices to show that
\[
\frac{k(\ln\ln x)^{2}}{\ln x}\leq \frac{c_{3}/2}{k^{3}\ln k}.
\]This inequality is equivalent to
\[
k^{4}\ln k (\ln \ln x)^{2} \leq \frac{c_{3}}{2} \ln x.
\]Since $k \leq (\ln x)^{1/5}$, we have
\[
k^{4}\ln k (\ln \ln x)^{2} \leq \frac{1}{5}(\ln x)^{4/5}(\ln \ln x)^{3} \leq \frac{c_{3}}{2} \ln x,
\]if $c(\varepsilon)$ is chosen large enough. Thus, inequalities in \eqref{L_w_n_for_T6_2} hold. We have (see \eqref{S2_low_est_for_T6})
\begin{align}
S_{3} &\leq c_{2} \frac{c_{3}^{2}}{k^{6}(\ln k)^{2}} k^{6}(\ln k)^{2}
\frac{B^{k}}{\varphi(B)^{k}}\mathfrak{S}_{B}(\mathcal{L})x (\ln R)^{k}I_{k}\leq\notag\\
&\leq \frac{1}{12}\frac{B^{k}}{\varphi(B)^{k}}\mathfrak{S}_{B}(\mathcal{L})x (\ln R)^{k}I_{k}\leq
\frac{m}{12}\frac{B^{k}}{\varphi(B)^{k}}\mathfrak{S}_{B}(\mathcal{L})x (\ln R)^{k}I_{k}\leq\notag\\
&\leq \frac{1}{6}S_{2}.\label{S3_1/6_S2}
\end{align}

Now we estimate
\begin{align*}
S_{4} &= k\sum_{n\in \mathcal{A}(x)}\sum_{\substack{1\leq b \leq 2\eta y\\ L=qt+b\notin \mathcal{L}}} \textbf{1}_{\mathcal{S}(\rho; B)}(L(n))w_{n} =\\
&= k \sum_{\substack{1\leq b \leq 2\eta y\\ L=qt+b\notin \mathcal{L}}}\sum_{n\in \mathcal{A}(x)}\textbf{1}_{\mathcal{S}(\rho; B)}(L(n))w_{n}.
\end{align*}Let $b$ to be in the range of summation of $S_{4}$. Then $L=qt+b\notin\mathcal{L}$ and
\[
\Delta_{L}=q^{k+1}\prod_{i=1}^{k}|(a+qb_{i}) - b| \neq 0.
\]Since $1\leq B \leq x^{1/5}$, we have (see \eqref{Maynard4})
\[
\sum_{n\in \mathcal{A}(x)}\textbf{1}_{\mathcal{S}(\rho; B)}(L(n))w_{n}\leq
 \frac{c_{4}}{\rho}\frac{\Delta_{L}}{\varphi(\Delta_{L})}\frac{B}{\varphi(B)}
 \frac{B^{k}}{\varphi(B)^{k}}\mathfrak{S}_{B}(\mathcal{L})x(\ln R)^{k-1}I_{k},
\]where $c_{4}>0$ is an absolute constant. Since $B/\varphi(B) \leq 2$ and $\rho$ lies in the interval \eqref{Interval_for_rho_T6}, we obtain
\begin{align*}
\sum_{n\in \mathcal{A}(x)}\textbf{1}_{\mathcal{S}(\rho; B)}(L(n))w_{n}&\leq
 \frac{4c_{4}}{c_{3}}k^{3}\ln k\frac{\Delta_{L}}{\varphi(\Delta_{L})} \frac{B^{k}}{\varphi(B)^{k}}\mathfrak{S}_{B}(\mathcal{L})x(\ln R)^{k-1}I_{k}=\\
 &= c_{5} k^{3}\ln k\frac{\Delta_{L}}{\varphi(\Delta_{L})} \frac{B^{k}}{\varphi(B)^{k}}\mathfrak{S}_{B}(\mathcal{L})x(\ln R)^{k-1}I_{k}.
\end{align*}Hence,
\begin{equation}\label{S4_est_up_for_T6}
S_{4} \leq c_{5} k^{4}\ln k \frac{B^{k}}{\varphi(B)^{k}}\mathfrak{S}_{B}(\mathcal{L})x(\ln R)^{k-1}I_{k}
  \sum_{\substack{1\leq b \leq 2\eta y\\ L=qt+b\notin \mathcal{L}}} \frac{\Delta_{L}}{\varphi(\Delta_{L})}.
\end{equation}We put
\begin{equation}\label{Def_c6_for_T6}
c_{6}=36Cc_{5},
\end{equation}where $C>0$ is the absolute constant in Lemma \ref{L:sum_Delta_L}. We put
\begin{equation}\label{Def_eta_for_T6}
\eta = \frac{1}{12 c_{6} k^{4}(\ln k)^{2}\ln\ln (q+2)}.
\end{equation}Let us show that
\begin{equation}\label{eta_eq_for_T6}
(\ln x)^{-9/10}\leq 2\eta\leq 1.
\end{equation} The second inequality in \eqref{eta_eq_for_T6} is equivalent to the inequality
\[
6c_{6} k^{4}(\ln k)^{2}\ln\ln (q+2) \geq 1.
\]We may assume that $\widetilde{C}_{0}\geq 3$; therefore $\ln k \geq 1$. We have
\[
6c_{6} k^{4}(\ln k)^{2}\ln\ln (q+2) \geq 6c_{6}(\ln\ln 3) k^{4}
\geq 6c_{6}(\ln\ln 3) {\widetilde{C}_{0}}^{4} \geq 1,
\]if $\widetilde{C}_{0}$ is chosen large enough. The first inequality in \eqref{eta_eq_for_T6} is equivalent to the inequality
\[
6c_{6} k^{4}(\ln k)^{2}\ln\ln (q+2) \leq (\ln x)^{9/10}.
\]Since $q \leq \ln x$ and $k\leq (\ln x)^{1/5}$, we have
\begin{align*}
6c_{6} k^{4}(\ln k)^{2}\ln\ln (q+2) &\leq 6c_{6} (\ln x)^{4/5}\frac{1}{25}(\ln \ln x)^{2}\ln\ln(\ln x+2)\leq\\
&\leq (\ln x)^{9/10},
\end{align*}if $c(\varepsilon)$ is chosen large enough. Thus, \eqref{eta_eq_for_T6} holds. We have $x \geq c$, where $c$ is the absolute constant in Lemma \ref{L:sum_Delta_L}, if $c(\varepsilon)$ is chosen large enough. Applying Lemma \ref{L:sum_Delta_L} and taking into account that $\ln (k+1) \leq 2 \ln k$, we have
\begin{align*}
&\sum_{\substack{1\leq b \leq 2\eta y\\ L=qt+b\notin \mathcal{L}}} \frac{\Delta_{L}}{\varphi(\Delta_{L})} \leq \sum_{\substack{1\leq b \leq 2\eta \ln x\\ L=qt+b\notin \mathcal{L}}} \frac{\Delta_{L}}{\varphi(\Delta_{L})}\leq
4C\ln\ln (q+2)(\ln k) \eta \ln x=\\
&= 36C \ln\ln (q+2)(\ln k) \eta \ln R.
\end{align*}Substituting this estimate into \eqref{S4_est_up_for_T6}, we obtain (see also \eqref{Def_c6_for_T6}, \eqref{Def_eta_for_T6} and \eqref{S2_low_est_for_T6})
\begin{align}
&S_{4} \leq 36Cc_{5} k^{4}(\ln k)^{2} \frac{B^{k}}{\varphi(B)^{k}}\mathfrak{S}_{B}(\mathcal{L})x(\ln R)^{k}I_{k}
   \eta  \ln\ln (q+2)= c_{6} k^{4}(\ln k)^{2} \cdot\notag\\
&\cdot  \frac{B^{k}}{\varphi(B)^{k}}\mathfrak{S}_{B}(\mathcal{L})x(\ln R)^{k}I_{k}
    \ln\ln (q+2) \frac{1}{12 c_{6} k^{4}(\ln k)^{2}\ln\ln (q+2)}=\notag\\
&= \frac{1}{12} \frac{B^{k}}{\varphi(B)^{k}}\mathfrak{S}_{B}(\mathcal{L})x(\ln R)^{k}I_{k}
\leq \frac{m}{12} \frac{B^{k}}{\varphi(B)^{k}}\mathfrak{S}_{B}(\mathcal{L})x(\ln R)^{k}I_{k}\leq
\frac{1}{6} S_{2}.\label{S4_1/6_S2}
\end{align}From \eqref{S3_1/6_S2} and \eqref{S4_1/6_S2} we obtain
\[
\widetilde{S}(\mathcal{L}) = S_{1}- S_{2}-S_{3}-S_{4}\geq S_{1} - \frac{4}{3}S_{2}.
\]We have (see \eqref{Maynard2})
\begin{align*}
S_{2} &= m\sum_{n\in \mathcal{A}(x)} w_{n} = m (1+ o(1))\frac{B^{k}}{\varphi(B)^{k}}\mathfrak{S}_{B}(\mathcal{L})x(\ln R)^{k}I_{k} \leq\\
&\leq \frac{3}{2}m \frac{B^{k}}{\varphi(B)^{k}}\mathfrak{S}_{B}(\mathcal{L})x(\ln R)^{k}I_{k},
\end{align*}if $c(\varepsilon)$ is chosen large enough. Applying \eqref{S1_LOW_for_T6}, where we replace $c$ by $3 c_{1}$, we obtain
\[
\widetilde{S}(\mathcal{L})\geq \frac{B^{k}}{\varphi(B)^{k}}\mathfrak{S}_{B}(\mathcal{L})x(\ln R)^{k}I_{k}(3c_{1}\ln k - 2m),
\]where $c_{1}>0$ is an absolute constant. We put
\begin{gather}
\widetilde{c} = \widetilde{C}_{0}+ \frac{1}{c_1},\label{c_wave_view_T6}\\
k=\lceil \textup{exp}(\widetilde{c}m)\rceil.\label{k_view_T6}
\end{gather}It is not hard to see that
\[
k\geq \widetilde{C}_{0}\quad\text{ and }\quad 3c_{1}\ln k - 2m \geq m.
\]Since $m$ is a positive integer, we see that
\[
3c_{1}\ln k - 2m \geq 1.
\]Hence,
\[
\widetilde{S}(\mathcal{L})\geq \frac{B^{k}}{\varphi(B)^{k}}\mathfrak{S}_{B}(\mathcal{L})x(\ln R)^{k}I_{k}.
\]Since $B^{k}/\varphi(B)^{k}\geq 1$, $\ln R = (1/9)\ln x$, $\mathfrak{S}_{B}(\mathcal{L})\geq \textup{exp}(-c_{2}k)$ and $I_{k}\geq c_{3}(2k\ln k)^{-k}$, where $c_{2}$ and $c_{3}$ are positive absolute constants (see \eqref{exp_Ineq_for_sigma_B_L} and \eqref{Maynard_I_k_est}), we obtain
\[
\widetilde{S}(\mathcal{L})\geq \frac{1}{9^{k}}c_{3}(2k\ln k)^{-k}\textup{exp}(-c_{2}k)x(\ln x)^{k}
\geq \textup{exp}(-k^{2}) x(\ln x)^{k},
\]if $\widetilde{C}_{0}$ is chosen large enough. We obtain
\begin{align}
S&=\sum_{\substack{1\leq b_{1}<\ldots< b_{k}\\ qb_{k}\leq \eta y\\ \mathcal{L}=\mathcal{L}(\textbf{b})\text{ admissible}}}\widetilde{S}(\mathcal{L}) \geq\notag\\
&\geq \textup{exp}(-k^{2}) x(\ln x)^{k} \sum_{\substack{1\leq b_{1}<\ldots< b_{k}\\ qb_{k}\leq \eta y\\ \mathcal{L}=\mathcal{L}(\textbf{b})\text{ admissible}}} 1=\notag\\
&= \textup{exp}(-k^{2}) x(\ln x)^{k} S^{\prime}.\label{S_S_prime_Ineq_T6}
\end{align}

Now we obtain a lower bound for $S^{\prime}$. First let us show that
\begin{equation}\label{Ineq_for_Omega_T6}
2 \leq k \leq \Bigl[\frac{\eta y}{q}\Bigr]/2.
\end{equation}The first inequality obviously holds, since we may assume that $\widetilde{C}_{0} \geq 2$. To prove the second inequality it suffices to show that
\begin{equation}\label{Prepare_Ineq_for_Omega_T6}
2 k \leq \frac{\eta y}{q}.
\end{equation}We have (see \eqref{q_zona_T6}, \eqref{Def_eta_for_T6})
\[
\frac{\eta y}{q} \geq \eta y^{\varepsilon}= \frac{c_{4} y^{\varepsilon}}{k^{4}(\ln k)^{2}\ln\ln (q+2)},
\]where $c_{4}>0$ is an absolute constant. Thus, to prove \eqref{Prepare_Ineq_for_Omega_T6} it suffices to show that
\[
2k^{5}(\ln k)^{2}\ln\ln (q+2)\leq c_{4} y^{\varepsilon}.
\]In particular, from \eqref{q_zona_T6} it follows that $q\leq y$. Applying \eqref{k_zona_T6}, we have
\[
2k^{5}(\ln k)^{2}\ln\ln (q+2) \leq 2 y^{5\varepsilon/14}\frac{\varepsilon^{2}}{196}(\ln y)^{2}\ln\ln (y+2)
\leq c_{4} y^{\varepsilon},
\]if $c(\varepsilon)$ is chosen large enough. Thus, \eqref{Ineq_for_Omega_T6} is proved.

 We put
\[
\Omega = \Bigl\{ 1 \leq n \leq \Bigl[\frac{\eta y}{q}\Bigr]:\ (n,p)=1\ \forall p\leq k\Bigr\}.
\]Applying Lemma \ref{L_about_Phi(x,z)}, we have
\[
\#\Omega = \Phi \Bigl(\Bigl[\frac{\eta y}{q}\Bigr], k\Bigr)\geq c_{0} \frac{[\eta y/q]}{\ln k},
\]where $c_{0}>0$ is an absolute constant. In particular, from \eqref{Ineq_for_Omega_T6} it follows that $\eta y/ q \geq 4$, and hence
\[
\Bigl[\frac{\eta y}{q}\Bigr]\geq \frac{\eta y}{q} - 1\geq \frac{\eta y}{2q}.
\]We obtain
\begin{equation}\label{N_Omega_T6}
\#\Omega \geq c_{5} \frac{\eta y}{q\ln k},
\end{equation}where $c_{5}>0$ is an absolute constant. Let us show that
\begin{equation}\label{N_Omega_Ineq_2k_T6}
c_{5} \frac{\eta y}{q\ln k}\geq 2 k.
\end{equation}Applying \eqref{q_zona_T6} and \eqref{Def_eta_for_T6}, we have
\[
c_{5} \frac{\eta y}{q\ln k}\geq \frac{c_{6} y^{\varepsilon}}{k^{4}(\ln k)^{3}\ln\ln (q+2)},
\]where $c_{6}>0$ is an absolute constant. Therefore it suffices to show that
\[
2k^{5}(\ln k)^{3}\ln\ln (q+2) \leq c_{6} y^{\varepsilon}.
\]Applying \eqref{k_zona_T6} and taking into account that $q\leq y$, we have
\[
2k^{5}(\ln k)^{3}\ln\ln (q+2)\leq 2 y^{5\varepsilon/14}\Bigl(\frac{\varepsilon}{14}\Bigr)^{3}(\ln y)^{3}\ln\ln (y+2)
\leq c_{6} y^{\varepsilon},
\]if $c(\varepsilon)$ is chosen large enough. Thus, \eqref{N_Omega_Ineq_2k_T6} is proved.

  Let $b_{1} <\ldots < b_{k}$ be positive integers from the set $\Omega$. Let us show that for any prime $p$ with $p\nmid q$ there is an integer $m_{p}$ such that $m_{p}\not \equiv b_{i}\ (\text{mod }p)$ for all $1\leq i \leq k$. Let $p$ be a prime with $p\nmid q$. If $p > k$, then the statement is obvious. If $p \leq k$, then we may put $m_p=0$; from the definition of the set $\Omega$ it follows that $b_{i}\not\equiv 0$ (mod $p$) for all $1\leq i \leq k$. Thus, the statement is proved. By Lemma \ref{L:admissible_set}, $\mathcal{L}(\textbf{b})$ is an admissible set. Hence (see also Lemma \ref{L_Binom_C_n_k}, \eqref{N_Omega_T6}, \eqref{N_Omega_Ineq_2k_T6} and \eqref{Def_eta_for_T6}),
 \begin{align*}
 S^{\prime} &\geq \binom{\#\Omega}{k}\geq k^{-k}(\#\Omega - k)^{k} \geq k^{-k}\Bigl(c_{5} \frac{\eta y}{q\ln k} - k\Bigr)^{k}\geq\\
 &\geq k^{-k}\Bigl(\frac{c_{5}}{2} \frac{\eta y}{q\ln k}\Bigr)^{k}=
 k^{-k}\Bigl(c_{6} \frac{ y}{q\ln\ln(q+2)k^{4}(\ln k)^{3}}\Bigr)^{k}=\\
 &= \Bigl(\frac{y}{q\ln\ln (q+2)}\Bigr)^{k}\Bigl(\frac{c_{6}}{k^{5}(\ln k)^{3}}\Bigr)^{k},
 \end{align*}where $c_{6}>0$ is an absolute constant. We have
  \[
 \Bigl(\frac{c_{6}}{k^{5}(\ln k)^{3}}\Bigr)^{k} \geq \textup{exp}(-k^{2}),
 \]if $\widetilde{C}_{0}$ is chosen large enough. Hence,
 \[
 S^{\prime} \geq \Bigl(\frac{y}{q\ln\ln (q+2)}\Bigr)^{k} \textup{exp}(-k^{2}).
 \]Substituting this estimate into \eqref{S_S_prime_Ineq_T6}, we obtain
 \begin{align}
 S&\geq \textup{exp}(-2k^{2})x (\ln x)^{k} \Bigl(\frac{y}{q\ln\ln (q+2)}\Bigr)^{k}\geq\notag\\
 &\geq \textup{exp}(-2k^{5})x (\ln x)^{k} \Bigl(\frac{y}{q\ln\ln (q+2)}\Bigr)^{k}.\label{S_LOW_T6}
 \end{align}

 Now we obtain an upper bound for $S$. Applying \eqref{S_Basic_T6} and \eqref{w_n_UP_T6}, we have
 \begin{align*}
 &S\leq \sum_{\substack{1\leq b_{1}<\ldots< b_{k}\\ qb_{k}\leq \eta y\\ \mathcal{L}=\mathcal{L}(\textbf{b})\text{ admissible}}} \sum_{\substack{n\in \mathcal{A}(x):\\ A_{n}(\mathcal{L})>0}}A_{n}(\mathcal{L})w_{n}(\mathcal{L})\leq\\
 &\leq Ck (\ln R)^{2k} e^{(4k)/\rho}\sum_{\substack{1\leq b_{1}<\ldots< b_{k}\\ qb_{k}\leq \eta y\\ \mathcal{L}=\mathcal{L}(\textbf{b})\text{ admissible}}} \sum_{\substack{n\in \mathcal{A}(x):\\ A_{n}(\mathcal{L})>0}} 1.
 \end{align*} We have (see parts 1) -- 3), i) and ii) above)
 \begin{gather*}
 \sum_{\substack{1\leq b_{1}<\ldots< b_{k}\\ qb_{k}\leq \eta y\\ \mathcal{L}=\mathcal{L}(\textbf{b})\text{ admissible}}} \sum_{\substack{n\in \mathcal{A}(x):\\ A_{n}(\mathcal{L})>0}} 1 \leq\\
 \leq\#\{ x\leq n < 2x:\text{ there are $\geq (m+1)$ consecutive primes}\\
  \text{all congruent to $a$ (mod $q$) in the interval $[qn+1, qn + 2\eta y]$}\}\leq\\
  \leq\#\{ x\leq n < 2x:\text{ there are $\geq (m+1)$ consecutive primes}\\
  \text{all congruent to $a$ (mod $q$) in the interval $[qn+1, qn + y]$}\}:=N_{1}.
   \end{gather*}Hence,
 \[
 S \leq Ck (\ln R)^{2k} e^{(4k)/\rho} N_{1}.
 \]Since $\rho$ lies in the interval \eqref{Interval_for_rho_T6}, we have
 \[
 \frac{4k}{\rho} \leq \frac{8 k^{4}\ln k}{c_{3}}= c_{4}k^{4}\ln k,
 \]where $c_{4}>0$ is an absolute constant. Since $\ln R = (1/9)\ln x$, we have
 \[
 Ck (\ln R)^{2k} e^{(4k)/\rho}\leq C\frac{k}{9^{2k}}\textup{exp}(c_{4}k^{4}\ln k) (\ln x)^{2k}
 \leq \textup{exp}(k^{5})(\ln x)^{2k},
 \]if $\widetilde{C}_{0}$ is chosen large enough. Hence,
 \begin{equation}\label{S_UP_T6}
 S \leq \textup{exp}(k^{5})(\ln x)^{2k} N_{1}.
 \end{equation}From \eqref{S_LOW_T6} and \eqref{S_UP_T6} we obtain
 \begin{equation}\label{N1_LOW_T6}
 N_{1}\geq x \Bigl(\frac{y}{\ln x}\Bigr)^{k}  \Bigl(\frac{1}{q\ln\ln (q+2)}\Bigr)^{k}\textup{exp}(-3k^{5}).
 \end{equation}

 We define
 \begin{align*}
  \Omega_{1} &= \{ x\leq n \leq 2x-1:\text{ there are $\geq (m+1)$ consecutive primes}\\
  &\qquad\text{all congruent to $a$ (mod $q$) in the interval $[qn+1, qn + y]$}\};&&\\
 \Omega_{2} &= \{qx+1\leq p_{n} \leq q(2x-1)+y:\ p_{n}\equiv \ldots \equiv p_{n+m} \equiv
 a\ \text{(mod $q$)},\\
 &\qquad p_{n+m} - p_{n} \leq y\}.
 \end{align*}We put $N_{2}=\#\Omega_{2}$. Since $x$ is a positive integer, we have $N_{1}=\#\Omega_{1}$. Let us show that
  \begin{equation}\label{N1_N2_T6}
  N_{1} \leq (\lceil y\rceil +1)N_{2}.
  \end{equation}
   Let $n\in \Omega_{1}$. Then there are $\geq (m+1)$ consecutive primes all congruent to $a$ (mod $q$) in the interval $[qn+1, qn + y]$. Let $p$ be the first of them. Then $p\in \Omega_{2}$. We put
  \[
  \Lambda = \{j\in \mathbb{Z}:\ qj+ 1 \leq p \leq qj+y\}.
  \]We claim that
  \begin{equation}\label{Lambda_UP_T6}
  \#\Lambda \leq \lceil y \rceil +1.
  \end{equation}We put $I_{j}= [qj+1, qj+y]$, $j\in \mathbb{Z}$. Since $p\in I_{n}$, we have $\Lambda\neq \varnothing$. Let $l$ be the minimal element in $\Lambda$. We put $t=\lceil y\rceil + 1$. Then $t>y$. We have
    \[
  q(l+t)+1> q(l+t)=ql+qt\geq ql +t> ql+y\geq p.
  \]Hence, $p\notin I_{j}$ for $j \geq l+t$ and $j\leq l-1$. We obtain $\#\Lambda \leq t$. Thus, \eqref{Lambda_UP_T6} is proved; \eqref{N1_N2_T6} follows from \eqref{Lambda_UP_T6}. We have
  \[
  \lceil y\rceil + 1\leq y+2 \leq 2y,
  \]if $c(\varepsilon)$ is chosen large enough. Since
  \begin{align*}
  &N_{2} \leq \#\{qx+1\leq p_{n} \leq 2qx+y:\ p_{n}\equiv \ldots \equiv p_{n+m} \equiv
 a\ \text{(mod $q$)},\\
 &\qquad\qquad p_{n+m} - p_{n} \leq y\}=:N_{3},
  \end{align*}we obtain (see \eqref{N1_LOW_T6})
  \begin{equation}\label{N3_LOW_SUPER_T6}
  N_{3}\geq \frac{1}{2}\frac{x}{y}  \Bigl(\frac{y}{\ln x}\Bigr)^{k}  \Bigl(\frac{1}{q\ln\ln (q+2)}\Bigr)^{k}\textup{exp}(-3k^{5}).
  \end{equation}We put
  \begin{align}
  N_{4} &= \#\{qx< p_{n} \leq 2qx - 5q:\ p_{n}\equiv \ldots \equiv p_{n+m} \equiv
 a\ \text{(mod $q$)},\notag\\
 &\qquad\quad p_{n+m} - p_{n} \leq y\};\label{N4_view_T6}\\
  N_{5} &= \#\{2qx -5q< p_{n} \leq 2qx + y:\ p_{n}\equiv \ldots \equiv p_{n+m} \equiv
 a\ \text{(mod $q$)},\notag\\
 &\qquad\quad p_{n+m} - p_{n} \leq y\}.\notag
  \end{align}Then
  \begin{equation}\label{N3=N4+N5_T6}
  N_{3}= N_{4}+N_{5}.
  \end{equation}Since $q \leq y$, we have
  \begin{equation}\label{N5<2y_T6}
  N_{5} \leq 5q+ [y]\leq 5q+ y\leq 6 y.
  \end{equation} Let us show that
  \begin{equation}\label{y<T1_T6}
  y \leq \frac{1}{24}\frac{x}{y}  \Bigl(\frac{y}{\ln x}\Bigr)^{k}  \Bigl(\frac{1}{q\ln\ln (q+2)}\Bigr)^{k}\textup{exp}(-3k^{5})
  :=T_{1}.
  \end{equation}Since $q\leq y^{1-\varepsilon}\leq y$ and $k \leq y^{\varepsilon/14}$, we have
  \[
  T_{1} \geq \frac{1}{24}\frac{x}{y}  \Bigl(\frac{y^{\varepsilon}}{\ln x\ln\ln(y+2)}\Bigr)^{k} \textup{exp}(-3y^{5\varepsilon/14}).
  \]Therefore to prove \eqref{y<T1_T6} it suffices to show that
  \[
  y \leq \frac{1}{24}\frac{x}{y}  \Bigl(\frac{y^{\varepsilon}}{\ln x\ln\ln(y+2)}\Bigr)^{k} \textup{exp}(-3y^{5\varepsilon/14}).
  \]Taking logarithms, we obtain
  \[
  \ln y \leq -\ln 24 + \ln x-\ln y + k(\varepsilon \ln y-\ln\ln x-\ln\ln\ln (y+2)) - 3y^{5\varepsilon/14}
  \]or, that is equivalent,
  \[
  T_{2}:= 2\ln y + \ln 24 - \varepsilon k \ln y + k\ln\ln x + k\ln\ln\ln (y+2) + 3y^{5\varepsilon/14} \leq \ln x.
  \] Since $y \leq \ln x$, $0< \varepsilon < 1$, we have $k \leq (\ln x)^{\varepsilon / 14}\leq (\ln x)^{1/14}$. We have
  \begin{align*}
  &T_{2} \leq 2 \ln\ln x + \ln 24 + (\ln x)^{1/14}\ln\ln x+ (\ln x)^{1/14} \ln\ln\ln(\ln x+2)+\\
  &+ 3 (\ln x)^{5/14} \leq \ln x,
  \end{align*}if $c(\varepsilon)$ is chosen large enough. Thus, \eqref{y<T1_T6} is proved. From \eqref{N5<2y_T6} and \eqref{y<T1_T6} it follows that
  \begin{equation}\label{N5_UP_SUPER_T6}
  N_{5} \leq \frac{1}{4}\frac{x}{y}  \Bigl(\frac{y}{\ln x}\Bigr)^{k}  \Bigl(\frac{1}{q\ln\ln (q+2)}\Bigr)^{k}\textup{exp}(-3k^{5}).
  \end{equation} Applying \eqref{N3_LOW_SUPER_T6}, \eqref{N3=N4+N5_T6} and \eqref{N5_UP_SUPER_T6}, we obtain
  \begin{equation}\label{N4>T3_for_T6}
  N_{4} \geq \frac{1}{4}\frac{x}{y}  \Bigl(\frac{y}{\ln x}\Bigr)^{k}  \Bigl(\frac{1}{q\ln\ln (q+2)}\Bigr)^{k}\textup{exp}(-3k^{5})=: T_{3}.
  \end{equation}We have (see \eqref{Notation:Pi_func_est})
  \[
  \pi(2qx) \leq c_{1} \frac{2qx}{\ln(2qx)}\leq c_{1} \frac{2qx}{\ln x}=c_{2}\frac{qx}{\ln x},
  \]where $c_{1}>0$ and $c_{2}=2c_{1}>0$ are absolute constants. We obtain
  \begin{align*}
  &T_{3} = \frac{qx}{\ln x}\frac{\ln x}{qx} \frac{1}{4}\frac{x}{y}  \Bigl(\frac{y}{\ln x}\Bigr)^{k}  \Bigl(\frac{1}{q\ln\ln (q+2)}\Bigr)^{k}\textup{exp}(-3k^{5})\geq\\
  &\geq \frac{1}{4c_{2}} \pi(2qx)\Bigl(\frac{y}{\ln x}\Bigr)^{k-1}\frac{1}{q^{k+1} (\ln\ln (q+2))^{k}}\textup{exp}(-3k^{5}).
  \end{align*}Using the inequality $\ln (1+x)\leq x$, $x>0$, we obtain
  \[
  \ln\ln(q+2)\leq \ln(1+q)\leq q.
  \]Hence,
  \[
  \frac{1}{q^{k+1} (\ln\ln (q+2))^{k}}\geq \frac{1}{q^{2k+1}}
  \geq \frac{1}{q^{3k^{5}}}.
  \]We have
  \[
  4c_{2} \leq 2^{3k^{5}},
  \]if $\widetilde{C}_{0}$ is chosen large enough. We have
  \[
  T_{3} \geq  \pi(2qx)\Bigl(\frac{y}{\ln x}\Bigr)^{k-1}\frac{1}{(2eq)^{3k^{5}}}.
  \]We have
  \[
  3k^{5}\leq k^{6},
  \]if $\widetilde{C}_{0}$ is chosen large enough. Hence,
  \[
  \frac{1}{(2eq)^{3k^{5}}}\geq \frac{1}{(2eq)^{k^{6}}}.
  \]We have
  \[
  (2e)^{k^{6}} \leq 2^{k^{7}},
  \]if $\widetilde{C}_{0}$ is chosen large enough. It is clear that
  \[
  q^{k^{6}} \leq q^{k^{7}}.
  \] We obtain
  \[
  \frac{1}{(2eq)^{k^{6}}} \geq \frac{1}{(2q)^{k^{7}}}.
  \] Since (see \eqref{y_zona_T6})
  \[
  0< \frac{y}{\ln x} \leq 1,
  \]we have
  \[
  \Bigl(\frac{y}{\ln x}\Bigr)^{k-1} \geq \Bigl(\frac{y}{\ln x}\Bigr)^{k^{7}}.
  \]We obtain
  \[
  T_{3} \geq \pi(2qx)\Bigl(\frac{y}{2q\ln x}\Bigr)^{k^{7}}.
  \]We have (see \eqref{c_wave_view_T6}, \eqref{k_view_T6})
  \begin{equation}\label{k<exponent_T6}
  k= \lceil \textup{exp}(\widetilde{c}m)\rceil \leq \textup{exp}(\widetilde{c}m) +1 \leq
  \textup{exp}(2\widetilde{c}m),
  \end{equation}if $\widetilde{C}_{0}$ is chosen large enough. Therefore
  \[
  k^{7} \leq \textup{exp}(14\widetilde{c}m).
  \]Since $\widetilde{C}_{0}$ is a positive absolute constant, we see from \eqref{c_wave_view_T6} that $\widetilde{c}$ is a positive absolute constant. We have
  \begin{equation}\label{T3_low_estim_T6}
  T_{3} \geq \pi(2qx)\Bigl(\frac{y}{2q\ln x}\Bigr)^{\textup{exp}(14\widetilde{c}m)}=
  \pi(2qx)\Bigl(\frac{y}{2q\ln x}\Bigr)^{\textup{exp}(Cm)},
  \end{equation}where $C = 14\widetilde{c} >0$ is an absolute constant. From \eqref{N4_view_T6}, \eqref{N4>T3_for_T6} and \eqref{T3_low_estim_T6} we obtain
  \begin{align*}
  &\#\{qx< p_{n} \leq 2qx - 5q:\ p_{n}\equiv \ldots \equiv p_{n+m} \equiv
 a\ \text{(mod $q$)},\notag\\
 &\qquad p_{n+m} - p_{n} \leq y\} \geq \pi(2qx)\Bigl(\frac{y}{2q\ln x}\Bigr)^{\textup{exp}(Cm)}.
  \end{align*} Applying \eqref{k<exponent_T6}, we see that the inequality
  \[
  k\leq y^{\varepsilon/14}
  \]holds if
  \[
  \textup{exp}(2\widetilde{c}m) \leq y^{\varepsilon/14}.
  \]This inequality is equivalent to
   \[
  m \leq \frac{\varepsilon}{28 \widetilde{c}} \ln y =c\cdot \varepsilon \ln y,
  \]where $c= 1/(28 \widetilde{c})>0$ is an absolute constant. Let us denote $c(\varepsilon)$ by $c_{0}(\varepsilon)$. Lemma \ref{Basic_Lemma_for_T6} is proved.

  \begin{lemma}\label{Lemma_REAL_for_T6}
  There are positive absolute constants $c$ and $C$ such that the following holds. Let $\varepsilon$ be a real number with $0< \varepsilon < 1$. Then there is a number $c_{0}(\varepsilon)>0$, depending only on $\varepsilon$, such that if $x\in \mathbb{R}$, $y\in\mathbb{R}$, $m\in \mathbb{Z}$, $q\in \mathbb{Z}$, $a\in \mathbb{Z}$ are such that $c_{0}(\varepsilon) \leq y \leq \ln x$, $1 \leq m \leq c\cdot\varepsilon \ln y$, $1 \leq q \leq y^{1-\varepsilon}$ and $(a,q)=1$, then
  \begin{align*}
&\#\{qx< p_{n}\leq 2qx:\ p_{n}\equiv\ldots\equiv p_{n+m}\equiv a\ (\textup{\text{mod }}q),\ p_{n+m} - p_{n} \leq y\}\\
&\geq \pi(2qx)\left(\frac{y}{2q\ln x}\right)^{\textup{exp} (Cm)}.
\end{align*}
 \end{lemma}
 \textsc{Proof of Lemma \ref{Lemma_REAL_for_T6}.} Let $c$, $C$, $c_{0}(\varepsilon)$ be the quantities in Lemma \ref{Basic_Lemma_for_T6}. We choose a quantity $\widetilde{c}_{0}(\varepsilon)$ and an absolute constant $\widetilde{C}$ later; they will be large enough. Let $\widetilde{c}_{0}(\varepsilon) \geq c_{0}(\varepsilon)$ and $\widetilde{C} \geq C$. Let $x\in \mathbb{R}$, $y\in \mathbb{R}$, $m\in \mathbb{Z}$, $q\in \mathbb{Z}$, $a\in \mathbb{Z}$ are such that  $\widetilde{c}_{0}(\varepsilon)\leq y \leq \ln x$, $1\leq m \leq c\cdot \varepsilon\ln y$, $1\leq q \leq y^{1-\varepsilon}$, $(a,q)=1$. We put $l=\lceil x\rceil$. Then, by Lemma \ref{Basic_Lemma_for_T6}, we have
 \begin{align}
&N_{1}=\#\{ql< p_{n}\leq 2ql - 5q:\ p_{n}\equiv\ldots\equiv p_{n+m}\equiv a\ (\textup{\text{mod }}q),\notag\\
 &\qquad\qquad p_{n+m} - p_{n} \leq y\}\geq \pi(2ql)\left(\frac{y}{2q\ln l}\right)^{\textup{exp} (Cm)}=:T_{1}.\label{L65_N1>T1}
\end{align}Since $x \leq l < x+1$, we have
\begin{gather*}
ql\geq qx,\\
2ql - 5q\leq 2q(x+1)-5q=2qx - 3q< 2qx.
\end{gather*}Therefore
\begin{align}
&N_{1}\leq\#\{qx< p_{n}\leq 2qx:\ p_{n}\equiv\ldots\equiv p_{n+m}\equiv a\ (\textup{\text{mod }}q),\notag\\
 &\qquad\qquad p_{n+m} - p_{n} \leq y\}=: N_{2}.\label{L65_N1<N2}
\end{align}We have
\[
x+1 \leq x^{2},
\]if $\widetilde{c}_{0}(\varepsilon)$ is chosen large enough. Hence,
\[
\ln l \leq \ln (x+1)\leq 2\ln x.
\]Since
\[
\pi(2ql)\geq \pi (2qx),
\]we have
\[
T_{1} \geq \pi(2qx)\left(\frac{y}{4q\ln x}\right)^{\textup{exp} (Cm)}=
\pi(2qx)\left(\frac{y}{q\ln x}\right)^{\textup{exp} (Cm)} \left(\frac{1}{4}\right)^{\textup{exp} (Cm)}.
\]We have
\[
2\textup{exp}(\widetilde{C}m)\leq \textup{exp}(2\widetilde{C}m),
\]if $\widetilde{C}$ is chosen large enough. Since $\widetilde{C}\geq C$, we have
\[
\left(\frac{1}{4}\right)^{\textup{exp} (Cm)} \geq \left(\frac{1}{4}\right)^{\textup{exp} (\widetilde{C}m)}=
\left(\frac{1}{2}\right)^{2\,\textup{exp} (\widetilde{C}m)}\geq
 \left(\frac{1}{2}\right)^{\textup{exp} (2\widetilde{C}m)}.
\]Since
\[
0<\frac{y}{q\ln x}\leq 1,
\]we have
\[
\left(\frac{y}{q\ln x}\right)^{\textup{exp} (Cm)} \geq
\left(\frac{y}{q\ln x}\right)^{\textup{exp} (2\widetilde{C}m)}.
\]Hence,
\begin{equation}\label{L65_T1_EST}
T_{1} \geq \pi(2qx)\left(\frac{y}{2q\ln x}\right)^{\textup{exp} (2\widetilde{C}m)}.
\end{equation}From \eqref{L65_N1>T1}, \eqref{L65_N1<N2} and \eqref{L65_T1_EST} we obtain
\begin{align*}
&\#\{qx< p_{n}\leq 2qx:\ p_{n}\equiv\ldots\equiv p_{n+m}\equiv a\ (\textup{\text{mod }}q),\notag\\
 &\qquad p_{n+m} - p_{n} \leq y\}\geq \pi(2qx)\left(\frac{y}{2q\ln x}\right)^{\textup{exp} (2\widetilde{C}m)}.
\end{align*}Let us denote $\widetilde{c}_{0}(\varepsilon)$ by $c_{0}(\varepsilon)$ and $2\widetilde{C}$ by $C$. Lemma \ref{Lemma_REAL_for_T6} is proved.

Let us complete the proof of Theorem \ref{T6}. Let $c_{0}(\varepsilon)$, $c$, $C$ be the quantities in Lemma \ref{Lemma_REAL_for_T6}. We choose a quantity $\widetilde{c}_{0}(\varepsilon)$ and an absolute constant $\widetilde{C}$ later; they will be large enough. Let $\widetilde{c}_{0}(\varepsilon) \geq c_{0}(\varepsilon)$ and $\widetilde{C}\geq C$.

I) Let us prove the following statement.

 \emph{Let $\varepsilon$ be a real number with $0< \varepsilon <1$. Let $t\in \mathbb{R}$, $y\in\mathbb{R}$, $m\in\mathbb{Z}$, $q\in\mathbb{Z}$, $a\in\mathbb{Z}$ be such that
\begin{gather*}
t\geq 100,\quad \widetilde{c}_{0}(\varepsilon) \leq y \leq \ln \Bigl(\frac{t}{2\ln t}\Bigr),\\
1\leq m \leq c\cdot\varepsilon \ln y,\quad 1\leq q \leq y^{1-\varepsilon},\quad (a,q)=1.
\end{gather*}Then
\begin{align}
&\#\{t/2< p_{n}\leq t:\ p_{n}\equiv\ldots\equiv p_{n+m}\equiv a\ (\textup{\text{mod }}q),\ p_{n+m} - p_{n} \leq y\}\notag\\
&\geq \pi(t)\Bigl(\frac{y}{2q\ln t}\Bigr)^{\textup{exp} (\widetilde{C}m)}.\label{I_NUMBER_T6}
\end{align}}Indeed, since $t\geq 100$, we have $2\ln t \geq 1$. Hence,
\[
y\leq \ln \Bigl(\frac{t}{2\ln t}\Bigr) \leq \ln t.
\]We have
\[
q\leq y^{1-\varepsilon} \leq y \leq \ln t.
\]Therefore
\[
y\leq \ln \Bigl(\frac{t}{2\ln t}\Bigr) \leq \ln \Bigl(\frac{t}{2q}\Bigr).
\]We put
\[
x=\frac{t}{2q}.
\]We have $x\in \mathbb{R}$, $y\in\mathbb{R}$, $m\in\mathbb{Z}$, $q\in\mathbb{Z}$, $a\in\mathbb{Z}$,
\begin{gather*}
c_{0}(\varepsilon)\leq  y \leq \ln x,\\
1\leq m \leq c\cdot\varepsilon \ln y,\quad 1\leq q \leq y^{1-\varepsilon},\quad (a,q)=1.
\end{gather*}By Lemma \ref{Lemma_REAL_for_T6}, we have
\begin{align*}
&\#\{qx< p_{n}\leq 2qx:\ p_{n}\equiv\ldots\equiv p_{n+m}\equiv a\ (\textup{\text{mod }}q),\ p_{n+m} - p_{n} \leq y\}\\
&\geq \pi(2qx)\left(\frac{y}{2q\ln x}\right)^{\textup{exp} (Cm)} \geq \pi(2qx)\left(\frac{y}{2q\ln x}\right)^{\textup{exp} (\widetilde{C}m)}\geq\\
&\geq  \pi(2qx)\left(\frac{y}{2q\ln (2qx)}\right)^{\textup{exp} (\widetilde{C}m)}.
\end{align*}Returning to the variable $t$, we obtain \eqref{I_NUMBER_T6}.

II) Let us prove the following statement.

\emph{Let $\varepsilon$ be a real number with $0< \varepsilon <1$. Let $t\in \mathbb{R}$, $m\in\mathbb{Z}$, $q\in\mathbb{Z}$, $a\in\mathbb{Z}$ be such that
\begin{gather*}
t\geq 100,\quad \widetilde{c}_{0}\Bigl(\frac{\varepsilon}{2}\Bigr) \leq \ln \Bigl(\frac{t}{2\ln t}\Bigr),\\
1\leq m \leq \frac{c}{4}\cdot\varepsilon \ln\ln t,\quad 1\leq q \leq (\ln t)^{1-\varepsilon},\quad (a,q)=1.
\end{gather*}Then
\begin{align*}
&\#\Bigl\{t/2< p_{n}\leq t:\ p_{n}\equiv\ldots\equiv p_{n+m}\equiv a\ (\textup{\text{mod }}q),\\
 &\quad\ \  p_{n+m} - p_{n} \leq \ln\Bigl(\frac{t}{2\ln t}\Bigr)\Bigr\}\geq \pi(t)\Bigl(\frac{1}{4q}\Bigr)^{\textup{exp} (\widetilde{C}m)}.
\end{align*}}
\textsc{Proof.} We need the following

\begin{lemma}\label{Lemma_PREPARATORY_T6}
Let $t$ be a real number with $t \geq 100$. Then
\begin{align*}
&1)\ 2\ln t \leq \sqrt{t};\\
&2)\ \ln\Bigl(\frac{t}{2\ln t}\Bigr)\geq \frac{1}{2}\ln t;\\
&3)\ \ln\ln \Bigl(\frac{t}{2\ln t}\Bigr)\geq \frac{1}{2}\ln\ln t;\\
&4)\ 1-\frac{\ln (2\ln t)}{\ln t} \geq \frac{1}{2}.
\end{align*}
\end{lemma}
 \textsc{Proof of Lemma \ref{Lemma_PREPARATORY_T6}.}

 1) Let us consider the function
 \[
 f(x)=\sqrt{x} - 2\ln x.
 \] Then
 \[
 f^{\prime}(x)=\frac{1}{2\sqrt{x}} - \frac{2}{x}=\frac{1}{2x}(\sqrt{x}-4)>0,
 \]if $x>16$. Hence, the function $f(x)$ is increasing on the interval $(16, +\infty)$. Since
 \[
 f(100)=0.78\ldots >0,
 \]we see that $f(x)>0$ on $[100, +\infty)$. Hence, $f(t)>0$. The inequality 1) is proved.

 2) Taking into account that $t\geq 100$ and applying the inequality 1), we have
 \[
 \frac{t}{2\ln t}\geq \sqrt{t}.
 \]Therefore
 \[
 \ln \Bigl(\frac{t}{2\ln t}\Bigr)\geq \ln \sqrt{t}=\frac{1}{2}\ln t .
 \]The inequality 2) is proved.

 3) If $z \geq 4$, then
  \[
 \frac{z}{2}\geq \sqrt{z}.
 \]Since
 \[
 \ln t\geq \ln 100 = 4.605\ldots >4,
 \]we have
 \[
 \frac{1}{2}\ln t\geq \sqrt{\ln t}.
 \]Applying the inequality 2), we have
 \[
 \ln \Bigl(\frac{t}{2\ln t}\Bigr)\geq \sqrt{\ln t}.
 \]Taking logarithms, we obtain
 \[
 \ln\ln\Bigl(\frac{t}{2\ln t}\Bigr)\geq \ln\sqrt{\ln t}=\frac{1}{2}\ln\ln t.
 \]The inequality 3) is proved.

 4) Let us consider the function
 \[
 f(x)=\frac{x}{2}-\ln (2x).
 \]We have
 \[
 f^{\prime}(x)=\frac{1}{2}-\frac{1}{x}>0,
 \]if $x>2$. Hence, the function $f(x)$ is increasing on the interval $(2, +\infty)$. Since
 \[
 f(4.5)=0.05\ldots >0,
 \]we see that $f(x)>0$ on $[4.5, +\infty)$. Since
 \[
 \ln t\geq \ln 100= 4.605\ldots > 4.5,
 \]we have
 \[
 \frac{\ln t}{2} - \ln(2\ln t)\geq 0.
 \]We obtain
 \[
 \frac{\ln (2\ln t)}{\ln t}\leq \frac{1}{2}.
 \]Hence,
 \[
 1-\frac{\ln (2\ln t)}{\ln t}\geq \frac{1}{2}.
 \]The inequality 4) is proved. Lemma \ref{Lemma_PREPARATORY_T6} is proved.

 We put
\[
y=\ln\Bigl(\frac{t}{2\ln t}\Bigr).
\]
 Since $t\geq 100$, we have (see Lemma \ref{Lemma_PREPARATORY_T6}, the inequality 3))
 \[
 \ln y = \ln\ln \Bigl(\frac{t}{2\ln t}\Bigr)\geq \frac{1}{2}\ln\ln t.
 \]Therefore
 \[
 1\leq m\leq c\cdot\frac{\varepsilon}{2}\ln y.
 \]We may assume that
 \[
 \widetilde{c}_{0}(\varepsilon) \geq 2^{1/\varepsilon}.
 \]Since $t \geq 100$, we have
  \[
 \frac{t}{2\ln t} \leq t.
 \]We have
 \[
 \widetilde{c}_{0}\Bigl(\frac{\varepsilon}{2}\Bigr) \leq \ln \Bigl(\frac{t}{2\ln t}\Bigr) \leq \ln t.
 \]Hence,
 \[
 t \geq \textup{exp}\Bigl(\widetilde{c}_{0}\Bigl(\frac{\varepsilon}{2}\Bigr)\Bigr)
 \geq \textup{exp}(2^{2/\varepsilon}).
 \]Therefore
 \begin{equation}\label{II_ineq_T6}
 \frac{1}{2} (\ln t)^{1- \varepsilon/2}\geq (\ln t)^{1-\varepsilon}.
 \end{equation}From \eqref{II_ineq_T6} and the inequality 4) of Lemma \ref{Lemma_PREPARATORY_T6} we have
 \begin{align*}
 &y^{1 - \varepsilon/2} = \left(\ln \Bigl(\frac{t}{2\ln t}\Bigr)\right)^{1- \varepsilon/2}
 =(\ln t)^{1- \varepsilon/2} \left(1- \frac{\ln (2\ln t)}{\ln t}\right)^{1-\varepsilon/2}\geq\\
 &\geq (\ln t)^{1- \varepsilon/2} \left(\frac{1}{2}\right)^{1-\varepsilon/2}\geq
 \frac{1}{2}(\ln t)^{1- \varepsilon/2} \geq (\ln t)^{1-\varepsilon}.
 \end{align*}Since
 \[
 1 \leq q \leq (\ln t)^{1-\varepsilon},
 \]we have
 \[
 1 \leq q \leq y^{1 - \varepsilon/2}.
 \] Applying the statement of part I) with $\varepsilon /2$ and the inequality 2) of Lemma \ref{Lemma_PREPARATORY_T6}, we have
 \begin{align*}
&\#\Bigl\{t/2< p_{n}\leq t:\ p_{n}\equiv\ldots\equiv p_{n+m}\equiv a\ (\textup{\text{mod }}q),\\
&\quad p_{n+m} - p_{n} \leq \ln \Bigl(\frac{t}{2\ln t}\Bigr)\Bigr\}\geq \pi(t)\Bigl(\frac{\ln (t/(2\ln t))}{2q\ln t}\Bigr)^{\textup{exp} (\widetilde{C}m)} \geq\\
&\geq \pi(t)\Bigl(\frac{1}{4q}\Bigr)^{\textup{exp} (\widetilde{C}m)}.
\end{align*}The statement is proved.

III) \emph{Let $\varepsilon$ be a real number with $0< \varepsilon <1$. Let $t\in \mathbb{R}$, $y\in \mathbb{R}$, $m\in\mathbb{Z}$, $q\in\mathbb{Z}$, $a\in\mathbb{Z}$ be such that
\begin{gather*}
t\geq 100,\quad \widetilde{c}_{0}\Bigl(\frac{\varepsilon}{2}\Bigr) \leq \ln \Bigl(\frac{t}{2\ln t}\Bigr) \leq y \leq \ln t,\\
1\leq m \leq \frac{c}{4}\cdot\varepsilon \ln y,\quad 1\leq q \leq y^{1-\varepsilon},\quad (a,q)=1.
\end{gather*}Then
\begin{align*}
&\#\{t/2< p_{n}\leq t:\ p_{n}\equiv\ldots\equiv p_{n+m}\equiv a\ (\textup{\text{mod }}q),\\
 &\quad\ \  p_{n+m} - p_{n} \leq y\}\geq \pi(t)\Bigl(\frac{y}{4q\ln t}\Bigr)^{\textup{exp} (\widetilde{C}m)}.
\end{align*}}
\textsc{Proof.} Since
\[
y \leq \ln t,
\]we have
\[
1\leq m \leq \frac{c}{4}\cdot\varepsilon \ln \ln t,\quad 1\leq q \leq (\ln t)^{1-\varepsilon}.
\]Applying the statement of part II), we have
\begin{align*}
&\#\{t/2< p_{n}\leq t:\ p_{n}\equiv\ldots\equiv p_{n+m}\equiv a\ (\textup{\text{mod }}q),\  p_{n+m} - p_{n} \leq y\}\\
&\geq \#\Bigl\{t/2< p_{n}\leq t:\ p_{n}\equiv\ldots\equiv p_{n+m}\equiv a\ (\textup{\text{mod }}q),\\
  &\qquad p_{n+m} - p_{n} \leq \ln\Bigl(\frac{t}{2\ln t}\Bigr)\Bigr\}\geq \pi(t)\Bigl(\frac{1}{4q}\Bigr)^{\textup{exp} (\widetilde{C}m)}\geq\\
  &\geq \pi(t)\Bigl(\frac{y}{4q\ln t}\Bigr)^{\textup{exp} (\widetilde{C}m)}.
\end{align*}The statement is proved.

IV) For $0< \varepsilon < 1$ we define the quantity $t_{0}(\varepsilon)$ as follows:
\begin{align*}
&1)\ t_{0}(\varepsilon) \geq 100;\\
&2)\ \ln\Bigl(\frac{t}{2\ln t}\Bigr) \geq \max \left(\widetilde{c}_{0}\Bigl(\frac{\varepsilon}{2}\Bigr), \widetilde{c}_{0}(\varepsilon)\right) \text{ for all $t \geq t_{0}(\varepsilon)$}.
\end{align*}Let us prove the following statement.

\emph{Let $\varepsilon$ be a real number with $0< \varepsilon <1$. Let $t\in \mathbb{R}$, $y\in \mathbb{R}$, $m\in\mathbb{Z}$, $q\in\mathbb{Z}$, $a\in\mathbb{Z}$ be such that
\begin{gather*}
t\geq t_{0}(\varepsilon),\quad \max \left(\widetilde{c}_{0}\Bigl(\frac{\varepsilon}{2}\Bigr), \widetilde{c}_{0}(\varepsilon)\right) \leq  y \leq \ln t,\\
1\leq m \leq \frac{c}{4}\cdot\varepsilon \ln y,\quad 1\leq q \leq y^{1-\varepsilon},\quad (a,q)=1.
\end{gather*}Then
\begin{align*}
&\#\{t/2< p_{n}\leq t:\ p_{n}\equiv\ldots\equiv p_{n+m}\equiv a\ (\textup{\text{mod }}q),\\
 &\quad\ \  p_{n+m} - p_{n} \leq y\}\geq \pi(t)\Bigl(\frac{y}{4q\ln t}\Bigr)^{\textup{exp} (\widetilde{C}m)}.
\end{align*}}
\textsc{Proof.} Let us consider two cases.
\[
1)\  \ln \Bigl(\frac{t}{2\ln t}\Bigr) < y \leq \ln t.
\]Then
\begin{gather*}
t\geq 100,\quad \widetilde{c}_{0}\Bigl(\frac{\varepsilon}{2}\Bigr) \leq \ln \Bigl(\frac{t}{2\ln t}\Bigr) \leq y \leq \ln t,\\
1\leq m \leq \frac{c}{4}\cdot\varepsilon \ln y,\quad 1\leq q \leq y^{1-\varepsilon},\quad (a,q)=1.
\end{gather*}Applying the statement of part III), we have
\begin{align*}
&\#\{t/2< p_{n}\leq t:\ p_{n}\equiv\ldots\equiv p_{n+m}\equiv a\ (\textup{\text{mod }}q),\\
 &\quad\ \  p_{n+m} - p_{n} \leq y\}\geq \pi(t)\Bigl(\frac{y}{4q\ln t}\Bigr)^{\textup{exp} (\widetilde{C}m)}.
\end{align*}Now we consider the second case:
\[
2)\ y \leq \ln \Bigl(\frac{t}{2\ln t}\Bigr).
\]Then
\begin{gather*}
t\geq 100,\quad \widetilde{c}_{0}(\varepsilon) \leq y \leq \ln \Bigl(\frac{t}{2\ln t}\Bigr),\\
1\leq m \leq c\cdot\varepsilon \ln y,\quad 1\leq q \leq y^{1-\varepsilon},\quad (a,q)=1.
\end{gather*}Applying the statement of part I), we have
\begin{align*}
&\#\{t/2< p_{n}\leq t:\ p_{n}\equiv\ldots\equiv p_{n+m}\equiv a\ (\textup{\text{mod }}q),\\
&\quad p_{n+m} - p_{n} \leq y\}\geq \pi(t)\Bigl(\frac{y}{2q\ln t}\Bigr)^{\textup{exp} (\widetilde{C}m)}
\geq \pi(t)\Bigl(\frac{y}{4q\ln t}\Bigr)^{\textup{exp} (\widetilde{C}m)}.
\end{align*}The statement is proved.

V) For $0< \varepsilon <1$ we put
\[
\rho(\varepsilon) = \max \left(\widetilde{c}_{0}\Bigl(\frac{\varepsilon}{2}\Bigr), \widetilde{c}_{0}(\varepsilon)\right)  + t_{0}(\varepsilon).
\]Let us prove the following statement.

\emph{Let $\varepsilon$ be a real number with $0< \varepsilon <1$. Let $t\in \mathbb{R}$, $y\in \mathbb{R}$, $m\in\mathbb{Z}$, $q\in\mathbb{Z}$, $a\in\mathbb{Z}$ be such that
\begin{gather*}
\rho(\varepsilon) \leq  y \leq \ln t,\\
1\leq m \leq \frac{c}{4}\cdot\varepsilon \ln y,\quad 1\leq q \leq y^{1-\varepsilon},\quad (a,q)=1.
\end{gather*}Then
\begin{align*}
&\#\{t/2< p_{n}\leq t:\ p_{n}\equiv\ldots\equiv p_{n+m}\equiv a\ (\textup{\text{mod }}q),\\
 &\quad\ \  p_{n+m} - p_{n} \leq y\}\geq \pi(t)\Bigl(\frac{y}{2q\ln t}\Bigr)^{\textup{exp} (2\widetilde{C}m)}.
\end{align*}}
\textsc{Proof.} We have
\[
\max \left(\widetilde{c}_{0}\Bigl(\frac{\varepsilon}{2}\Bigr), \widetilde{c}_{0}(\varepsilon)\right)
\leq y \leq \ln t
\]and
\[
t\geq \textup{exp}(\rho(\varepsilon))\geq \rho(\varepsilon)\geq t_{0}(\varepsilon).
\]Applying the statement of part IV), we have
\begin{align}
&\#\{t/2< p_{n}\leq t:\ p_{n}\equiv\ldots\equiv p_{n+m}\equiv a\ (\textup{\text{mod }}q),\notag\\
 &\quad\ \  p_{n+m} - p_{n} \leq y\}\geq \pi(t)\Bigl(\frac{y}{4q\ln t}\Bigr)^{\textup{exp} (\widetilde{C}m)}.\label{V_1_for_T6}
\end{align}We may assume that
\[
\widetilde{C} \geq 2.
\]Therefore
\[
\textup{exp} (\widetilde{C}m)\geq \widetilde{C}m\geq \widetilde{C} \geq 2.
\]Hence,
\[
2\,\textup{exp} (\widetilde{C}m)\leq \textup{exp} (2\widetilde{C}m).
\]We have
\[
\left(\frac{1}{4}\right)^{\textup{exp} (\widetilde{C}m)}=
\left(\frac{1}{2}\right)^{2\,\textup{exp} (\widetilde{C}m)}\geq
\left(\frac{1}{2}\right)^{\textup{exp} (2\widetilde{C}m)}.
\]Since
\[
0< \frac{y}{q\ln t}\leq 1,
\]we have
\[
\Bigl(\frac{y}{q\ln t}\Bigr)^{\textup{exp} (\widetilde{C}m)} \geq
\Bigl(\frac{y}{q\ln t}\Bigr)^{\textup{exp} (2\widetilde{C}m)}.
\]We obtain
\begin{equation}\label{V_2_for_T6}
\Bigl(\frac{y}{4q\ln t}\Bigr)^{\textup{exp} (\widetilde{C}m)} \geq
\Bigl(\frac{y}{2q\ln t}\Bigr)^{\textup{exp} (2\widetilde{C}m)}.
\end{equation}From \eqref{V_1_for_T6} and \eqref{V_2_for_T6} we obtain
\begin{align*}
&\#\{t/2< p_{n}\leq t:\ p_{n}\equiv\ldots\equiv p_{n+m}\equiv a\ (\textup{\text{mod }}q),\\
 &\quad\ \  p_{n+m} - p_{n} \leq y\}\geq \pi(t)\Bigl(\frac{y}{2q\ln t}\Bigr)^{\textup{exp} (2\widetilde{C}m)}.
\end{align*}The statement is proved. Let us denote $\rho(\varepsilon)$ by $c_{0}(\varepsilon)$, $c/4$ by $c$ and $2\widetilde{C}$ by $C$. Theorem \ref{T6} is proved.

\textsc{Proof of Corollary \ref{C5}.} Let $c_{0}(\varepsilon)$, $c$ and $C$ be the quantities in Theorem \ref{T6}. We put
\[
C_{1}= \max\left(\frac{2}{c}, c_{0}\Bigl(\frac{1}{2}\Bigr), C\right).
\]Let $m$ be a positive integer. Let $x\in \mathbb{R}$, $y\in \mathbb{R}$ be such that
\[
\textup{exp}(C_{1}m) \leq y \leq \ln x.
\]Then
\begin{gather*}
y\geq \textup{exp}(C_{1}m) \geq C_{1}m\geq C_{1} \geq c_{0}\Bigl(\frac{1}{2}\Bigr),\\
y\geq \textup{exp}(C_{1}m)\geq  \textup{exp}\Bigl(\frac{2}{c}m\Bigr).
\end{gather*}From the last inequality we obtain
\[
m\leq c\cdot\frac{1}{2}\ln y.
\]We put $q=1$, $a=1$. We have
\begin{gather*}
c_{0}\Bigl(\frac{1}{2}\Bigr)\leq y \leq \ln x,\\
1\leq m \leq c\cdot\frac{1}{2}\ln y,\quad 1\leq q\leq y^{1/2},\quad (a,q)=1.
\end{gather*}Applying Theorem \ref{T6} with $\varepsilon = 1/2$, we have
\begin{align*}
&\#\{x/2< p_{n}\leq x:\ p_{n+m} - p_{n} \leq y\}=\\
&=\#\{x/2< p_{n}\leq x:\ p_{n}\equiv\ldots\equiv p_{n+m}\equiv a\ (\textup{\text{mod }}q),\\
 &\quad\ \  p_{n+m} - p_{n} \leq y\}\geq \pi(x)\Bigl(\frac{y}{2q\ln x}\Bigr)^{\textup{exp} (Cm)}=
 \pi(x)\Bigl(\frac{y}{2\ln x}\Bigr)^{\textup{exp} (Cm)}\geq\\
 &\geq \pi(x)\Bigl(\frac{y}{2\ln x}\Bigr)^{\textup{exp} (C_{1}m)}.
\end{align*}Let us denote $C_{1}$ by $C$. Corollary \ref{C5} is proved.

\section{Acknowledgements}\label{S_Acknowledgements}

The author is deeply grateful to Sergei Konyagin and Maxim Korolev for their attention to this work and useful comments.

The author also expresses gratitude to Mikhail Gabdullin and Pavel Grigor'ev for useful comments and suggestions.

\end{document}